\documentclass[12pt,twoside]{amsart}\usepackage{amssymb}
\usepackage{enumerate}

\newtheorem{thmspec}{\relax}

\newtheorem{theorem}{Theorem}[section]
\newtheorem{thm}[theorem]{Theorem}
\newtheorem{lem}[theorem]{Lemma}

\newtheorem{prop}[theorem]{Proposition}
\newtheorem{defi}[theorem]{Definition}

\theoremstyle{definition}

\theoremstyle{remark}

\numberwithin{equation}{section}
\def \Bbb{\mathbb}

\def\onto{{\kern3pt\to\kern-8pt\to\kern3pt}}

\def\<{\langle}
\def\>{\rangle}
\def\|{{\ |\ }}

\def\onto{\twoheadrightarrow}

\def\-{\underline}

\def\dist{\operatorname{dist}}

\def\mes{\operatorname{mes}}

\def\arg{\operatorname{arg}}

\def\End{\operatorname{End}}
\def\Alim{\mathcal{A}-\lim}
\def\Alimsup{\mathcal{A}-\limsup}
\def\N{\Bbb N}

\def\R{\Bbb R}
\def\C{\Bbb C}
\def\P{\Bbb P}
\def\B{\Bbb B}

\def\X{\Bbb X}

\def\<{\langle}
\def\>{\rangle}

\catcode`\@=11
\def\serieslogo@{\relax}
\def\@setcopyright{\relax}
\catcode`\@=12

\title[A unified approach]
{A unified approach to the theory of   separately holomorphic mappings}

\begin{document}

\author{Vi{\^e}t-Anh  Nguy\^en}
\address{Vi{\^e}t-Anh  Nguy\^en\\
Mathematics Section\\
The Abdus Salam international centre 
 for theoretical physics\\
Strada costiera, 11\\
34014 Trieste, Italy}
\email{vnguyen0@ictp.trieste.it}

\subjclass[2000]{Primary 32D15, 32D10}
\date{}

\keywords{Hartogs' theorem,  holomorphic extension,   Poletsky theory of discs,
Rosay Theorem on holomorphic discs.}

\begin{abstract}
We extend the theory of  separately holomorphic mappings between complex analytic spaces.
Our method is based on   Poletsky theory of discs,
Rosay Theorem on holomorphic discs and our recent joint-work with Pflug
on boundary cross theorems in dimension $1.$ It also  relies on our new technique of conformal mappings and 
 a generalization of Siciak's  relative extremal function.
 Our  approach illustrates the unified character:
``From local informations to global extensions". Moreover,
 it avoids systematically the use of the classical method of doubly orthogonal bases of Bergman type.
 \end{abstract}
\maketitle

\section{Introduction}

 In this article all complex manifolds are supposed to be  of  finite  dimension and  countable at infinity, 
 and all
  complex analytic spaces  are supposed to be reduced,
  irreducible, of finite dimension  and  countable at infinity.
For a subset $S$ of a topological space $M,$  $\overline{S}$   denotes the closure of $S$
in $M,$  and the set  $\partial S:=\overline{S}\cap\overline{M\setminus S}$  denotes, as usual, the boundary of $S$
in $M.$ 

\smallskip

The main purpose of this work is to investigate the following 

\medskip

\noindent{\bf  PROBLEM.} {\it
 Let $X,$ $Y$ be two complex manifolds, let  $D$ (resp. $G$) be an open subset of $X$  (resp.  $Y$), let
  $A$ (resp. $B$) be a  subset of  $\overline{D}$ (resp.
  $\overline{G}$)  and let  $Z$ be
  a complex analytic  space.
 Define the cross
\begin{equation*}
 W:=\big(A\times (G\cup B)\big) \bigcup \big((D\cup A)\times B\big).
 \end{equation*}
   We  want to determine  the ``envelope of holomorphy" of the cross $W, $ that is,
   an ``optimal" open subset  of $ X\times Y,$ denoted by  $\widehat{\widetilde{W}},$ which is characterized by the following properties:

  Let $f:\    W\longrightarrow Z$ be a mapping that satisfies, in essence,
   the following two conditions: 
   \begin{itemize}
   \item[$\bullet$] $f(a,\cdot)$ is holomorphic on $G$ for all $a\in A,$
   $f(\cdot,b)$ is holomorphic on $D$  for all $b\in B;$
   \item[$\bullet$] $f(a,\cdot)$ is continuous on $G\cup B$ for all $a\in A,$ 
   $f(\cdot,b)$ is continuous on $D\cup A$  for all $b\in B.$
    \end{itemize}
    Then there is  a   holomorphic  mapping
  $\hat{f}$ defined on   $\widehat{\widetilde{W}}$  such that for every $(\zeta,\eta)\in W,$
  $\hat{f}(z,w)$ tends to $f(\zeta,\eta)$ as $(z,w)\in\widehat{\widetilde{W}}$ tends, in some sense, to $(\zeta,\eta).$ 
  }
  
  \medskip
   
Now we  recall    briefly  the main   developments around this problem.
 All the results  obtained  so far  may be divided into two directions.
The first direction  investigates  the results in the  ``interior" context:
 $A\subset D$ and  $B\subset G,$  while 
the second one explores the ``boundary"  context:  $A\subset\partial  D$  and  $B\subset\partial G.$

 \smallskip

The first  fundamental result in the field of separate holomorphy is 
   the well-known Hartogs extension theorem for separately holomorphic functions
(see \cite{ha}). In the language  of the PROBLEM  
 the following case   $X=\C^n,\  Y=\C^m,\  A=D,\  B=G,  \ Z=\C$  has been solved  and the
result is  $\widehat{\widetilde{W}}=D\times G.$
In particular, this theorem may be considered as the first main result in the first direction.
 In his   famous article \cite{ber}   Bernstein obtained some  positive results for the  PROBLEM 
in certain  cases where
$A\subset D,$  $B\subset G,$ $X= Y=\C$ and  $Z=\C.$

\smallskip

More than  60 years later,  a next  important impetus  was  made by Siciak  (see \cite{si1,si2}) in  1969--1970,
 where he established some significant  generalizations  of
the Hartogs
extension theorem. In fact,    Siciak's   formulation of these generalizations gives rise to
the above  PROBLEM:    to determine the envelope of holomorphy for  separately
holomorphic functions defined on some {\it cross sets}  $W.$ The theorems
obtained under this formulation are often called {\it cross theorems.}
 Using  the so-called {\it relative
extremal function,}  Siciak completed the PROBLEM for  the case where
$A\subset D,$  $B\subset G,$ $X= Y=\C$ and  $Z=\C.$

\smallskip

The next deep steps were initiated by Zahariuta in 1976 (see \cite{za}) when he started to use
the method of common bases of Hilbert spaces. This original approach permitted
him to obtain new  cross theorems for some cases where $A\subset D,\  B\subset G$ and $D=X,$  $G=Y$ are  Stein manifolds. 
As a consequence, he was able to generalize
the result of Siciak in higher dimensions.

\smallskip

Later, Nguy\^en  Thanh V\^an and Zeriahi (see \cite{nz1,nz2,nz3}) developed the method of doubly orthogonal bases of
Bergman type in order to generalize the result of
Zahariuta. This  is a significantly simpler and    more constructive version
of Zahariuta's original  method.  Nguy\^en  Thanh V\^an and Zeriahi have recently achieved an elegant
 improvement  of their method  (see \cite{ng2}, \cite{zer}).

\smallskip

Using  Siciak's method, Shiffman  (see \cite{sh1}) was the first to generalize some
Siciak's results  to separately holomorphic mappings with values in a complex analytic space  $Z.$
Shiffman's result (see \cite{sh2}) shows that the natural ``target spaces" for obtaining
satisfactory generalizations of cross theorems are the ones which possess {\it the Hartogs extension property}
(see Subsection  \ref{HEP} below for more explanations).



\smallskip

In 2001 Alehyane and Zeriahi   solved  the PROBLEM for the case where $A\subset D,$  $B\subset G$ and $X,\ Y$  are Stein manifolds,
and $Z$ is a complex analytic space which possesses the Hartogs extension property
(see Theorem 2.2.4 in  \cite{az}).

\smallskip

In a recent  work (see \cite{nv}) we complete, in some sense,  the PROBLEM for the case  where  $A\subset D,$  $B\subset G$ and $X,\ Y$  are   arbitrary
complex manifolds. The main ingredients in our approach are  Poletsky theory of   discs developed in \cite{po1, po2},
    Rosay's Theorem on holomorphic discs (see \cite{ro}), the above mentioned result of  Alehyane--Zeriahi 
    and the technique of level sets of
the plurisubharmonic measure which was previously  introduced in our 
joint-work with Pflug (see \cite{pn1}). 

\smallskip

To conclude   the first direction  of research we mention the  survey
articles by  Nguy\^en Thanh V\^an \cite{ng} and Peter  Pflug \cite{pf}  which  give  nice accounts on this subject.

\smallskip

The  first result in the second direction (i.e.  ``boundary context")  was  established in  the  work of  Malgrange--Zerner \cite{zen} in the  1960s. 
Further  results in this direction were obtained by
Komatsu \cite{ko} and Dru\.{z}kowski \cite{dr}, but only for some special
cases.
Recently, Gonchar \cite{go1,go2} has proved  a  more
general result where the following case  has    been solved:  $X=Y=\C,$ $D$ and $G$ are Jordan domains,
  $ A$ (resp. $B$) is an open boundary  subset of $\partial D$ (resp. $\partial G$), and   $ Z=\C.$
   It should be noted that  
  Airapetyan and Henkin published  a general version of
 the edge-of-the-wedge theorem for CR manifolds (see  \cite{ah1}  for a brief version
 and  \cite{ah2} for a complete proof).
   Gonchar's result could be deduced from the latter works.
  In  our joint-articles with Pflug (see \cite{pn1,pn2,pn3}),  Gonchar's result has been generalized
  considerably. More precisely, 
    the work in \cite{pn3} treats the case where   the ``source spaces" $X,Y$ are arbitrary  complex manifolds,
 $ A$ (resp. $B$) is an open boundary  subset of $\partial D$ (resp. $\partial G$), and   $ Z=\C.$  
 The work in \cite{pn2}  solves the case where   the ``source spaces" $X,Y$ are Riemann surfaces,
 $ A$ (resp. $B$) is a measurable (boundary)  subset of $\partial D$ (resp. $\partial G$), and   $ Z=\C.$   

\smallskip

The main purpose of this article  is to give a new version  of the Hartogs extension theorem which unifies
all results   up to now. Namely, we are able to  give a reasonable solution to the   PROBLEM
when  the  ``target space" $Z$  possesses  the Hartogs extension property. Our method is based on 
 a systematic application of Poletsky theory of discs,
Rosay Theorem on holomorphic discs and our  joint-work with Pflug
on boundary  cross theorems in dimension $1$  (see \cite{pn2}). It also relies on  our new technique of conformal mappings and  a generalization of
Siciak's  relative extremal function. The  approach illustrates the unified character
in the theory of extension of holomorphic mappings:

\smallskip

{\bf One can deduce the  global extension  from local informations.}

\smallskip

 Moreover,
 the novelty of this new approach is 
 that one does not use  the classical method of doubly orthogonal bases of Bergman type.

\smallskip

  We close  the introduction with a brief outline of the paper to follow.

\smallskip

 In Section 2 we formulate  the main results.

\smallskip

The tools which are needed for the proof of the main results are developed
in Section 3, 4, 5  and 7.

\smallskip

   The proof of the    main results    is divided into three
   parts, which correspond to Section 6, 8 and 9.
 Section 10 concludes the article with  various applications of our results.

 \smallskip

\smallskip

\indent{\it{\bf Acknowledgment.}}
 The paper was written while  the  author was visiting  the  Abdus Salam International Centre 
 for Theoretical Physics
in Trieste. He wishes to express his gratitude to this organization.

\section{Preliminaries and statement of the main result}
First  we  develop  some  new notions  such  as  system of  approach  regions for an open set
in a  complex  manifold,
and  the corresponding plurisubharmonic measure. These  will provide the  framework
for an  exact  formulation of the PROBLEM and for  our solution. 
 \subsection{Approach regions, local pluripolarity and plurisubharmonic measure}\label{Subsection_approach_regions}
\begin{defi}\label{defi_approach_region}
  Let $X$ be a complex manifold  and let $D\subset X$ be an open subset.     
  A  {\rm system of approach regions} for $D$ is a collection
  $\mathcal{A}=\big(\mathcal{A}_{\alpha}(\zeta)\big)_{\zeta\in\overline{ D},\  \alpha\in I_{\zeta}}$
  of open subsets of $D$
  with the following properties:
  \begin{itemize}
  \item[(i)] For all $\zeta\in D,$ the system $\big(\mathcal{A}_{\alpha}(\zeta)\big)_{ \alpha\in I_{\zeta}}$
  forms a basis of open neighborhoods of $\zeta$  (i.e., for any  open neighborhood $U$  of a point  $\zeta\in D,$
  there is $ \alpha\in I_{\zeta}$ such that  $\zeta\in \mathcal{A}_{\alpha}(\zeta)\subset U$).
  \item[(ii)] For  all $\zeta\in\partial D$ and
     $\alpha\in I_{\zeta},$   
  $\zeta\in \overline{\mathcal{A}_{\alpha}(\zeta)}.$
  \end{itemize}
 $\mathcal{A}_{\alpha}(\zeta)$ is often
  called an {\rm approach region} at $\zeta.$

 $\mathcal{A}$  is
said to be {\rm canonical} if it satisfies (i) and the following  property (which is stronger than (ii)):
  \begin{itemize}
  \item[(ii')]
For every point  $\zeta\in \partial D,$
  there is  a basis of open neighborhoods  $(U_{\alpha})_{\alpha\in I_{\zeta}}$ of $\zeta$  in $X$ such that 
   $ \mathcal{A}_{\alpha}(\zeta)=U_{\alpha}\cap D,$  $\alpha\in I_{\zeta}.$
   \end{itemize}
  \end{defi}
   It is possible that  $I_{\zeta}=\varnothing$ for some $\zeta\in\partial D.$


Various  systems of approach regions which one often encounters in Complex Analysis
 will be described in the  next  subsection.   Systems of approach regions  for $D$ are used to deal with the limit  at points in $\overline{D}$
 of  mappings defined on  some  open subsets of $D.$ Consequently, we deduce from 
 Definition  \ref{defi_approach_region} that 
 the subfamily $\big(\mathcal{A}_{\alpha}(\zeta)\big)_{\zeta\in D,\ \alpha\in I_{\zeta}}$
 is, in  a certain sense, independent of the choice of a system of approach regions  $\mathcal{A}.$
 In addition, any two canonical systems of approach regions are, in some sense, equivalent.
 These observations lead us to use, throughout the paper,  the following {\bf convention}:

{\it We fix, for every open set $D\subset X,$  a {\bf canonical system of approach regions.}
 When we want to define a  system of approach regions $\mathcal{A}$ for an open set $D\subset X,$ 
 we only need to specify the subfamily   $\big(\mathcal{A}_{\alpha}(\zeta)\big)_{\zeta\in\partial 
  D,\ \alpha\in I_{\zeta}}.$}
 
  In what follows we fix an open subset $D\subset X$ and a  system of approach regions   $\mathcal{A}=\big(\mathcal{A}_{\alpha}(\zeta)\big)_{
  \zeta\in\overline{ D},\ 
   \alpha\in I_{\zeta}}$ for $D.$ 

 For every function $u:\ D\longrightarrow [-\infty,\infty),$ let
\begin{equation*}
 (\Alimsup u)(z):=
\begin{cases}
 \sup\limits_{\alpha\in I_{z}}\limsup\limits_{w\in \mathcal{A}_{\alpha}(z),\ w\to z}u(w), & z\in\overline{D},\  I_z\not=\varnothing,\\
 \limsup\limits_{w\in D,\  w\to  z} u(w), & z \in \partial D,\ I_z=\varnothing.
\end{cases}
\end{equation*}
By  Definition \ref{defi_approach_region} (i),
 $(\Alimsup u)|_D$ coincides with the usual 
{\it upper semicontinuous regularization}  of $u.$
  
For a set  $A\subset \overline{D}$ put
\begin{equation*}
h_{A,D}:=\sup\left\lbrace u\ :\  u\in\mathcal{PSH}(D),\ u\leq 1\ \text{on}\ D,\
   \Alimsup u\leq 0\ \text{on}\ A    \right\rbrace,
\end{equation*}
where $\mathcal{PSH}(D)$ denotes the cone  of all functions  plurisubharmonic
on $D.$

$A$ is said to be {\it pluripolar} in $D$ if there is  $u\in
\mathcal{PSH}(D)$ such that $u$ is not identically $-\infty$ on every
connected component of $D$ and $A\subset \left\lbrace z\in D:\ u(z)=-\infty\right\rbrace.$
$A$ is said to be {\it locally  pluripolar} in $D$ if  for any
$z\in A,$ there is an open  neighborhood $V\subset D$ of $z$ such that $A\cap V$ is
pluripolar in $V.$ $A$ is said  to be {\it nonpluripolar} (resp. {\it non locally  pluripolar}) if it is not
pluripolar (resp. not locally pluripolar).
According to a classical result of Josefson and Bedford  (see \cite{jo},
\cite{be}), if $D$ is a Riemann domain over a Stein manifold,
then   $A\subset D$ is   locally  pluripolar if and only if it is
pluripolar.


\begin{defi}\label{defi_relative_extremal}
 The
  {\rm relative extremal function of $A$  relative to $D$} is
 the function $ \omega(\cdot,A,D)$
 defined by  
 \begin{equation*}
 \omega(z,A,D)= \omega_{\mathcal{A}}(z,A,D):= (\Alimsup  h_{A,D})(z),\qquad  z\in \overline{D}.\quad \footnotemark
\end{equation*}
\end{defi}
\footnotetext{Observe that this function depends on the system of
approach regions.}
 Note that when $A\subset D,$  Definition \ref{defi_relative_extremal}  
   coincides with the classical
definition of Siciak's  relative extremal function.

 Next, we say that a  set  $A\subset \overline{D}$ is {\it
locally pluriregular at a point $a\in \overline{A}$}   if $\omega(a,A\cap U,D\cap U)=0$
  for  all open neighborhoods $U$ of $a.$
Moreover, $ A$ is said to be {\it locally pluriregular } if it is locally
pluriregular at all points $a\in A.$ 
It should be noted from  Definition  \ref{defi_approach_region} that  if $a\in \overline{A}\cap D$ then  the property of local pluriregularity
of $A$ at $a$  does not depend  on any particular choices of  a system of approach regions $\mathcal{A},$
while the situation is different when  $a\in \overline{A}\cap \partial D:$ the  property    does depend on  $\mathcal{A}.$

 We denote by $A^{\ast}$ the following set
 \begin{equation*}
  (A\cap\partial D)\bigcup\left\lbrace   a\in \overline{A}\cap D:\ A\ \text{is locally pluriregular at}\ a \right\rbrace.
  \end{equation*}
If $A\subset D$ is non locally pluripolar, then  a classical result of Bedford and Taylor
(see \cite{be,bt})
says that    $A^{\ast}$ is locally
pluriregular    
 and  $A\setminus A^{\ast}$ is locally  pluripolar.
Moreover, $A^{\ast}$ is locally  of type $\mathcal{G}_{\delta},$  that  is,
for every $a\in A^{\ast}$ there is an open  neighborhood $U\subset D$ of $a$
such that $A^{\ast}\cap U$ is  a countable intersection of open sets.

Now we are  in the position to formulate the following  version
of the  plurisubharmonic measure.
\begin{defi}\label{defi_pluri_measure}
For a set $A\subset \overline{D},$ let 
$\widetilde{A} =\widetilde{A}(\mathcal{A}):=\bigcup\limits_{P\in \mathcal{E}(A)} P,$ where
\begin{equation*}
\mathcal{E}(A)=\mathcal{E}(A,\mathcal{A}):=\left\lbrace P\subset \overline{D}:\  P\ \text{is  locally pluriregular,}\
\overline{P}\subset  A^{\ast}      \right\rbrace,
\end{equation*}
 The {\rm  plurisubharmonic measure of $A$  relative to $D$} is
 the function $\widetilde{\omega}(\cdot,A,D)$
 defined by
\begin{equation*}
\widetilde{\omega}(z,A,D):=  \omega(z,\widetilde{A},D),\qquad  z\in D.
\end{equation*}
\end{defi}

It is worthy to remark that  $\widetilde{\omega}(\cdot,A,D)\in\mathcal{PSH}(D)$ and
$0\leq \widetilde{\omega}(z,A,D)\leq 1,\ z\in  D.$  Moreover,
\begin{equation}\label{eq_defi_pluri_measure}
\Big(\Alimsup \widetilde{\omega}(\cdot,A,D)\Big)(z)=0,\qquad z\in \widetilde{A}.
\end{equation}
 An example  in  \cite{ah} shows  that     in general, 
  $\omega(\cdot,A,D)\not=\widetilde{\omega}(\cdot,A,D)$ on $D.$
Section  \ref{section_Applications}  below is devoted  to the study of  $ \widetilde{\omega}(\cdot,A,D)$
in some important cases.

Now  we compare the plurisubharmonic measure $ \widetilde{\omega}(\cdot,A,D)$ with 
Siciak's relative  extremal function $ \omega(\cdot,A,D).$   We only  consider two 
important special cases:  $A\subset D$ and  $A\subset\partial D.$
For the moment, we only focus  on the case  where  $A\subset D.$ 
The latter one  will be discussed  in 
Section   \ref{section_Applications}  below.

  If $A$ is an open subset of an arbitrary complex manifold  $D$, then it is  easy to see that
\begin{equation*}
\widetilde{\omega}(z,A,D)=  \omega(z,A,D),\qquad  z\in D.
\end{equation*}
If $A$  is  a (not necessarily open) subset of an arbitrary complex manifold $D,$
then we will prove  in Proposition \ref{prop_coincidence_two_definitions} below 
that
  \begin{equation*}
\widetilde{\omega}(z,A,D)=  \omega(z,A^{\ast},D),\qquad  z\in D.
\end{equation*}
On the other hand,  if, morever,  $D$  is  a  bounded  open subset of  $\C^n$   then we  have (see, for  example,
  Lemma 3.5.3 in \cite{jp1}) $
\omega(z,A,D)=  \omega(z,A^{\ast},D),$ $ z\in D.$
Consequently, under  the last assumption,       
\begin{equation*}
\widetilde{\omega}(z,A,D)=  \omega(z,A,D),\qquad  z\in D.
\end{equation*}
 Our discussion  shows that at least  in the case  where  $A\subset D$, the notion of  the plurisubharmonic measure
is  a  good candidate  for  generalizing  Siciak's relative  extremal function  to the manifold context in the theory
of separate holomorphy.

 For a good background of the pluripotential
theory, see the books  \cite{jp1} or  \cite{kl}.
%
%


%
%
%
%
%
%
%
\subsection{Examples of systems of  approach  regions} \label{Subsection_Examples}
  There  are  many  systems    of  approach  regions which are very useful 
in Complex  Analysis. In this  subsection  we present some of them.

\noindent {\bf 1. Canonical  system   of  approach regions.}
It has been given by   Definition \ref{defi_approach_region} (i)--(ii').

\noindent {\bf 2. System of  angular (or  Stolz)    approach regions for the open unit  disc.}
Let $E$  be the  open unit  disc of  $\C.$ Put
\begin{equation*}
\mathcal{A}_{\alpha}(\zeta):=
 \left\lbrace          t\in E:\ \left\vert
 \arg\left(\frac{\zeta-t}{\zeta}\right)
 \right\vert<\alpha\right\rbrace   ,\qquad  \zeta\in\partial E,\  0<\alpha<\frac{\pi}{2},
\end{equation*}
where  $\arg:\ \C\longrightarrow (-\pi,\pi]$ is  as usual the argument function. 
$ \mathcal{A}=\left(\mathcal{A}_{\alpha}(\zeta)\right)_{\zeta\in\partial E,\ 
 0<\alpha<\frac{\pi}{2}}$ is   referred to as   {\it   the system of  angular (or  Stolz) 
    approach regions for $E.$}
In this  context  $\Alim$ is also called {\it angular limit}.

\noindent {\bf 3. System of   angular    approach regions for   certain ``good" open subsets of Riemann surfaces.}
Now we  generalize the previous construction (for the open unit disc) to  a  global  situation. More  precisely, we
will use as  the  local  model the  system    of   angular    approach regions for  $E.$
Let $X$ be  a  complex  manifold of  dimension $1,$ in other words, $X$ is  a Riemann surface,
and   $D\subset X$  an  open  set. Then $D$   is   said  to be  {\it
good at a  point $\zeta\in\partial D$}\footnote{ In the  work \cite{pn2}
 we  use  the more  appealing  word  {\it Jordan-curve-like} for this  notion.}
 if there is a  Jordan domain 
   $U\subset X$  such that   $\zeta\in U$ and $U\cap
\partial D$ is the  interior of  a  Jordan curve.  

Suppose that   $D$   is  good  at $\zeta.$ This point is   said  to be 
 {\it of type 1} if there  is a neighborhood $V$ of
$\zeta$ such that $V_0=V\cap D$ is  a Jordan domain.  Otherwise,
$\zeta$ is   said  to be {\it   of type 2}. We  see  easily  that if 
$\zeta$ is of  type 2, then  there  are an open  neighborhood    $V$ of $\zeta$ and  two disjoint  Jordan domains  $V_1,$ $V_2$ 
 such that  $V\cap
D=V_1\cup V_2.$  Moreover,    $D$  is   said  to be   {\it
good  on  a  subset  $A$}  of  $\partial D$ if  $D$ is
good at all points  of   $A.$

Here  is  a simple  example   which may clarify the  above definitions.
Let $G$ be  the  open  square in $\C$
with vertices $1+i,$ $-1+i,$ $-1-i,$ and  $1-i.$
Define the domain
\begin{equation*}
D:=G\setminus \left [-\frac{1}{2},\frac{1}{2}\right].
\end{equation*}
Then   $D$ is  good on $\partial G\cup  \left
(-\frac{1}{2},\frac{1}{2}\right).$ All points  of $\partial G$ are of  type
1 and all  points of $\left(-\frac{1}{2},\frac{1}{2}\right)$ are of type
2.

Suppose  now  that   $D$  is  good  
on a  nonempty subset  $A$ of $\partial D.$ 
 We define {\it the system    of angular approach regions   supported on  $A$:}
   $\mathcal{A}=\big(\mathcal{A}_{\alpha}(\zeta)\big)_{\zeta\in\overline{ D},\  \alpha\in I_{\zeta}}$
as  follows:
  \begin{itemize}
  \item[$\bullet$] If $\zeta\in  \overline{ D}\setminus A,$  then $\big(\mathcal{A}_{\alpha}(\zeta)\big)_{ \alpha\in I_{\zeta}}$
coincide with the canonical approach  regions. 
 \item[$\bullet$] If $\zeta\in  A,$ then  by using  a conformal  mapping $\Phi$ from $V_0$  (resp. $V_1$ and $V_2$)
 onto $E$
 when $\zeta$ is  of  type 1  (resp. 2), we can ``transfer'' the  angular  approach regions
at the  point $\Phi(\zeta)\in \partial E:$ $\left(\mathcal{A}_{\alpha}(\Phi(\zeta))\right)_{ 
 0<\alpha<\frac{\pi}{2}}$ to those at the point $\zeta\in  \partial D$  (see \cite{pn2} for more  detailed explanations).  
\end{itemize}
Making use  of conformal  mappings in  a  local  way, we  can   transfer, in  the  same  way, many
notions which  exist   on  $E$  (resp.  $\partial E$)  to those  on  $D$  (resp.  $\partial D$).  

\noindent {\bf 4. System     of conical approach regions.}

Let $D\subset\C^n$ be a  domain  and $A\subset\partial D.$ Suppose in  addition that
for every  point $\zeta\in A$ there exists the (real) tangent space $T_{\zeta}$ to $\partial D$ at $\zeta.$   
 We  define {\it the system   of   conical approach regions supported on $A$:}
   $\mathcal{A}=\big(\mathcal{A}_{\alpha}(\zeta)\big)_{\zeta\in\overline{ D},\  \alpha\in I_{\zeta}}$
as  follows:
  \begin{itemize}
  \item[$\bullet$] If $\zeta\in  \overline{ D}\setminus A,$  then $\big(\mathcal{A}_{\alpha}(\zeta)\big)_{ \alpha\in I_{\zeta}}$
coincide with the  canonical    approach regions.
 \item[$\bullet$] If $\zeta\in  A,$ then
 \begin{equation*}
 \mathcal{A}_{\alpha}(\zeta):=\left\lbrace  z\in D:\  \vert  z-\zeta\vert <\alpha\cdot \dist(z,T_{\zeta})  \right\rbrace,
 \end{equation*}
where $I_{\zeta}:=(1,\infty)$  and $ \dist(z,T_{\zeta})$ denotes the Euclidean    distance  from the  point $z$ to $T_{\zeta}.$
\end{itemize}

We can also   generalize the previous construction to  a  global situation:

{\it  $X$ is an  arbitrary complex manifold, $D\subset  X$ is an open set and  $A\subset \partial D$ is  a  subset
with  the property that
at every point $\zeta\in A$   
 there exists the  (real)  tangent  space  $T_{\zeta}$ to $\partial D.$} 

We can  also    formulate  the  notion of points of  type 1  or  2 in this  general  context  in the same way
as we have  already done
in Paragraph 3  above.

\subsection{Cross and  separate holomorphicity and $\mathcal{A}$-limit.} \label{subsection_cross}
Let $X,\ Y$  be two complex manifolds,
  let $D\subset X,$ $ G\subset Y$ be two nonempty open sets, let
  $A\subset \overline{D}$   and  $B\subset \overline{G}.$
  Moreover, $D$  (resp.  $G$) is equipped with a 
  system of approach regions 
  $\mathcal{A}(D)=\big(\mathcal{A}_{\alpha}(\zeta)\big)_{\zeta\in\overline{D},\  \alpha\in I_{\zeta}}$
  (resp.  $\mathcal{A}(G)=\big(\mathcal{A}_{\alpha}(\eta)\big)_{\eta\in\overline{G},\  \alpha\in I_{\eta}}$).
 We define
a {\it $2$-fold cross} $W,$  its {\it  interior} $W^{\text{o}}$ and   its {\it  regular part} $\widetilde{W}$
(with  respect to $\mathcal{A}(D)$ and $\mathcal{A}(G)$)
as
\begin{eqnarray*}
W &=&\X(A,B; D,G)
:=\big((D\cup A)\times B\big)\bigcup\big (A\times(B\cup G)\big),\\
W^{\text{o}} &=&\X^{\text{o}}(A,B; D,G)
:= (A\times  G)\cup (D\times B),\\
\widetilde{W} &=&\widetilde{\X}(A,B;D,G)
:= \big((D\cup \widetilde{A})\times \widetilde{B}\big)\bigcup \big(\widetilde{A}\times(G\cup \widetilde{B})\big).
\end{eqnarray*}
Moreover, put
\begin{eqnarray*}
\omega(z,w)&:=&\omega(z,A,D)+\omega(w,B,G),\qquad
(z,w)\in D\times G,\\
\widetilde{\omega}(z,w)&:=&\widetilde{\omega}(z,A,D)+\widetilde{\omega}(w,B,G),\qquad
(z,w)\in D\times G.
\end{eqnarray*}

For a $2$-fold cross $W :=\X(A,B; D,G)$
let 
\begin{eqnarray*}
\widehat{W}&:=&\widehat{\X}(A,B;D,G)
=\left\lbrace (z,w)\in D\times G:\ \omega(z,w)  <1
\right\rbrace,\\
 \widehat{\widetilde{W}} &:=&\widehat{\X}(\widetilde{A},\widetilde{B};D,G)
 =\left\lbrace (z,w)\in D\times G :\  \widetilde{\omega}(z,w)<1
\right\rbrace.
\end{eqnarray*}

Let $Z$ be a complex analytic space.
We say that a mapping $f:W^{\text{o}}\longrightarrow Z$ is {\it separately holomorphic}
  and write $f\in\mathcal{O}_s(W^{\text{o}},Z),$   if,  
 for any $a\in A $ (resp.  $b\in B$)
 the restricted mapping $f(a,\cdot)$  (resp.  $f(\cdot,b)$)  is holomorphic  on $G$  (resp. on $D$).

 We say that a mapping $f:\  W\longrightarrow Z$
   is  {\it separately continuous}
and write 
 $f\in \mathcal{C}_s\Big( W,Z  \Big)$ 
 if,  
 for any $a\in A$ (resp.  $b\in B$)
 the restricted mapping $f(a,\cdot)$  (resp.  $f(\cdot,b)$)  is continuous   on $G\cup B$ 
  (resp. on $D\cup A$).

  
 In  virtue  of (\ref{eq_defi_pluri_measure}), for every $(\zeta,\eta)\in  \widetilde{W}$ and 
 every $ \alpha\in I_{\zeta},$  $ \beta\in I_{\eta},$
  there are open neighborhoods $U$ of $\zeta$ and $V$ of $\eta$ such that
  \begin{equation*}
  \Big(U\cap \mathcal{A}_{\alpha}(\zeta)\Big)   \times \Big(V\cap \mathcal{A}_{\beta}(\eta)
  \Big)\subset \widehat{\widetilde{W}}.
  \end{equation*}
   Then a mapping 
  $f:\ \widehat{\widetilde{W}}\longrightarrow Z$ is said to {\it admit $\mathcal{A}$-limit $\lambda$ at $
(\zeta,\eta)\in \widetilde{W},$}  and one writes
\begin{equation*}
(\Alim f)(\zeta,\eta)=\lambda,\qquad\footnote{Note that here $\mathcal{A}=\mathcal{A}(D)\times\mathcal{A}(G).$}
\end{equation*}
 if, for all $\alpha\in I_{\zeta},\  \beta\in I_{\eta},$
\begin{equation*}
\lim\limits_{ \widehat{\widetilde{W}}\ni 
(z,w)\to (\zeta,\eta),\  z\in \mathcal{A}_{\alpha}(\zeta),\  w\in  \mathcal{A}_{\beta}(\eta)}f(z,w)=\lambda.
\end{equation*}

Throughout the paper,
 for a topological space  $\mathcal{M},$ $\mathcal{C}(\mathcal{M},Z)$ denotes the set of  all continuous mappings  
 $f:\ \mathcal{M}\longrightarrow  Z.$
If, moreover, $Z=\C,$ then  $\mathcal{C}(\mathcal{M},\C)$  is equipped with the ``sup-norm''
   $\vert f\vert_{\mathcal{M}}:=\sup_{\mathcal{M}} \vert f\vert\in [0,\infty].$
A mapping  $f:\ \mathcal{M}\longrightarrow  Z$ is said to be {\it bounded} if
there exist an open neighborhood $U$ of $f(\mathcal{M})$ in $Z$ and
a holomorphic embedding $\phi$ of $ U $ into a polydisc of $ \C^k$  such that
$\phi(U)$ is an analytic set in this polydisc.
      $f$ is said to be {\it locally bounded along} $\mathcal{N}\subset \mathcal{M}$ if
for every point $z\in \mathcal{N},$ there is an open neighborhood $U$ of $z$ (in $\mathcal{M}$) such that  $f|_{U}:\ U
\longrightarrow Z$ is
bounded.   $f$ is said to be {\it locally bounded} if  it is so for  $\mathcal{N}= \mathcal{M}.$  
It is clear that  if $Z=\C$ then the above notions of boundedness coincide with the usual ones.

\subsection{Hartogs extension property.}\label{HEP}
The following example (see Shiffman  \cite{sh2}) shows that an additional 
hypothesis on the ``target space" $Z$ is necessary in order that the  PROBLEM  makes sense. Consider the mapping $f:\
\C^2\longrightarrow \P^1$ given by
\begin{equation*}
  f(z,w):=
\begin{cases}
[(z+w)^2: (z-w)^2],
  & (z,w)\not=(0,0) ,\\
 [1:1], &  (z,w)=(0,0).
\end{cases}
\end{equation*}
Then $f\in\mathcal{O}_s\Big(\X^{\text{o}}(\C,\C;\C,\C),\P^1\Big),$ but $f$ is not
continuous at $(0,0).$

 We recall here  the following notion
(see,  for  example, Shiffman \cite{sh1}).
 Let $p\geq 2$ be an integer. For $0<r<1,$  the {\it Hartogs
 figure} in dimension $p,$ denoted by $H_p(r),$ is given by
 \begin{equation*}
H_p(r):=\left\lbrace (z^{'},z_p)\in E^p: \  \Vert z^{'}\Vert<r \ \ \text{or}\ \ \vert z_p\vert >1-r \right\rbrace,
 \end{equation*}
where $E$ is the open  unit disc of $\C$ and $z^{'}=(z_1,\ldots,z_{p-1}),$
$\Vert z^{'}\Vert:=\max\limits_{1\leq j\leq p-1} \vert z_j\vert.$

\begin{defi}\label{defi_HEP}
A complex analytic space $Z$ is said to {\rm possess the Hartogs extension property
  in dimension $p$} if  every  mapping $f\in\mathcal{O}(H_p(r) , Z)$ extends to a  mapping
  $\hat{f}\in\mathcal{O}(E^p , Z).$  Moreover, $Z$ is said to {\rm possess the Hartogs extension property}
 if it does in any dimension $p\geq 2.$
\end{defi}
It is a classical result of Ivashkovich (see \cite{iv1}) that if $Z$ possesses
the Hartogs extension property  in dimension 2, then it does in all dimensions $p\geq 2.$
 Some typical  examples of complex analytic  spaces  possessing the Hartogs extension property are  the
complex Lie groups (see \cite{asy}), the taut spaces (see \cite{wu}), the Hermitian manifold with
negative holomorphic
sectional  curvature (see \cite{sh1}), the holomorphically convex K\"{a}hler  manifold without rational curves
(see \cite{iv1}).

Here we  mention  an important characterization due to Shiffman (see \cite{sh1}).
\begin{thm}\label{thm_Shiffman}
A complex analytic space $Z$   possesses the Hartogs extension property
    if and only if  for every domain $D$ of any Stein manifold $\mathcal{M},$ every mapping
    $f\in\mathcal{O}(D, Z)$ extends to a  mapping $\hat{f}\in
    \mathcal{O}(\widehat{D},Z) ,$   where $\widehat{D}$ is the envelope of holomorphy
of $D.$
\end{thm}
 
In the light of Definition  \ref{defi_HEP} and   Shiffman's Theorem,
 the natural ``target spaces" $Z$ for obtaining
satisfactory answers to the PROBLEM are the complex analytic spaces  which possess the Hartogs extension property.

\subsection{Statement of the  main results}
We are now ready to state the   main results.

\renewcommand{\thethmspec}{Theorem A}
  \begin{thmspec}
  Let $X,\ Y$  be two complex manifolds,
  let $D\subset X,$ $ G\subset Y$ be two open sets, let
  $A$ (resp. $B$) be a subset of  $\overline{ D}$ (resp.
  $\overline{ G}$).  $D$  (resp.  $G$) is equipped with a 
  system of approach regions 
  $\big(\mathcal{A}_{\alpha}(\zeta)\big)_{\zeta\in\overline{ D},\  \alpha\in I_{\zeta}}$
  (resp.  $\big(\mathcal{A}_{\beta}(\eta)\big)_{\eta\in\overline{ G},\  \beta\in I_{\eta}}$).
   Let $Z$ be a complex analytic  space possessing the  Hartogs extension property.
   Then,
   for every   mapping   $f:\ W\longrightarrow Z$
  which satisfies the following conditions: 
   \begin{itemize}
   \item[$\bullet$]    $f\in\mathcal{C}_s(W,Z)\cap \mathcal{O}_s(W^{\text{o}},Z);$
    \item[$\bullet$] $f$ is locally bounded along $\X\big (A\cap\partial D,B\cap\partial G;D,G\big );$\footnote{
    It follows from  Subsection \ref{subsection_cross}  that $$ \X\big (A\cap\partial D,B\cap\partial G;D,G\big )=
    \big((A\cap\partial D)\times ( G\cup B)\big)\bigcup
    \big( (D\cup A)\times (B\cap\partial G) \big).$$}
    \item[$\bullet$]          $f|_{A\times B}$ is continuous at all points of 
   $(A\cap\partial D)\times (B\cap \partial G),$
    \end{itemize}
     there exists a unique mapping  
$\hat{f}\in\mathcal{O}(\widehat{\widetilde{W}},Z)$ which
 admits  $\mathcal{A}$-limit $f(\zeta,\eta)$ at every point
  $(\zeta,\eta)\in  W\cap   \widetilde{W}.$        

 If, moreover, $Z=\C$ and  $\vert f\vert_W<\infty,$ then
\begin{equation*}
 \vert \hat{f}(z,w)\vert\leq \vert f\vert_{A\times B}^{1-\widetilde{\omega}(z,w)} \vert
 f\vert_W^{\widetilde{\omega}(z,w)},\qquad (z,w)\in\widehat{\widetilde{W}}.
\end{equation*}
\end{thmspec}

Theorem A has an important corollary. Before stating this, we need to
introduce a terminology. A complex manifold $\mathcal{M}$ is said to be
 a {\it Liouville manifold} if $\mathcal{PSH}(\mathcal{M})$ does not
 contain any non-constant bounded above functions. We see clearly that the class
 of Liouville manifolds contains the class of connected compact manifolds.

\renewcommand{\thethmspec}{Corollary B}
  \begin{thmspec}
We keep the hypothesis and the notation in Theorem A.
   Suppose in addition that  $G $ is  a Liouville manifold  and that
   $\widetilde{A},\widetilde{B}\not=\varnothing.$
Then,  
    for every   mapping   $f:\ W\longrightarrow Z$
  which satisfies the following conditions: 
   \begin{itemize}
   \item[$\bullet$]    $f\in\mathcal{C}_s(W,Z)\cap \mathcal{O}_s(W^{\text{o}},Z);$
    \item[$\bullet$] $f$ is locally bounded along  $\X\big (A\cap\partial D,B\cap\partial G;D,G\big );$
    \item[$\bullet$]          $f|_{A\times B}$ is continuous at all points of 
   $(A\cap\partial D)\times (B\cap \partial G),$
    \end{itemize}    
 there is a unique mapping
$\hat{f}\in\mathcal{O}(D\times G,Z)$
which  admits  $\mathcal{A}$-limit $f(\zeta,\eta)$ at every point
  $(\zeta,\eta)\in     W\cap   \widetilde{W}.$        
\end{thmspec}

Corollary B follows immediately from Theorem A since
$\widetilde{\omega}(\cdot,B,G)\equiv 0.$

We will see in Section 10 below that Theorem A and Corollary B generalizes all the results discussed in Section 1 above.
Moreover, they also  give  many new  results. 
Although our main  results have been stated only  for the case  of a $2$-fold cross, 
they can be formulated   for the general case of an $N$-fold cross  with $N\geq 2$  (see also \cite{nv,pn1}).

\section{Holomorphic  discs and a Two-Constant Theorem}
We recall here some elements of Poletsky theory of  discs,  some  background of the
pluripotential theory and auxiliary results     needed for the proof of Theorem A.
\subsection{Poletsky theory of  discs and Rosay Theorem on
 holomorphic  discs}
  Let $E$ denote as usual the open unit disc in $\C.$ %
For a complex
manifold $\mathcal{M},$ let $\mathcal{O}(\overline{E},\mathcal{M})$ denote
the set of all holomorphic mappings $\phi:\ E\longrightarrow \mathcal{M}$ which
extend holomorphically  to   a neighborhood of  $\overline{E}.$
Such a mapping $\phi$ is called a {\it holomorphic disc} on $\mathcal{M}.$ Moreover, for
a subset $A$ of $\mathcal{M},$ let
\begin{equation*}
 1_{  A,\mathcal{M}}(z):=
\begin{cases}
1,
  &z\in   A,\\
 0, & z\in \mathcal{M}\setminus A.
\end{cases}
\end{equation*}

In the work \cite{ro}  Rosay proved the following remarkable result.
\begin{thm}\label{Rosaythm}
Let $u$ be an upper semicontinuous function on a complex manifold
$\mathcal{M}.$ Then the Poisson functional of $u$  defined by
\begin{equation*}
\mathcal{P}[u](z):=\inf\left\lbrace\frac{1}{2\pi}\int\limits_{0}^{2\pi} u(\phi(e^{i\theta}))d\theta:  \
\phi\in   \mathcal{O}(\overline{E},\mathcal{M}), \ \phi(0)=z
\right\rbrace,
\end{equation*}
is plurisubharmonic on $\mathcal{M}.$
\end{thm}

Rosay Theorem may be viewed as an important development in Poletsky
theory of   discs. Observe that special cases of Theorem
\ref{Rosaythm} have been considered by Poletsky (see \cite{po1,po2}),
L\'arusson--Sigurdsson (see \cite{ls}) and Edigarian (see \cite{ed}).

The following Rosay type result  gives the connections between  the Poisson functional of the characteristic function
 $1_{\mathcal{M}\setminus A,\mathcal{M}}$
and holomorphic discs.
\begin{lem}\label{Rosaylem}
Let $\mathcal{M}$ be a complex manifold and let $A$ be a nonempty open
subset of $\mathcal{M}.$   Then
for any $\epsilon>0$ and any $z_0\in \mathcal{M},$  there are an open neighborhood $U$ of
$z_0,$ an open subset $T$ of $\C,$    and a  family of holomorphic discs
 $(\phi_z)_{z\in U}\subset \mathcal{O}(\overline{E} ,\mathcal{M})$ with the
following properties:
\begin{itemize}
\item[(i)] $\Phi\in\mathcal{O}(U\times E ,\mathcal{M}),$ where
$\Phi(z,t):=\phi_z(t),$  $(z,t)\in U\times E ;$
\item[(ii)] $\phi_{z}(0)=z,\qquad z\in U;$
\item[(iii)] $\phi_z(t)\in  A, \qquad t\in T\cap  \overline{E},\ z\in U;$
\item[(iv)] $  \frac{1}{2\pi}\int\limits_{0}^{2\pi}  1_{ \partial E \setminus T,\partial E } (e^{i\theta})d\theta
< \mathcal{P}[1_{\mathcal{M}\setminus A,\mathcal{M}}](z_0 )+\epsilon.$
\end{itemize}
\end{lem}
\begin{proof} See Lemma 3.2 in \cite{nv}.
\end{proof}
The next result describes the situation  in dimension $1.$  It will  be very useful later on.
\begin{lem}\label{lem_circle}
Let $T$ be an open subset of $\overline{E}.$   Then
\begin{equation*}
\omega(0,T\cap E,E)\leq \frac{1}{2\pi}\int\limits_{0}^{2\pi}1_{  \partial E\setminus T,T}
(e^{i\theta})d\theta.
\end{equation*}
\end{lem}
\begin{proof} See, for example,  Lemma 3.3 in \cite{nv}.
\end{proof}
The last result, which is an important consequence of Rosay's Theorem,  gives the connection between
  the Poisson functional       and the plurisubharmonic measure.
\begin{prop}\label{prop_Rosay}
Let $\mathcal{M}$ be a complex manifold and   $A$   a nonempty open
subset of $\mathcal{M}.$    Then
   $ \omega(z,A,\mathcal{M}) = \mathcal{P}[1_{\mathcal{M}\setminus A,\mathcal{M}}](z),$ $z\in\mathcal{M}.$
\end{prop}
\begin{proof} See, for example, the proof of Proposition 3.4 in \cite{nv}.
\end{proof}

\subsection{Level sets of the relative extremal functions and a Two-Constant Theorem}
%
%
\label{subsection_level_sets}
Let $X$ be a complex manifold and  $D\subset X$ an open set.
Suppose that $D$ is equipped with a system of approach regions
 $\mathcal{A}=\big(\mathcal{A}_{\alpha}(\zeta)\big)_{\zeta\in\overline{D},\  \alpha\in I_{\zeta}}.$
For every   open subset $G$ of $D,$ there is a natural  system of approach regions for $G$
which is called  the {\it induced system 
of approach regions $\mathcal{A}^{'}=\big(\mathcal{A}^{'}_{\alpha}(\zeta)\big)_{\zeta\in\overline{ G},\ 
 \alpha\in I^{'}_{\zeta}}$ of $\mathcal{A}$ onto $G$}. It is given by
\begin{equation*}
\mathcal{A}^{'}_{\alpha}(\zeta):=\mathcal{A}_{\alpha}(\zeta)\cap G,\qquad \zeta\in\overline{G},\ 
\alpha\in I^{'}_{\zeta},
\end{equation*}
where $I^{'}_{\zeta}:=\left\lbrace  \alpha\in I_{\zeta}:\ \zeta\in\overline{\mathcal{A}_{\alpha}(\zeta)\cap G}
\right\rbrace.$

\begin{prop}\label{prop_levelsets}
Under the  above  hypothesis  and  notation,
let  $A\subset \overline{D}$  be a    locally pluriregular set (relative to $\mathcal{A}$).
  For $0<\delta<1,$ define the {\rm $\delta$-level
set of $D$ relative to $A$} as follows
\begin{equation*}
 D_{\delta,A}:=
\left\lbrace z\in  D:\ \omega(z,A,D)<1-\delta \right\rbrace.
\end{equation*}
We equip  $D_{\delta,A}$  with the induced system of approach regions
$\mathcal{A}^{'}$ of $\mathcal{A}$ onto $D_{\delta,A}$  (see Subsection 2.1 above).
Then   $A\subset\overline{D_{\delta,A}}$ 
 and 
 \begin{equation}\label{eq3.levelsets}
\omega(z,A,D_{\delta,A})=\frac{\omega(z,A,D)}{1-\delta},\qquad
z\in  D_{\delta,A}.
\end{equation}
Moreover,  $A$ is locally pluriregular relative to $\mathcal{A}^{'}.$
\end{prop}
\begin{proof}
Since  $A$ is locally pluriregular, we see that
\begin{equation}\label{eq3.levelsets.1}
 \big(\Alimsup\omega(\cdot ,A,D)\big)(z)=0,\qquad  z\in A. 
\end{equation}
Therefore, for every $z\in A$ and $ \alpha\in I_{z},$ there is an open  neighborhood  $U$ of $z$ such that
$\varnothing\not=\mathcal{A}_{\alpha}(z)\cap U\subset D_{\delta,A}.$ Hence,  $A\subset\overline{D_{\delta,A}}.$

Next, we turn to the proof of  identity (\ref{eq3.levelsets}). 
Observe that
$0\leq \frac{\omega (\cdot,A,D)}{1-\delta}\leq 1$ on
$D_{\delta,A}$   by definition. This, combined with (\ref{eq3.levelsets.1}), implies that 
\begin{equation} \label{eq3.levelsets.2}
\frac{\omega (z,A,D)}{1-\delta}\leq  \omega(z,A,D_{\delta,A} ),
\qquad z\in D_{\delta,A}.
\end{equation}
To prove the converse inequality of (\ref{eq3.levelsets.2}), let
$u\in\mathcal{PSH}(D_{\delta,A})$ be such that $u\leq 1$ on
$D_{\delta,A}$
and $\mathcal{A}^{'}-\limsup u \leq 0$ on $ A .$
Consider  the following function
\begin{equation}\label{eq3.levelsets.3}
 \hat{u}(z):=
\begin{cases}
\max\left\lbrace (1-\delta)u(z),\omega (z,A,D)\right\rbrace,
  &z\in   D_{\delta,A},\\
 \omega (z,A,D), & z \in  D\setminus D_{\delta,A}.
\end{cases}
\end{equation}
It can be checked that  $\hat{u}\in\mathcal{PSH}(D)$   and $0\leq \hat{u}\leq 1.$ Moreover,
in virtue of the assumption on  $u$ and (\ref{eq3.levelsets.1}) and  (\ref{eq3.levelsets.3}),   we have that
\begin{equation*} (\Alimsup \hat{u})(a)
\leq \max\left\lbrace   (1-\delta)(\mathcal{A}^{'}-\limsup u)(a),
\big (\Alimsup  \omega  (\cdot,A,D)\big)(a) \right\rbrace=0 
\end{equation*}
for all $a\in A.$
Consequently,    $\hat{u}\leq\omega (\cdot,   A,D).$
In particular, one gets  from   (\ref{eq3.levelsets.3}) that
\begin{equation*}
u(z)\leq \frac{ \omega(z,A,D)}{1-\delta},\qquad z\in D_{\delta,A}.
\end{equation*}
Since $u$ is arbitrary, we deduce from the latter estimate that the
converse inequality of (\ref{eq3.levelsets.2}) also holds. This, combined with  (\ref{eq3.levelsets.2}),
  completes the
proof of (\ref{eq3.levelsets}).

To prove the last conclusion of the proposition,
fix a point $a\in A$ and  an open neighborhood $U$ of $a.$ Then we have
\begin{multline*}
\big(\Alimsup\omega(\cdot,A\cap U,D_{\delta,A}\cap U)\big)(a)\leq 
\big(\Alimsup\omega(\cdot,A\cap U,(D\cap U)_{\delta,A\cap U})\big)(a)\\
=\frac{1}{1-\delta}\cdot \big(\Alimsup\omega(\cdot,A\cap U,D\cap U)\big)(a)=0,
\end{multline*}
where the first equality
 follows  from identity (\ref{eq3.levelsets}) and the second one from
 the hypothesis that $A$ is locally pluriregular.
\end{proof}

The following Two-Constant Theorem for plurisubharmonic  functions
will play an important  role in  the proof of the estimate in Theorem A.
\begin{thm}\label{two-constant}
Let $X$ be a   complex manifold and  $D\subset X$  an open subset.  
Suppose that  $D$   is equipped with a 
  system of approach regions 
  $\big(\mathcal{A}_{\alpha}(\zeta)\big)_{\zeta\in\overline{D},\  \alpha\in I_{\zeta}}.$
   Let  $A\subset \overline{D}$  be a    locally pluriregular set.
     Let  $m,M\in\R$  and   $u\in\mathcal{PSH}(D)$ such that
$u(z)\leq M$ for $z\in D,$  and $(\Alimsup u)(z)\leq m$ for $z\in
A.$ 
Then
\begin{equation*}
u(z)\leq m(1- \omega(z,A,D))+M\cdot \omega(z,A,D),\qquad z\in  D.
\end{equation*}
\end{thm}
\begin{proof}
It follows immediately from 
Definition \ref{defi_relative_extremal}. 
\end{proof}

\begin{thm}\label{unique1}
We keep the hypotheses and notation of Theorem \ref{two-constant}.
Let $f$ be a  bounded function in $ \mathcal{O}(D,\C)$  such that
$(\Alim f)(\zeta)=0,$  $\zeta\in A.$
Then $f(z)= 0$  for all $z\in D$ such that  $\omega(z,A,D)\not=1.$  
\end{thm}
\begin{proof}
 Fix a finite  positive constant $M$ such that $\vert f\vert_D<M.$ 
Consequently, the desired conclusion follows from applying  Theorem \ref{two-constant}
to the function $u:=\log\vert f\vert.$
\end{proof}

%
%
%
%
%
\subsection{Construction of discs}
In this subsection we present the construction of discs \`a la Poletsky (see \cite{po2}).
This is one of the main ingredients in the proof of Theorem A. 


Let $\mes$ denote the Lebesgue measure on the unit circle $\partial E.$
For a bounded mapping $\phi\in\mathcal{O}(E,\C^n)$  and $\zeta\in \partial E,$
$f(\zeta)$ denotes the angular limit value of $f$ at $\zeta$ if it exists.
A classical theorem of Fatou  says that $\mes\left( \{\zeta\in\partial E:\ \exists f(\zeta)\} \right)=2\pi.$
 For $z\in\C^n$ and $r>0,$ let $\B(z,r)$ denote the open ball centered  at $z$ with radius $r.$

\begin{thm}\label{Poletsky}
Let $D$ be a bounded  open set in $\C^n,$  $A\subset \overline{D},$  $z_0\in D$  and $\epsilon>0.$
Let  $\mathcal{A}$ be a system of approach regions for $D.$
Suppose in addition that  $A$ is locally pluriregular (relative to $\mathcal{A}$).
Then  there exist a  bounded  mapping $ \phi\in\mathcal{O}(E,\C^n)$ and a measurable subset $\Gamma_0\subset \partial E$ with the following  properties:
\begin{itemize}
\item[1)] $\Gamma_0$ is pluriregular (with respect to  the system of angular  
approach regions), $\phi(0)=z_0,$ $\phi(E)\subset  \overline{D},$ $\Gamma_0 \subset \left\lbrace \zeta\in\partial E:\ \phi(\zeta)\in \overline{A}
\right\rbrace,$  and  
\begin{equation*}
1-\frac{1}{2\pi}\cdot\mes (\Gamma_0 )<\omega(z_0,A,D)+\epsilon.
\end{equation*} 
\item[2)] Let  $f\in\mathcal{C}(D\cup \overline{A},\C)\cap \mathcal{O}(D,\C)$  be such that $f(D)$
 is bounded.  Then there exist a bounded function $g\in\mathcal{O}(E,\C)$ such that
  $g=f\circ \phi$ in  a neighborhood of $0\in E$  and\footnote{ Note  here that  by Part 1),  $(f\circ \phi)(\zeta)$ exists
   for all $\zeta\in \Gamma_0.$} 
 $g(\zeta)=(f\circ \phi)(\zeta)$  for all $\zeta\in \Gamma_0.$ Moreover, $g|_{\Gamma_0}\in\mathcal{C}(\Gamma_0,\C).$
 \end{itemize}
\end{thm}

This theorem motivates the following 
\begin{defi}
\label{candidate}
We keep the hypothesis and notation of Theorem \ref{Poletsky}. Then every pair $(\phi,\Gamma_0)$
satisfying the conclusions 1)--2) of this theorem
 is said to be an {\rm $\epsilon$-candidate for the triplet $(z_0,A,D).$}
\end{defi}

Theorem \ref{Poletsky}  says that  there always exist $\epsilon$-candidates for all triplets $(z,A,D).$

\begin{proof}
First we will construct $\phi.$ To do this
we will construct  by induction a sequence $ (\phi_k)_{k=1}^{\infty}\subset\mathcal{O}(\overline{E},D)$ which approximates $\phi$ as
   $k\nearrow\infty.$ 
This will allow  to define the desired  mapping  as  $\phi:=\lim\limits_{k\to\infty} \phi_k.$
The construction of such a sequence  is divided  into three steps.

For $0<\delta,\ r<1$  
let
\begin{equation}\label{eq1_Poletsky}
\begin{split}
 D_{a,r}&:= D\cap \B(a,r), \qquad  a\in A. \\
A_{a,r,\delta}&:=\left\lbrace  z\in  D_{a,r}:\ \omega(z, A\cap \B(a,r),D_{a,r})<\delta    \right\rbrace,\qquad a\in A,\\
A_{r,\delta}&:=\bigcup\limits_{a\in A} A_{a,r,\delta},
\end{split}
\end{equation}
where  in the second ``$:=$"  $D_{a,r}$ is equipped with the induced system 
of approach regions of $\mathcal{A}$ onto $D_{a,r}$  (see  Subsection \ref{subsection_level_sets} above). 

Suppose without loss of generality that $D\subset \B(0,1).$

\noindent{\bf Step 1:} Construction of $\phi_1.$

Let  $\delta_0:=\frac{\epsilon}{3}$ and $r_0:=1.$  Fix $0<\delta_1 <\frac{\delta_0}{3}$  and $0<r_1<\frac{r_0}{3}.$
Applying Proposition \ref{prop_Rosay}, we obtain  $\phi_1\in \mathcal{O}(\overline{E},D)$ such that
$\phi_1(0)=z_0$ and 
\begin{equation*}
 1-\frac{1}{2\pi} \cdot \mes\Big(\partial E\cap \phi_1^{-1}( A_{r_1,\delta_1})\Big)\leq 
   \omega(z_0, A_{r_1,\delta_1},D)+\delta_0.
\end{equation*}
On the other hand, using (\ref{eq1_Poletsky}) and Definition   \ref{defi_relative_extremal}  and the  hypothesis that  $A$
is  locally  pluriregular, we  obtain 
\begin{equation*}
  \omega(z_0, A_{r_1,\delta_1},D)\leq \omega(z_0,A,D).
\end{equation*}
Consequently, we may choose a subset $\Gamma_1$ of $\Gamma_0:=\partial E\cap \phi_1^{-1}( A_{r_1,\delta_1})$  which consists of
 finite
disjoint closed arcs $(\Gamma_{1j})_{j\in J_1}$ so that  
\begin{equation}\label{eq2_Poletsky}
 1-\frac{1}{2\pi} \cdot \mes(\Gamma_1)< \omega(z_0, A_{r_1,\delta_1},D)+2\delta_0\leq  \omega(z_0,A,D)+2\delta_0,
\end{equation}
and 
\begin{equation*}
\sup\limits_{t,\tau\in \Gamma_{1j}}\vert t-\tau\vert <2\delta_1,\ \sup\limits_{t,\tau\in \Gamma_{1j}}\vert \phi_1(t)-\phi_1(\tau)\vert <2r_1,\qquad j\in J_1.
\end{equation*}
\noindent{\bf Step 2:}  Construction of $\phi_{k+1}$ from $\phi_k$  for all  $k\geq 1.$

By the inductive   construction we have 
 $0<\delta_{k}<\frac{\delta_{k-1}}{3}$ and 
 $0<r_{k}<\frac{r_{k-1}}{3}$
and   $\phi_{k}\in \mathcal{O}(\overline{E},D)$ such that 
$\phi_{k}(0)=z_0$ and  there exists a closed  subset $\Gamma_{k}$ of
 $\partial E\cap \phi_{k}^{-1}( A_{r_{k},\delta_{k}})  \cap  \Gamma_{k-1}$  which consists of finite
closed arcs       $(\Gamma_{k,j})_{j\in J_{k}}$      such that
$\Gamma_{k}$ is relatively compact in the interior of $\Gamma_{k-1},$ and 
\begin{equation}\label{eq3_Poletsky}
 1-\frac{1}{2\pi} \cdot \mes(\Gamma_{k})<  1-\frac{1}{2\pi} \cdot \mes(\Gamma_{k-1})+2\delta_{k-1},
\end{equation}
and
\begin{equation*}
\sup\limits_{t,\tau\in \Gamma_{k,j}}\vert t-\tau\vert <2\delta_{k},\  \sup\limits_{t,\tau\in \Gamma_{k,j}}\vert \phi_{k}(t)-\phi_{k}(\tau)\vert 
<2r_{k},\qquad j\in J_{k},
\end{equation*}
and \begin{equation*}
\vert \phi_{k}-\phi_{k-1}\vert_{\Gamma_{k}} <2r_{k-1}.
\end{equation*}
Here  we  make the convention that  the last  inequality is empty  when   $k=1.$

In particular, we have that  $\phi_k(\Gamma_k)\subset A_{r_k,\delta_k}.$ 
 Therefore, by (\ref{eq1_Poletsky}),
for every $\zeta\in \phi_k(\Gamma_k)$  there  is $a\in A$ such that
$\zeta\in A_{a,r_k,\delta_k},$ that is, 
\begin{equation*}
\omega(\zeta, A\cap\B(a,r_k),D_{a,r_k})<\delta_k.  
\end{equation*}
Using the hypothesis that $A$ is locally pluriregular and (\ref{eq1_Poletsky})  
we see that 
\begin{equation*}
\omega(z,A_{r,\delta}\cap D_{a,r_k},D_{a,r_k})\leq \omega(z, A\cap\B(a,r_k),D_{a,r_k}),
\qquad 0<\delta,r<1. 
\end{equation*}
Consequently, for every  $\zeta\in \phi_k(\Gamma_k)$  there  is $a\in A$ such that
\begin{equation*}
\omega(\zeta,A_{r,\delta}\cap D_{a,r_k},D_{a,r_k}) <\delta_k,
\qquad 0<\delta,r<1. 
\end{equation*}
Using the last estimate and arguing as in \cite[p. 120--121]{po2} (see also the proof of Theorem 1.10.7
in  \cite{jp5} for a nice presentation), we can choose $0<\delta_{k+1}<\frac{\delta_k}{3}$ and 
 $0<r_{k+1}<\frac{r_k}{3}$
and   $\phi_{k+1}\in \mathcal{O}(\overline{E},D)$ such that 
$\phi_{k+1}(0)=z_0,$ and  there exists a closed  subset $\Gamma_{k+1}$ of
 $\partial E\cap \phi_{k+1}^{-1}( A_{r_{k+1},\delta_{k+1}})  \cap  \Gamma_k$  which consists of finite
closed arcs       $(\Gamma_{k+1,j})_{j\in J_{k+1}}$      such that
$\Gamma_{k+1}$ is relatively compact in the interior of $\Gamma_k,$ and
\begin{equation}\label{eq4_Poletsky}
 1-\frac{1}{2\pi} \cdot \mes(\Gamma_{k+1})<  1-\frac{1}{2\pi} \cdot \mes(\Gamma_k)+2\delta_{k},
\end{equation}
and
\begin{equation*}
\sup\limits_{t,\tau\in \Gamma_{k+1,j}}\vert t-\tau\vert <2\delta_{k+1},\  \sup\limits_{t,\tau\in 
\Gamma_{k+1,j}}\vert \phi_{k+1}(t)-\phi_{k+1}(\tau)\vert <2r_{k+1},\qquad j\in J_{k+1},
\end{equation*}
and
\begin{equation*}
\vert \phi_{k+1}-\phi_k\vert_{\Gamma_{k+1}} <2r_k.
\end{equation*}

\noindent{\bf Step 3:}  Construction of $\phi$ from the sequence $(\phi_k)_{k=1}^{\infty}.$

In summary, we have constructed a decreasing sequence  $(\Gamma_k)_{k=1}^{\infty}$ of closed subsets
of $\partial E.$  Consider the new closed set
\begin{equation*}
\Gamma:= \bigcap\limits_{k=1}^{\infty} \Gamma_{k}.
\end{equation*}
 By  (\ref{eq3_Poletsky})--(\ref{eq4_Poletsky}), 
\begin{equation*}
\frac{1}{2\pi}\cdot \mes(\Gamma)=\frac{1}{2\pi} \mes(\Gamma_1)-2\sum\limits_{k=1}^{\infty}\delta_k
> \frac{1}{2\pi} \mes(\Gamma_1)- 3\delta_1. 
\end{equation*}
This, combined with  (\ref{eq2_Poletsky}), implies the following property\\
(i)
 \begin{equation*}
 1-\frac{1}{2\pi} \cdot \mes(\Gamma)   < 1-\frac{1}{2\pi} \cdot \mes(\Gamma_1) +3\delta_1
    \leq  \omega(z_0,A,D)+2\delta_0+3\delta_1  <  \omega(z_0,A,D)+\epsilon.
\end{equation*}
On the other hand, we recall from the above construction the following properties:\\
(ii) $\phi_k(\Gamma)\subset \phi_k(\Gamma_k)\subset A_{r_k,\delta_k}.$\\
(iii)  $\delta_0=\frac{\epsilon}{3},$ $r_0=1,$ $0< \delta_{k+1}< \frac{\delta_{k}}{3},$    $0< r_{k+1}< \frac{r_{k}}{3}$
and                       $    \vert \phi_{k+1}-\phi_k\vert_{\Gamma}\leq
                 \vert \phi_{k+1}-\phi_k\vert_{\Gamma_{k+1}} <2r_k.$\\
(iv)   $\sup\limits_{t,\tau\in \Gamma_{kj}}\vert t-\tau\vert <2\delta_k$ and  
$ \sup\limits_{t,\tau\in \Gamma_{k,j}}\vert \phi_{k}(t)-\phi_{k}(\tau)\vert <2r_{k},\qquad j\in J_{k}.$\\
(v) For every $\zeta\in \Gamma$ there exists a sequence $(j_k)_{k\geq 1}$ such that $j_k\in J_k,$ and
$\zeta$ is an interior point of $\Gamma_{k,j_{k}},$
and $\Gamma_{k+1,j_{k+1}}\Subset \Gamma_{k,j_{k}},$ and  $\zeta= \bigcap\limits_{k=1}^{\infty} \Gamma_{k,j_k}.$

Therefore, we are able to apply the Khinchin--Ostrowski Theorem (see \cite[Theorem 4, p. 397]{go})
to the sequence   $(\phi_k)_{k=1}^{\infty}.$ Consequently, this sequence converges uniformly
on compact subsets of $E$ to a mapping  $\phi\in \mathcal{O}(E,D).$ Moreover, $\phi$ admits (angular) boundary values
at all points of $\Gamma$ and $\phi(\Gamma)\subset \bigcap\limits_{k=1}^{\infty} A_{r_k,\delta_k}\subset \overline{A}.$

Observe that since $\phi_k(0)=\phi(0)=z_0\in D$  and $f\in\mathcal{C}(D\cup \overline{A},\C)\cap \mathcal{O}(D,\C),$ 
 the sequence  $( f\circ \phi_k)_{k=1}^{\infty}$
converges to $f\circ\phi$  uniformly on a neighborhood of $0\in E.$ On the other hand,
$f(D)$ is bounded by the hypothesis. Thus
by Montel Theorem, the  family  $( f\circ \phi_k)_{k=1}^{\infty}\subset  \mathcal{O}(E,\C)  $ is normal.
Consequently,  the sequence  $( f\circ \phi_k)_{k=1}^{\infty}$
converges   uniformly on compact subsets of  $ E.$   Let $g$ be the limit mapping. Then 
  $g\in\mathcal{O}(E,\C)$ and  $g=f\circ \phi$ in  a neighborhood of $0\in E.$  
Moreover, it follows from (i)--(iii) above and 
the hypothesis   $f\in\mathcal{C}(D\cup \overline{A},\C) $
 that $g(\zeta)=(f\circ \phi)(\zeta)$  for all $\zeta\in \Gamma.$
 We deduce from  (iii)--(v) above that  $g|_{\Gamma}\in\mathcal{C}(\Gamma,\C) $     
 Finally, applying Lemma \ref{lem4.1} below we may choose a locally pluriregular
 subset    $\Gamma_0\subset\Gamma$ (relative to the system of angular approach regions)
  such that $\mes(\Gamma_0)=\mes(\Gamma).$ 
Hence, the proof is finished.
\end{proof}

    It is worthy to remark  that
$\phi(E)\subset \overline{D};$ but in general,  $\phi(E)\not\subset D$ !

The last result of this section sharpens Theorem \ref{Poletsky}.
\begin{thm}\label{cor_Poletsky}
Let $D$ be a bounded  open set in $\C^n,$  $A\subset \overline{D},$    and $\epsilon>0.$
Let  $\mathcal{A}$ be a system of approach regions for $D.$
Suppose in addition that  $A$ is locally pluriregular (relative to $\mathcal{A}$).
Then there exists a Borel mapping $\Phi:\  D\times  E\longrightarrow \C^n$ with the following
property:
for every $z\in D,$ there is  a measurable subset $\Gamma_z$ of $\partial E$ such that
$(\Phi(z,\cdot),\Gamma_z)$ is an $\epsilon$-candidate for the triplet $(z,A,D).$
\end{thm}  

Roughly speaking, this result says that one can construct $\epsilon$-candidates for $(z,A,D)$ so that they depend
in a Borel-measurable  way on $z\in D.$
\begin{proof}
  Observe that in  Proposition  \ref{prop_Rosay} we can construct  $\epsilon$-candidates for $(z,A,\mathcal{M})$ so that they depend
in a Borel-measurable  way on $z\in \mathcal{M}.$ Here an {\it $\epsilon$-candidate for  $(z,A,\mathcal{M})$} is 
a holomorphic disc   $\phi\in\mathcal{O}(\overline{E},\mathcal{M}) $ such that  $\phi(0)=z$ and 
$\frac{1}{2\pi}\int\limits_{0}^{2\pi} 1_{\partial E\setminus \phi^{-1}(A), \partial E}(e^{i\theta}) d\theta<\mathcal{P}[1_{\mathcal{M}\setminus A,
\mathcal{M}}](z)+\epsilon.$

Using this we can  adapt the proof of Theorem  \ref{Poletsky} in order to obtain the desired result.
\end{proof}

\section{A mixed cross theorem}
Let  $E$  be  as usual  the open unit disc in $\C.$
Let $B$ be a measurable subset of $\partial E$ and $\omega(\cdot,B,E)$  the relative extremal
function  of $B$ relative to $E$ (with respect to the canonical system of approach regions).
 Then it is well-known (see \cite{ra}) that
\begin{equation}\label{eq4.1}
\omega(z,B,E)=\frac{1}{2\pi}\int\limits_{0}^{2\pi} \frac{1-\vert z\vert^2}{\vert e^{i\theta}-z\vert^2}\cdot
 1_{\partial
E\setminus B,\partial E} (e^{i\theta})d\theta.
\end{equation}
 
The following elementary lemma will be very useful.

\begin{lem}\label{lem4.1}
We keep the above hypotheses and notation.
\begin{itemize}
\item[1)] 
 Let
$u$ be a subharmonic function defined on $E$ with $u\leq 1$ and let $\alpha\in (0,\frac{\pi}{2})$
be such that
\begin{equation*}
 \limsup_{z\to \zeta,\ z\in\mathcal{A}_{\alpha}(\zeta)}
 u(z)\leq 0\qquad\text{for a.e.}\ \zeta\in B,
\end{equation*}
where $\mathcal{A}=(\mathcal{A}_{\alpha}(\zeta))$ is the system of angular approach regions defined in Subsection
\ref{Subsection_Examples}.
Then $u\leq \omega(\cdot,B,E)$ on $E.$
\item[2)] $\omega(\cdot,B,E)$ is also   the relative extremal
function  of $B$ relative to $E$ (with respect to the  system of angular approach regions).  
\item[3)] For all subsets $\mathcal{N}\subset\partial E$ with $\mes(\mathcal{N})=0,$
 $\omega(\cdot,B,E)= \omega(\cdot,B\cup\mathcal{N},E).$
\item[4)]   Let $B^{'}$ be the set of    all density points of $B.$ Then
\begin{equation*}
\lim\limits_{z\to \zeta,\ z\in
\mathcal{A}_{\alpha}(\zeta) }\omega(z,B,E)=0,\qquad \zeta\in B^{'},\ 0<\alpha<\frac{\pi}{2}.
\end{equation*}
In particular,   $B^{'}$ is locally pluriregular (with respect to the  system of angular approach regions).
\item[5)] $ \omega(\cdot,B,E)=\widetilde{\omega}_c(\cdot,B,E)=\widetilde{\omega}_a(\cdot,B,E)$ on $E,$  where
 $\widetilde{\omega}_c(\cdot,B,E)$ (resp. $\widetilde{\omega}_a(\cdot,B,E)$)
  is given by Definition  \ref{defi_pluri_measure} relative to the system of  canonical approach regions
  (resp. angular approach regions).
  \end{itemize} 
\end{lem}
\begin{proof}  It follows immediately from the explicit formula (\ref{eq4.1}).   
\end{proof}

The main ingredient in the proof of Theorem A  is the following  mixed cross theorem.
\begin{thm}\label{mixedcrossthm}
Let $D$ be a complex manifold and $E$  as usual the open  unit disc in $\C.$ 
 $D$  (resp. $E$) is equipped with the canonical system of approach regions
 (resp. the system of angular approach regions).
  Let $A$ be an open subset of  $D$  and $B$  a measurable subset of $\partial E$ such that
  $B$ is locally pluriregular  (relative to  the system of angular approach regions).  For $0\leq\delta<1$ put  
    $G:=\left\lbrace
   w\in E:\ \omega(w,B,E)<1-\delta\right\rbrace.$
   Let $W:= \X(A,B;D,G)$, $W^{\text{o}}:= \X^{\text{o}}(A,B;D,G),$
      and\footnote{ In fact, Theorem 4.10  in \cite{pn2} says that $\omega(\cdot,B,G)= \frac{\omega
     (\cdot,B,E)}{1-\delta} $ on $G,$  where   $\omega(\cdot,B,G)$ is  the  relative extremal  function
     with respect to the  system of angular approach  regions induced  onto $G.$} 
         \begin{equation*}
    \widehat{W}=
     \widehat{\X}(A,B;D,G):=\left\lbrace (z,w)\in D\times G:\ \omega(z,A,D)+\frac{\omega
     (w,B,E)}{1-\delta}<1 \right\rbrace.
     \end{equation*}
Let  $f:\ W\longrightarrow \C$ be   such that
\begin{itemize}
\item[(i)] $f\in\mathcal{O}_s(W^{\text{o}},\C);$
\item[(ii)]  $f$ is  locally  bounded on $W,$  $f|_{A\times B}$ is a Borel function;
\item[(iii)] for all $z\in A,$
\begin{equation*}
\lim\limits_{w\to \eta,\
w\in\mathcal{A}_{\alpha}(\eta) }f(z,w)=f(z,\eta),\qquad \eta\in B,\
0<\alpha<\frac{\pi}{2}.
\end{equation*}
\end{itemize}
Then there is a unique function
$\hat{f}\in\mathcal{O}(\widehat{W},\C)$ such that $\hat{f}=f$ on $A\times G.$ 
Moreover,   $\vert f\vert_W=\vert \hat{f}\vert_{\widehat{W}}.$ 
\end{thm}

The proof of this theorem will occupy the present
and the next sections. Our approach here  avoid completely the  classical
method  of doubly orthogonal  bases of Bergman type. 
For the proof we need the following ``measurable" version of Gonchar's Theorem.

\begin{thm}\label{Gonchar_thm}
  Let $D=G:=E$ be equipped  with the system of angular approach regions.  Let
  $A$ (resp. $B$) be a Borel measurable  subset of  $\partial D$ (resp.
  $\partial G$) such that 
$A$ and $B$ are locally pluriregular and that
  $\mes(A),\ \mes(B)>0.$ 
  Put $W:=\X(A,B;D,G)$  and define $W^{\text{o}},$ $\widehat{W},$ $\omega(z,w)$
  as in Subsection  \ref{subsection_cross}. Let  $f:\ W\longrightarrow \C$ be   such that:
\begin{itemize}
\item[ (i)] $f$ is locally bounded on $W$ and  $f\in \mathcal{O}_s(W^{\text{o}},\C);$
\item[(ii)] $f|_{A\times B}$ is a Borel function; 
\item[ (iii)]   for all $a\in A$  (resp.  $b\in B$), $f(a,\cdot)|_{G}$
(resp.   $f(\cdot,b)|_{D}$)  admits  $\mathcal{A}$-limit\footnote{ that is, the angular limit} $f(a,b)$ at all $b\in B$  (resp. at all $a\in A$).
\end{itemize}
Then  there exists a unique function
$\hat{f}\in\mathcal{O}(\widehat{W},\C)$
which  admits  $\mathcal{A}$-limit  $f(\zeta,\eta)$ at all points
  $(\zeta,\eta)\in W^{\text{o}}.$ If, moreover,   $\vert f\vert_W<\infty,$ then
\begin{equation*}
 \vert \hat{f}(z,w)\vert\leq \vert f\vert_{A\times B}^{1-\omega(z,w)} \vert
 f\vert_W^{\omega(z,w)},\qquad (z,w)\in\widehat{W}.
\end{equation*}
\end{thm}
\begin{proof}
It follows from  Steps 1--3 of Section 6  in \cite{pn2}.
\end{proof}

The above theorem is also true in the context of an $N$-fold cross  $W$ $(N\geq 2).$
We give here a version of a special $3$-fold cross   which is needed for the proof of 
Theorem  \ref{mixedcrossthm}.
\begin{thm}\label{Gonchar_thm_extension}
  Let $D=G:=E$ be equipped  with the system of angular approach regions.  Let
  $A$ (resp. $B$) be a Borel measurable  subset of  $\partial D$ (resp.
  $\partial G$) such that 
$A$ and $B$ are locally pluriregular and
that  $\mes(A),\ \mes(B)>0.$ 
   Define  $W,$ $W^{\text{o}},$ $\widehat{W}$
  as follows:
  \begin{eqnarray*}
   W& =&\X(A,\partial E,B;D,E,G):= A\times\partial E\times (G\cup B) \bigcup
  A\times\overline{ E}\times B\bigcup (D\cup A)\times\partial E\times B,\\
  W^{\text{o}}& =&\X^{\text{o}}(A,\partial E,B;D,E,G):= A\times\partial E\times G\bigcup
  A\times E\times B\bigcup D\times\partial E\times B,\\
  \widehat{W}&=&  \widehat{\X}(A,\partial E,B;D,E,G) :=    \left\lbrace (z,t,w)\in D\times E\times G:\
  \omega(z,A,D)+\omega(w,B,G)<1   \right\rbrace.
   \end{eqnarray*}
   Let  $f:\ W\longrightarrow \C$ be   such that:
\begin{itemize}
\item[ (i)] $f$ is locally bounded on $W$ and  $f\in \mathcal{O}_s(W^{\text{o}},\C)$\footnote{
This notation means that   for all $(a,\lambda)\in A\times\partial  E$  (resp.  $(a,b)\in A\times  B$)  (resp.
$(\lambda,b)\in\partial E\times B$), the function  $f(a,\lambda,\cdot)|_{G}$
(resp.   $f(a,\cdot,b)|_{E}$)  (resp.       $f(\cdot,\lambda,b)|_{D}$)  is holomorphic.};
\item[(ii)] $f|_{A\times\partial E\times  B}$ is a Borel function; 
\item[ (iii)]   for all $(a,\lambda)\in A\times\partial  E$  (resp.  $(a,b)\in A\times  B$)  (resp.
$(\lambda,b)\in\partial E\times B$), $f(a,\lambda,\cdot)|_{G}$
(resp.   $f(a,\cdot,b)|_{E}$)  (resp.       $f(\cdot,\lambda,b)|_{D}$)  admits   the angular limit $f(a,\lambda,b)$
 at all
 $b\in B$  (resp. at all $\lambda\in \partial E$)   (resp.   at all $a\in A$).
\end{itemize}
Then  there exists a unique function
$\hat{f}\in\mathcal{O}(\widehat{W},\C)$
such that
\begin{multline*}
 \lim\limits_{\widehat{W}\ni (z,t,w) \to (\zeta,\tau,\eta), w\in \mathcal{A}_{\alpha}(\eta)}\hat{f} (z,t,w)
 =     \frac{1}{2\pi i}\int\limits_{\partial E}\frac{ f(\zeta,\lambda,\eta) }{\lambda-\tau}d\lambda, \\
  (\zeta,\tau,\eta)\in D\times E\times B,\   0<\alpha<\frac{\pi}{2}.
\end{multline*}
   If, moreover,   $\vert f\vert_W<\infty,$ then
\begin{equation*}
 \vert \hat{f}(z,t,w)\vert\leq \vert f\vert_{A\times\partial E\times  B}^{1-\omega(z,A,D)-\omega(w,B,G)} \vert
 f\vert_W^{\omega(z,A,D)+\omega(w,B,G)},\qquad (z,t,w)\in\widehat{W}.
\end{equation*}
\end{thm}
\begin{proof}
We refer the reader to Subsections 5.2 and 5.3 in  \cite{pn2}.

Let  $\hat{\omega}(\cdot, A,D)$  (resp.  $\hat{\omega}(\cdot, B,G)$)  be the conjugate harmonic  function of
$\omega(\cdot,A,D)$ (resp. $\omega(\cdot,B,G)$ ) such that $\hat{\omega}(z_0,A,D)=0$  (resp.  $\hat{\omega}(w_0,B,G)=0$)
for a certain fixed point $z_0\in D$  (resp.  $w_0\in G$). Thus we  define the holomorphic functions
$g_1(z):=\omega(z,A,D)+i\hat{\omega}(z,A,D),$    $g_2(w):=\omega(w,B,G)+i\hat{\omega}(w,B,G),$         and  
\begin{equation*}
g(z,w):=g_1(z)+g_2(w),\qquad (z,w)\in
D\times G.
\end{equation*}

Each function $e^{-g_1}$  (resp.  $e^{-g_2}$) is bounded on $D$  (resp. on  $G$).
 Therefore, in virtue of \cite[p. 439]{go},  we may define  $e^{-g_1(a)}$
  (resp.   $e^{-g_2(b)}$)
for a.e. $a\in A$ (resp. for a.e. $b\in B$) to be the angular limit of
 $e^{-g_1}$ at $a$   (resp. $e^{-g_2}$  at $b$).

In virtue of (i), for each positive integer $N,$
we define, as in \cite{go1,go2} (see also \cite{pn2}), the {\it Gonchar--Carleman  operator} as follows
\begin{equation}\label{eq1_Gonchar}
K_N(z,t,w)=K_N[f](z,t,w):=\frac{1}{(2\pi
i)^2}\int\limits_{A\times B}e^{-N(g(a,b)-g(z,w))}\frac{f(a,t,b)da db}{(a-z)(b-w)}
 \end{equation}
for  $ (z,t,w)\in D\times\partial E\times G.$
Reasoning as in   \cite{go2} and using (i)--(iii) above, we see  that  the following limit
\begin{equation}
\label{eq2_Gonchar}
K(z,t,w)=K[f](z,t,w):=\lim\limits_{N\to\infty}K_N(z,t,w)
\end{equation}
exists for all points in the set
 $\left\lbrace (z,t,w):\   t\in\partial E,\ (z,w)\in  \widehat{\X}(A,B;D,G)\right\rbrace,$ and
its limit is uniform  on compact subsets of the latter set.

Observe that for $ n=0,1,2,\ldots,$  and $N=1,2,\ldots,$
\begin{equation*}
\int\limits_{\partial E}  t^nK_N(z,t,w) d t=
\frac{1}{(2\pi
i)^2}\int\limits_{A\times B}  \Big(\int\limits_{\partial E}  t^nf(a,t,b) d t\Big)
     \frac{ e^{-N(g(a,b)-g(z,w))}            da db}{(a-z)(b-w)}=0,
\end{equation*}
where the first equality follows from (\ref{eq1_Gonchar}),
the second one from the equality $\int\limits_{\partial E}  t^nf(a,t,b) d t=0$
which itself is an immediate consequence  of (i).
Therefore, we deduce from  (\ref{eq2_Gonchar}) that 
\begin{equation*}
\int\limits_{\partial E}  t^nK(z,t,w) d t=0,\qquad (z,w)\in  \widehat{\X}(A,B;D,G),\  n=0,1,2,\ldots. 
\end{equation*}
On the other hand,
\begin{equation*}
\widehat{W}=\left\lbrace (z,t,w):\   t\in E,\ (z,w)\in  \widehat{\X}(A,B;D,G)\right\rbrace .
\end{equation*}
Hence, we are able to define the desired extension function 
 \begin{equation*}
\hat{f}(z,t,w):=\frac{1}{2\pi i}\int\limits_{\partial E} \frac{ K(z,\lambda,w)}{\lambda-t} d \lambda,\qquad  
(z,t,w)\in \widehat{W}. 
\end{equation*}
Recall from Steps 1--3 of Section 6 in \cite{pn2} that
\begin{equation*}
 \lim\limits_{\widehat{W}\ni (z,w) \to (\zeta,\eta), w\in \mathcal{A}_{\alpha}(\eta)} K (z,t,w)
 =     f(\zeta,t,\eta), \qquad
  (\zeta,t,\eta)\in D\times \partial E\times B,\   0<\alpha<\frac{\pi}{2}.
\end{equation*}
Inserting this into  the  above formula of $\hat{f},$ 
 the desired conclusion
of the theorem  follows. 
\end{proof}

We break the proof of Theorem \ref{mixedcrossthm} into two cases.

\smallskip

\noindent {\bf CASE 1:  $\delta=0$}  (that is $G=E$).

\smallskip

We follow essentially the arguments presented in Section 4 of \cite{nv}. For the sake of
clarity and completeness we give here the most basic arguments.

 We begin the proof with the following lemma.
\begin{lem}
\label{lem_Part1}
We keep the hypothesis of Theorem \ref{mixedcrossthm}. For $j\in \{1,2\},$ let $\phi_j\in\mathcal{O}(
 \overline{E}, D) $ be a holomorphic disc,  and let $t_j\in E$ such that $\phi_1(t_1)=\phi_2(t_2)$ and
$  \frac{1}{2\pi}\int\limits_{0}^{2\pi}1_{D\setminus A,D} (\phi_j(e^{i\theta}))d\theta
<1.$ Then:
\begin{itemize}
\item[1)] For  $j \in \{
1,2\},$  the function $(t,w)\mapsto f(\phi(t),w)$   defined on 
$\X(\phi^{-1}_j(A)\cap \partial E,B;E,G)$    satisfies the hypothesis of Theorem \ref{Gonchar_thm}, where
$\phi^{-1}_j(A):=\lbrace t\in \overline{E}:\ \phi_j(t)\in A\rbrace.$     
\item[2)]  For  $j \in \{
1,2\},$ in virtue of Part 1),
 let $\hat{f}_j$ be the unique function in
$ \mathcal{O}\left(\widehat{\X}(\phi^{-1}_j(A)\cap \partial E,B;E,G),\C\right)$
given by   Theorem    \ref{Gonchar_thm}. Then
\begin{equation*}
\hat{f}_1(t_1,w)=\hat{f}_2(t_2,w),
\end{equation*}
for all $w\in G$ such that
$(t_j,w)\in\widehat{\X}\left(\phi^{-1}_j(A)\cap \partial E,B;E,G\right),$   $j \in \{1,2\}.$
\end{itemize}
\end{lem}

\smallskip
\noindent {\it Proof of Lemma \ref{lem_Part1}.}
 Part 1) follows immediately from the hypothesis. Therefore, it remains to prove Part
 2). To do this fix  $w_0\in G$ such that $(t_j,w_0)\in\widehat{\X}\left(\phi^{-1}_j(A)\cap E
 ,B;E,G\right)$
 for  $j \in \{1,2\}.$
We need to show that $\hat{f}_1(t_1,w_0)=\hat{f}_2(t_2,w_0).$  Observe
that both functions  $w\in\mathcal{G}\mapsto \hat{f}_1(t_1,w)$ and $w\in\mathcal{G}\mapsto\hat{f}_2(t_2,w)$ belong to  $\mathcal{O}
(\mathcal{G},\C),$
where $\mathcal{G}$ is the connected component which contains $w_0$ of the
following open set
\begin{equation*}
  \left\lbrace w\in G:\ \omega(w,B,G)<1-\max\limits_{ j \in \{1,2\}}  \omega(t_j,\phi^{-1}_j(A)\cap\partial E,E)
   \right\rbrace.
   \end{equation*}
 Since
   $\phi_1(t_1)=\phi_2(t_2),$ it follows from Theorem \ref{Gonchar_thm} and the hypothesis of  Part 2)  that
\begin{equation*}
(\Alim\hat{f}_1)(t_1,\eta)=f(\phi_1(t_1),\eta)=f(\phi_2(t_2),\eta)=(\Alim\hat{f}_2)(t_2,\eta) , \qquad \eta\in B.
\end{equation*}
     Therefore, by Theorem  \ref{unique1},   $\hat{f}_1(t_1,w)=\hat{f}_2(t_2,w),$
  $  w\in\mathcal{G}.$ Hence,  $\hat{f}_1(t_1,w_0)=\hat{f}_2(t_2,w_0),$
which completes the proof of the lemma.
\hfill $\square$

\medskip

Now we return to the proof of the theorem in CASE 1 which is divided into two steps.

\smallskip

\noindent{\bf Step 1:} {\it Construction of the extension  function $\hat{f}$ on $\widehat{W}$ and its uniqueness.}

\smallskip

\noindent{ \it Proof of Step 1.}
 We define $\hat{f}$ as follows: Let $\mathcal{W}$ be  the set of all pairs
$(z,w)\in D\times G$  with the property that there are a holomorphic disc
$\phi\in\mathcal{O}(\overline{E},D)$ and $t\in E$ such that $\phi(t)=z$
and  $(t,w)\in\widehat{\X}\left(\phi^{-1}(A)\cap \partial E,B;E,G\right).$  By Part 1) of Lemma \ref{lem_Part1} and Theorem  \ref{Gonchar_thm},
 let $\hat{f}_{\phi}$ be the unique function  in
$\mathcal{O}\left(\widehat{\X}
 (\phi^{-1}(A)\cap \partial E,B;E,G),\C\right)$
  such that 
  \begin{equation}\label{eq4.1.1}
  (\Alim \hat{f}_{\phi})(t,w)= f(\phi(t),w),\qquad  (t,w)\in
  \X^{\text{o}}\left(\phi^{-1}(A)\cap \partial E,B;E,G\right).
  \end{equation}
   Then  the desired extension function $\hat{f}$ is given by
\begin{equation}\label{eq4.1.2}
 \hat{f}(z,w):=\hat{f}_{\phi}(t,w) .
\end{equation}
In virtue of Part 2) of Lemma \ref{lem_Part1}, $\hat{f}$ is well-defined on $\mathcal{W}.$
We next prove that
\begin{equation}
 \label{eq4.1.3}
 \mathcal{W}=\widehat{W}  .
\end{equation}
Taking  (\ref{eq4.1.3}) for granted, then $\hat{f}$ is
well-defined on $\widehat{W}.$  

Now we return to (\ref{eq4.1.3}). To prove the inclusion $\mathcal{W}\subset\widehat{W},$
let  $(z,w)\in\mathcal{W}.$  By the above definition of $\mathcal{W},$ one may find a  holomorphic  disc    $\phi\in
\mathcal{O}(\overline{E}, D),$  a point  $t\in E$ such that $\phi(t)=z$
and  $(t,w)\in\widehat{\X}\left(\phi^{-1}(A)\cap \partial E,B;E,G\right).$
Since  $ \omega(\phi(t),A,D)\leq  \omega(t,\phi^{-1}(A)\cap \partial E,E),$ it follows that
\begin{equation*}
 \omega (z,A,D) +  \omega(w,B,G) \leq   \omega(t,\phi^{-1}(A)\cap \partial E,E)+
  \omega(w,B,G)<1,
\end{equation*}
Hence  $(z,w)\in\widehat{W}.$ This proves the above mentioned inclusion.

To finish the proof of  (\ref{eq4.1.3}), it suffices to show that  $\widehat{W}\subset\mathcal{W}.$
To do this, let $(z,w)\in\widehat{W}$ and fix any $\epsilon>0$ such that
\begin{equation}\label{eq4.1.4}
\epsilon<1-  \omega(z,A,D) -  \omega(w,B,G).
\end{equation}
 Applying  Theorem \ref{Rosaythm} and Proposition \ref{prop_Rosay}, there is a holomorphic  disc   $\phi\in
\mathcal{O}(\overline{E},  D)$   such that $\phi(0)=z$ and
\begin{equation}\label{eq4.1.5}
\frac{1}{2\pi} \int\limits_{0}^{2\pi} 1_{D\setminus
A,D}(\phi(e^{i\theta}))d\theta< \omega(z,A,D)+\epsilon.
\end{equation}
Observe that
\begin{eqnarray*}
  \omega(0,\phi^{-1}(A)\cap \partial E,E)+ \omega(w,B,G)&= &\frac{1}{2\pi} \int\limits_{0}^{2\pi} 1_{D\setminus
A,D}(\phi(e^{i\theta}))d\theta+ \omega(w,B,G)\\
 &< &  \omega(z,A,D) +  \omega(w,B,G)+\epsilon<1,
\end{eqnarray*}
where the equality follows from  (\ref{eq4.1}), the first inequality   holds  by (\ref{eq4.1.5}), and the last
 one by (\ref{eq4.1.4}). Hence,
 $(0,w)\in\widehat{\X}\left(\phi^{-1}(A)\cap \partial E,B;E,G\right),$ which
 implies that $(z,w)\in\mathcal{W}.$ This completes  the proof of  (\ref{eq4.1.3}).
Hence, the construction of  $\hat{f}$ on $\widehat{W}$ has been completed.

Next we show that    $\hat{f}=f$ on $A\times G.$ To this end   let $(z_0,w_0)$ 
be an arbitrary point of $ A\times G.$
   Choose the   holomorphic  disc
$\phi\in \mathcal{O}(\overline{E},D)$ given by $\phi(t):=z_0,$ $
t\in\overline{E}.$ Then  by formula (\ref{eq4.1.2}),
\begin{equation*} 
\hat{f}(z_0,w_0)=\hat{f}_{\phi}(0,w_0)=f(\phi(0),w_0)=f(z_0,w_0).
\end{equation*} 

If $g\in\mathcal{O}(\widehat{W},\C)$  satisfies $g=f$ on $A\times G,$ then we deduce from
  (\ref{eq4.1.1})--(\ref{eq4.1.2}) that $g=\hat{f}.$ This proves the uniqueness of $\hat{f}.$ 
\hfill $\square$

\smallskip


Finally, we conclude  the proof of CASE 1 by the following

\smallskip

\noindent{\bf Step 2:} {\it  Proof of the fact that
$\hat{f}\in\mathcal{O}(\widehat{W},\C).$}  

\smallskip

\noindent{\it Proof of Step 2.} Fix an arbitrary point   $(z_0,w_0)\in
\widehat{W}$ and let  $\epsilon>0$ be so small  such that
\begin{equation}\label{eq4.3.1}
2\epsilon<1-  \omega(z_0,A,D) -  \omega(w_0,B,G).
\end{equation}
Since $ \omega(\cdot,B,G)\in\mathcal{PSH}(G),$   one may find an open neighborhood $V $ of $w_0
$ such that
\begin{equation}\label{eq4.3.2}
   \omega(w,B,D) <  \omega(w_0,B,G)+\epsilon,\qquad  w\in V.
\end{equation}
Let $n$ be the dimension of $D$ at the point $z_0.$
 Applying Lemma \ref{Rosaylem} and Proposition \ref{prop_Rosay},
  we obtain an open set $T$ in $\C,$ an open  neighborhood $U$ of $z_0,$ 
  and a  family of holomorphic discs $(\phi_z)_{z\in U}\subset \mathcal{O}(\overline{E},D)$ with the
following properties:
\begin{eqnarray} \label{eq4.3.3}
\text{the mapping}\ (z,t)\in U\times E &\mapsto& \phi_z(t)\ \text{is
holomorphic};  \\\label{eq4.3.4}
\phi_{z}(0)&=&z,\qquad z\in U;\\
\label{eq4.3.5}\phi_z(t)&\in & A, \qquad t\in T\cap\overline{E},\ z\in U;\\
  \label{eq4.3.6}\frac{1}{2\pi}\int\limits_{0}^{2\pi}1_{\partial E\setminus T,\partial E  } (e^{i\theta})d\theta
&<& \omega(z_0,A,D)+\epsilon.
\end{eqnarray}
By shrinking $U$ (if necessary), we may assume without loss of generality that
in  a chart, $z_0=0\in\C^n$ and 
\begin{equation}\label{eq4.3.7}
U=\left\{z=(z^1,\ldots,z^n)=(z^{'},z^n)\in\C^n:\  z^{'}\in S,\ \vert z^n\vert <2\right\},
\end{equation}
where $S\subset\C^{n-1}$ is an open set. 

Consider the $3$-fold cross (compared with the notation in Theorem     \ref{Gonchar_thm_extension})  
\begin{multline*}
\X\left(T\cap \partial E,U,B;E,U, G\right):=(T\cap \partial E)\times U\times (G\cup B)\\
\bigcup (T\cap \partial E)\times U\times  B \bigcup (E\cup(T\cap\partial E))\times U\times B,
\end{multline*}
and the function $g:\ \X\left(T\cap \partial E,U,B;E,U, G\right)\longrightarrow \C$ given  by
\begin{equation}\label{eq4.3.8}
g(t,z,w):=f(\phi_z(t),w),\qquad (t,z,w)\in\X\left(T\cap \partial E,U,B;E,U, G\right).
\end{equation}
We make the following observations:

  Let $t\in T\cap \partial E.$   Then, in virtue of (\ref{eq4.3.5}) we have $\phi_z(t)\in
 A$ for  $z\in U.$ Consequently, in virtue of  (\ref{eq4.3.3}), (\ref{eq4.3.8}) and the hypothesis
 $f\in\mathcal{O}_s(W^{\text{o}},\C),$ we conclude that
 $g(t,z,\cdot)|_G\in\mathcal{O}(G,\C)$ $\Big($resp.
 $g(t,\cdot,w)|_U\in\mathcal{O}(U,\C)\Big)$ for any $z\in U$ (resp. $w\in B).$
Analogously, for any $z\in U,$  $w\in B,$    we can show that
 $g(\cdot,z,w)|_E\in\mathcal{O}(E,\C).$

 In summary, we have shown that  $g$ is separately holomorphic.
 In addition, it follows from hypothesis (ii) and (\ref{eq4.3.3})--(\ref{eq4.3.5}) that  $g$ is locally bounded and
 $g|_{(T\cap \partial E)\times U\times B}$ is a Borel function. 
 

For $z^{'}\in S$ write  $E_{z^{'}}:= \left\{z=(z^{'},z_n)\in\C^n:\   \vert z_n\vert <1\right\}.$
Then by  (\ref{eq4.3.7}), $\bigcup\limits_{z^{'}\in S} E_{z^{'}}\subset U.$
  Consequently, for all $z^{'}\in S,$ using  hypothesis (iii) we are able to apply
 Theorem  \ref{Gonchar_thm_extension} to $g$ in order to obtain a unique function $\hat{g}\in
\mathcal{O}\left(
 \widehat{\X}\left(T\cap\partial  E,\partial E_{z^{'}},B;E,E_{z^{'}}, G\right), \C\right)$\footnote{
 In fact, we identify $E_{z^{'}}$  with $E$ in an obvious way.}
 such that
\begin{multline*}
\lim\limits_{(t,z,w)\to (\tau,\zeta,\eta),\ w\in\mathcal{A}_{\alpha}(\eta)}   \hat{g}(t,z,w)= \frac{1}{2\pi i}
\int\limits_{\partial  E_{z^{'}}} \frac{g(\tau,\zeta^{'},\lambda,\eta)}{\lambda-\zeta_n} d\lambda,
\\ (\tau,\zeta,\eta)\in  E\times  E_{z^{'}}\times B,\ z^{'}\in S,\  0<\alpha<\frac{\pi}{2}.
\end{multline*}
Using  (\ref{eq4.3.3}) and (\ref{eq4.3.7})--(\ref{eq4.3.8}) and  the Cauchy's formula, we  see that
the right hand side  is  equal to $g( \tau,\zeta,\eta).$ 
Hence, we have  shown that
\begin{equation}\label{eq4.3.9}
\lim\limits_{(t,z,w)\to (\tau,\zeta,\eta),\ w\in\mathcal{A}_{\alpha}(\eta)}   \hat{g}(t,z,w)= 
g (\tau,\zeta,\eta),
\qquad (\tau,\zeta,\eta)\in  E\times  E_{z^{'}}\times B,\ z^{'}\in S,\   0<\alpha<\frac{\pi}{2}.
\end{equation}
Observe that
\begin{equation*}
  \widehat{\X}\left(T\cap \partial E,\partial E_{z^{'}},B;E,E_{z^{'}}, G\right)=\left\lbrace 
   (t,z,w)\in E\times E_{z^{'}}\times G:\
  \omega(t,T\cap \partial E,E)+ \omega(w,B,G)<1\right\rbrace.
\end{equation*}
On the other hand, for any $w\in V,$
\begin{equation}\label{eq4.3.10}
\begin{split}
  \omega(0,T\cap\partial  E,E)+ \omega(w,B,G)&\leq \frac{1}{2\pi} \int\limits_{0}^{2\pi}
 1_{\partial E\setminus T,\partial E}(e^{i\theta})d\theta+ \omega(w_0,B,G)+\epsilon\\
 &<   \omega(z_0,A,D) + \omega(w_0,B,G)+2\epsilon<1,
 \end{split}
\end{equation}
where the first inequality follows from (\ref{eq4.1}) and (\ref{eq4.3.2}),
the second one from (\ref{eq4.3.6}), and the last one from (\ref{eq4.3.1}).
 Consequently,
\begin{equation}\label{eq4.3.11}
 (0,z,w)\in\widehat{\X}\left(T\cap\partial E,\partial E_{z^{'}},B;E,E_{z^{'}}, G\right),
 \qquad (z,w)\in E_{z^{'}}\times V,\ z^{'}\in S.
\end{equation}
 It  follows from   (\ref{eq4.1.2}),  (\ref{eq4.3.4}),  (\ref{eq4.3.5}) and  (\ref{eq4.3.10}) that, for $z^{'}\in S$ and
  $z \in E_{z^{'}} ,$
  $\hat{f}_{\phi_z}$ is well-defined and holomorphic on $\widehat{\X}(T\cap\partial E,B;E,
G),$ and
\begin{equation}\label{eq4.3.12}
 \hat{f}(z,w)=\hat{f}_{\phi_z}(0,w),\qquad w\in V.
\end{equation}
On the other hand, it follows from  (\ref{eq4.1.1}), (\ref{eq4.3.8}) and  (\ref{eq4.3.9}) that
\begin{multline*}
\lim\limits_{(t,w)\to (\tau,\eta),\ w\in\mathcal{A}_{\alpha}(\eta)}\hat{f}_{\phi_z}(t,w)
=\lim\limits_{(t,w)\to (\tau,\eta),\ w\in\mathcal{A}_{\alpha}(\eta)}\hat{g}(t,z,w),\\
(\tau,\eta)\in  E\times B,\ z\in E_{z^{'}},
\ z^{'}\in S,\  0<\alpha<\frac{\pi}{2} .
\end{multline*}
 Since, for fixed $z\in E_{z^{'}},$ the restricted functions
 $(t,w)\mapsto \hat{g}(t,z,w)$ and $\hat{f}_{\phi_z}$  are holomorphic on
  $\widehat{\X}(T\cap \partial E,B;E, G),$
 we deduce from the latter equality and the uniqueness of Theorem \ref{Gonchar_thm} that
\begin{equation*}
\hat{g}(t,z,w)=\hat{f}_{\phi_z}(t,w),\qquad (t,w)\in\widehat{\X}\left(T\cap \partial E,B;E, G\right),
 z \in   E_{z^{'}},\ z^{'}\in S.
\end{equation*}
In particular, using  (\ref{eq4.1.2}), (\ref{eq4.3.11}) and (\ref{eq4.3.12}),
\begin{equation*}
\hat{g}(0,z,w)=\hat{f}_{\phi_z}(0,w)=\hat{f}(z,w),\qquad   (z,w) \in E_{z^{'}}\times V,\  z^{'}\in S .
\end{equation*}
Since we know from (\ref{eq4.3.11}) that $\hat{g}$ is holomorphic  in the variables $z^n$ and $w$ on a neighborhood of $(0,z_0,w_0),$
it follows that $\hat{f}$ is holomorphic in the variables $z^n$ and $w$ on a neighborhood of $(z_0,w_0).$
Exchanging the role of  $z^n$ and any other variable $z^j,$ $j=1,\ldots,n-1,$ we see that
  $\hat{f}$ is separetely holomorphic  on a neighborhood of $(z_0,w_0).$ In addition, $\hat{f}$ is locally bounded.  
Consequently, we  conclude, by the classical  Hartogs extension Theorem, that $\hat{f}$
 is holomorphic on a neighborhood of $(z_0,w_0).$
 Since $(z_0,w_0)\in\widehat{W}$  is arbitrary,
it follows that  $\hat{f}\in\mathcal{O}(\widehat{W},\C).$
 \hfill $\square$

\smallskip

 Combining Steps 1--2, CASE 1 follows.   \hfill $\square$

 \section{Completion of the proof of Theorem \ref{mixedcrossthm}}

In this section  we introduce the new {\it technique of conformal mappings}.
This technique  will allow us to
pass from CASE 1 to the general case.
We recall a notion  from   Definition 4.8 in \cite{pn2} which will be relevant for our further study.
\begin{defi}\label{endpoint}
Let $\mathcal{A}$ be the system of angular approach regions for $E,$
let $\Omega$ be an open subset  of the unit disc $E$  and  $\zeta$    a point in $\partial E.$  Then
the point $\zeta$ is said to be an {\rm end-point} of
$\Omega$ if, for every $0<\alpha<\frac{\pi}{2},$ there is an open
neighborhood $U=U_{\alpha}$ of $\zeta$ such that
$U\cap \mathcal{A}_{\alpha}(\zeta)\subset \Omega.$
The set of all end-points of $\Omega$  is denoted by $\End(\Omega).$
\end{defi}


The main idea of the technique  of conformal mappings is described below.

\begin{prop}
\label{prop_conformal}
 Let $B$  be   a measurable subset of $\partial E$  with $\mes(B)>0.$  For $0\leq\delta<1$ put  
    $G:=\left\lbrace
   w\in E:\ \omega(w,B,E)<1-\delta\right\rbrace.$
Let $\Omega$ be  an arbitrary connected component of  $G.$
Then
\begin{itemize}
\item[1)] $\End(\Omega) $  is a  
 measurable subset of $\partial E$ and
$\mes(\End(\Omega))>0.$  Moreover,  $\Omega$ is   a  simply connected domain.

In virtue of Part 1) and the Riemann mapping theorem,  let $\Phi$ be a conformal mapping of $\Omega$ onto  $E.$
\item[2)] For every  $\zeta\in \End(\Omega),$ there is  $\eta\in\partial E$ such that  
\begin{equation*}
\lim\limits_{z\to\zeta,\ z\in \Omega\cap \mathcal{A}_{\alpha}(\zeta)}\Phi( z)=\eta, \qquad  0<\alpha<\frac{\pi}{2}.
\end{equation*}
$\eta$ is called {\rm the limit of $\Phi$ at the end-point $\zeta$} and it is denoted by $\Phi(\zeta).$
Moreover, $\Phi|_{\End(\Omega)}$  is  one-to-one. 
\item[3)] Let $f$ be a bounded holomorphic function on  $\Omega,$  $\zeta\in\End(\Omega)$ and $\lambda\in\C$ such that
$\lim\limits_{z\to\zeta,\ z\in \Omega\cap \mathcal{A}_{\alpha}(\zeta)}f(z)=\lambda$ for
some $0<\alpha<\frac{\pi}{2}.$ Then $f\circ\Phi^{-1}\in\mathcal{O}(E,\C)$ admits the
angular limit $\lambda$ at $\Phi(\zeta).$
\item[4)] Let $\Delta$ be a subset of  $\End(\Omega)$ such that $\mes(\Delta)=\mes(\End(\Omega)).$ Put
$\Phi(\Delta):=\{\Phi(\zeta),\ \zeta\in \Delta\},$ where $\Phi(\zeta)$ is given by
Part 2). Then $\Phi(\Delta)$ is a measurable  subset  of  of $\partial E$ with $\mes\big(\Phi(\Delta)\big)>0.$
and
\begin{equation*}
\omega(\Phi(z),\Phi(\Delta),E)=\frac{\omega(z,B,E)}{1-\delta},\qquad z\in \Omega.
\end{equation*}
\end{itemize}
\end{prop}
\begin{proof}
The first  assertions of Part 1)  follows from Theorem  4.9  in \cite{pn2}.
To show that  $\Omega$ is simply connected, take an arbitrary  Jordan domain $D$ such that
$\partial D\subset \Omega.$ We  need to prove that $D\subset \Omega.$
Observe that $D\subset E$ and $\omega(z,B,E)<1-\delta$ for all $z\in\partial D\subset \Omega\subset G.$
By the Maximum Principle, we deduce that  $\omega(z,B,E)<1-\delta$ for all $z\in D.$
Hence, $D\subset G,$ which, in turn, implies that  $D\subset \Omega.$ This completes Part 1).

Part 2) follows from the ``end-point" version of Theorem 4.4.13 in \cite{ra}
(that is, we replace the hypothesis ``accessible point" therein by end-point).

  Applying the classical Lindel\"{o}f's Theorem to $f\circ \Phi^{-1}\in\mathcal{O}(E,\C),$ Part 3) follows.

It remains to prove Part 4). A straightforward argument shows that
  $\Phi(\Delta)$ is a 
   measurable subset of $\partial E.$
Next, we  show that
\begin{equation}
\label{eq5.1.1}
 \omega(\Phi(z),\Phi(\Delta),E)\leq\frac{\omega(z,B,E)}{1-\delta},\qquad z\in
\Omega.
\end{equation}
 To do this pick any $u\in\mathcal{PSH}(E)$ such that
$u\leq 1$ and
\begin{equation*}
\limsup\limits_{w\to \eta} u(w)\leq 0,\qquad \eta\in \Phi(\Delta).
\end{equation*} 
Consequently,  Part 2)  gives that
\begin{equation}
\label{eq5.1.2}
\limsup\limits_{z\to\zeta,\
z\in\Omega\cap\mathcal{A}_{\alpha}(\zeta)} u\circ \Phi(z)=0,\qquad \zeta\in \Delta,\ 0<\alpha<\frac{\pi}{2}.
\end{equation}
Next, consider the following function
\begin{equation}\label{eq5.1.3}
 \tilde{u}(z):=
\begin{cases}
\max\{(1-\delta)\cdot(u\circ \Phi)(z),\omega(z,B,E)\},
  &z\in    \Omega,\\
 \omega(z,B,E), & z \in E\setminus \Omega.
\end{cases}
\end{equation}
Then it can be checked that $\tilde{u}$ is subharmonic   and $\tilde{u}\leq 1$ in $E.$ In addition,
  for every density point $\zeta$ of $B$ such that $\zeta\not\in \End(\Omega),$ we know from Theorem 4.9 in \cite{pn2} that
  there is  a connected component  $\Omega_{\zeta}$ of $G$ other than $\Omega$ such that  $\zeta\in \End(\Omega_{\zeta}).$
  Consequently,  
     Part 4) of  Lemma \ref{lem4.1}    gives, for such a point $\zeta,$ that
\begin{equation*}
\limsup\limits_{z\to\zeta,\
z\in \mathcal{A}_{\alpha}(\zeta)} \tilde{u} (z)=
\limsup\limits_{z\to\zeta,\
z\in  \mathcal{A}_{\alpha}(\zeta)}\omega(z,B,E)=0, \qquad 0<\alpha<\frac{\pi}{2} .
\end{equation*}
This, combined with (\ref{eq5.1.2}), implies that
\begin{equation*}
\limsup\limits_{z\to \zeta,\ z\in
\mathcal{A}_{\alpha}(\zeta) }\tilde{u}(z)=0, \qquad 0<\alpha<\frac{\pi}{2},\ \text{for a.e.}\ \zeta\in B.
\end{equation*}
Consequently, applying Part 1) of Lemma \ref{lem4.1} yields that $\tilde{u}\leq \omega(\cdot,B,E)$  on $E.$
Hence, by (\ref{eq5.1.3}), $(u\circ\Phi)(z)\leq \frac{\omega(z,B,E)}{1-\delta},$  $z\in \Omega,$  which completes the proof of
 (\ref{eq5.1.1}).
 In particular, we obtain that $\mes\left(\Phi(\Delta)\right)>0.$

To prove the opposite inequality of (\ref{eq5.1.1}), let $u$ be an arbitrary function
in  $\mathcal{PSH}(E)$ such that
$u\leq 1$ and
\begin{equation*}
\limsup\limits_{z\to \zeta} u(z)\leq 0,\qquad \zeta\in B.
\end{equation*}
  Applying Part 3) to the function $f(z):=z$, we obtain  that 
\begin{equation*}
\limsup\limits_{w\to \eta,\ w\in
\mathcal{A}_{\alpha}(\eta) }\frac{\left(u\circ\Phi^{-1}\right)(w)}{1-\delta}\leq 0, \qquad \eta\in \Phi(\Delta),
\ 0<\alpha<\frac{\pi}{2}.
\end{equation*}
On the other hand, since $u\leq \omega(\cdot,B,E)$ on $E,$ one gets that $\frac{\left(u\circ\Phi^{-1}\right)(w)}
{1-\delta}\leq 1,$
$w\in E.$ 
Therefore, applying Part 1) of Lemma \ref{lem4.1} yields that
\begin{equation*}
\frac{\left(u\circ\Phi^{-1}\right)(w)}{1-\delta}\leq
 \omega(w,\Phi(\Delta),E),\quad w\in E,
\end{equation*}
which, in turn,  implies
 the converse inequality of (\ref{eq5.1.1}). 
Hence,   the proof of Part 4) is complete.
\end{proof}

Now  we are in the position to complete  the proof of Theorem  \ref{mixedcrossthm}:

\smallskip

\noindent {\bf CASE 2:}  $0<\delta<1.$

\smallskip

Let  $(G_k)_{k\in K}$  be  the family of all  connected components of $G,$  where $K$ is an (at most) countable  index set.
By  Proposition \ref{prop_conformal}, we may fix a conformal mapping $\Phi_k$ from $G_k$ onto $E$ for every $k\in K.$
 Put
\begin{equation}\label{eq5.2.1}
\begin{split}
B_k&:=\Big [\Phi_k\big(\End(G_k)\cap B\big)\Big ]^{'},\  W_k:=\X(A,B^{'}_k;D,E),\\
  W_k^{\text{o}}&:=\X^{\text{o}}(A,B^{'}_k;D,E),
\  \widehat{W}_k:=\widehat{\X}(A,B^{'}_k;D,E),\qquad k\in K.
\end{split}
\end{equation}
where $ [ T ]^{'}$ (or simply  $T^{'}$) for $T\subset \partial E$ is,  following the notation of Lemma \ref{lem4.1},
 the set of all density points of $T.$
 
 Recall from the  hypotheses of Theorem  \ref{mixedcrossthm} that for every fixed $z\in A,$
  the  holomorphic  function $f(z,\cdot)|_{G}$  is bounded
and that for every $\eta\in B,$
\begin{equation*}
\lim\limits_{w\to \eta,\ w\in\Omega\cap \mathcal{A}_{\alpha}(\eta)}f(z,w)=f(z,\eta),\qquad 0<\alpha<\frac{\pi}{2}.
\end{equation*}
Consequently,  Part 3) of   Proposition \ref{prop_conformal}, applied to $f(z,\cdot)|_{G_k}$ with $k\in K,$
implies that for every fixed $z\in D,$ $f(z,\Phi_k^{-1}(\cdot))\in\mathcal{O}(E,\C)$  admits  the angular limit $f(z,\eta)$  at $\Phi_k(\eta)$
for all $\eta\in B\cap \End(G_k).$ By Part 1) of that proposition, we know that $\mes\big( B\cap \End(G_k)\big)>0.$
This  discussion and  the hypothesis 
allow us to apply   the result of CASE 1  to the function
$g_k:\  W_k\longrightarrow\C$  defined by
\begin{equation}\label{eq2_new_Imomkulov_thm}
g_k(z,w):=\begin{cases}
f(z,\Phi_k^{-1}(w)),  & (z,w)\in  D\times G_k,\\
f(z,\Phi_k^{-1}(w)) & (z,w)\in  D\times B^{'}_k,
\end{cases}
\end{equation}
where in the second line  we have  used the definition  of $\Phi_k|_{\End(G_k)}$ and  its one-to-one property  proved 
by Part 2) of    Proposition \ref{prop_conformal}.

Consequently, we obtain an extension function $\hat{g}_k\in\mathcal{O}(\widehat{W}_k,\C)$
such that
\begin{equation}\label{eq5.2.3}
  \hat{g}_k(z,w)=g_k(z,w),\qquad    (z,w)\in A\times E .
\end{equation}
Put
\begin{equation*}
\widehat{\mathcal{W}}_k:=\left\lbrace(z,\Phi^{-1}_k(w)),\  (z,w)\in  \widehat{W}_k 
 \right\rbrace,\qquad  k\in K.
\end{equation*}
Observe that the open sets $ (\widehat{\mathcal{W}}_k)_{k\in K}$  are pairwise disjoint. Moreover,
by (\ref{eq5.2.1}),
\begin{multline*}
\bigcup\limits_{k\in K} \widehat{\mathcal{W}}_k= \left\lbrace (z,w)\in D\times E:\  w\in  G_k\  \text{and}\right.\\
\left. \omega(z,A,D)+\omega\Big(\Phi_k(w),
\Phi_k(\End(G_k)),E\Big) <1  
\  \text{for some}\ k\in K        \right\rbrace\\
= \left\lbrace (z,w)\in D\times E:\  w\in  G_k\  \text{and}\ \omega(z,A,D)+\frac{\omega(w,
B,E)}{1-\delta} <1 \  \text{for some}\ k\in K        \right\rbrace\\
=\widehat{W},
\end{multline*}
where the second equality follows from Part 4) of Proposition \ref{prop_conformal}.
Therefore, we can define the desired extension function  $\hat{f}\in\mathcal{O}(\widehat{W},\C)$
  by the formula
\begin{equation*}
\hat{f}(z,w):=\hat{g}_k(z,\Phi_k(w)),\qquad   (z,w)\in  \widehat{\mathcal{W}}_k,\  k\in K. 
\end{equation*} 
This, combined with (\ref{eq5.2.1})--(\ref{eq5.2.3}), implies that  $\hat{f}=f$ on $A\times G.$
The uniqueness of $\hat{f}$ follows from that of $\hat{g}_k,$  $k\in K.$
Hence, the proof of the theorem is complete.
\hfill  $\square$
\section{A local version of Theorem A}
The main  purpose  of the section is to  prove  the following  result.

 \begin{thm}  \label{thm_local_euclidean}
 Let $D\subset \C^n,$ $ G\subset \C^m$ be bounded open sets.
  $D$  (resp.  $G$) is equipped with a 
  system of approach regions 
  $\big(\mathcal{A}_{\alpha}(\zeta)\big)_{\zeta\in\overline{D},\  \alpha\in I_{\zeta}}$
  (resp.  $\big(\mathcal{A}_{\alpha}(\eta)\big)_{\eta\in\overline{G},\  \alpha\in I_{\eta}}$  ).
 Let
  $A$ (resp. $B$) be a subset of  $\overline{ D}$ (resp.
  $\overline{ G}$) such that  $A$ and $B$ are locally pluriregular.
  Put 
    \begin{equation*}
    \begin{split}
    W&:= \X(A,B;D,G),\qquad \overline{W}:=\X(\overline{A},\overline{B};D,G), \\
    \overline{W}^{\text{o}}&:=\X^{\text{o}}(\overline{A},\overline{B};D,G),\qquad
     \widehat{W}:=\widehat{\X}(A,B;D,G).
     \end{split}
    \end{equation*}
   Then,
   for every bounded function  $f:\ \overline{W}\longrightarrow\C$ such that
      $f\in\mathcal{C}_s(\overline{W},\C)\cap \mathcal{O}_s(\overline{W}^{\text{o}},\C)$
      and that $f|_{\overline{A}\times \overline{B}}$ is continuous at all points of $(\overline{A}\cap\partial D)
      \times (\overline{B}\cap\partial G),$
     there exists a unique bounded function 
$\hat{f}\in\mathcal{O}(\widehat{W},\C)$ which
  admits  $\mathcal{A}$-limit $f(\zeta,\eta)$ at all points
  $(\zeta,\eta)\in    W.$        Moreover, 
\begin{equation}\label{eq_thm_local_euclidean}
 \vert \hat{f}(z,w)\vert\leq \vert f\vert_{A\times B}^{1-\omega(z,w)} \vert
 f\vert_{W}^{\omega(z,w)},\qquad (z,w)\in\widehat{W}.
\end{equation}
\end{thm}
The core  of  our unified approach will be  presented in the proof  below.
Our idea  is to use  Theorem \ref{Poletsky} in order  to reduce Theorem  \ref{thm_local_euclidean}
to the case  of bidisk, that is, the case of Theorem   \ref{Gonchar_thm}.
This  reduction is  based  on  Theorem  \ref{mixedcrossthm} and on the technique of level sets. 
\begin{proof}  
It is divided into four steps.

\smallskip

\noindent {\bf Step 1:} Construction of  the desired function $\hat{f}\in\mathcal{O}(\widehat{W},\C)$
and  proof of the estimate
$\vert \hat{f}\vert_{\widehat{W}}\leq\vert f\vert_{\overline{W}}.$ 

\smallskip

\noindent{\it Proof of Step 1.}
We define $\hat{f}$  at an arbitrary point $(z,w)\in \widehat{W}$ as follows:
 Let $\epsilon>0$ be such that
 \begin{equation}\label{eq6_Step1_1}
 \omega(z,A,D)+\omega(w,B,G)+2\epsilon<1.
 \end{equation}
   By Theorem \ref{Poletsky} and Definition \ref{candidate}, there is
 an $\epsilon$-candidate $(\phi,\Gamma)$  (resp. $(\psi,\Delta)$)  for $(z,A,D)$  (resp.  $(w,B,G)$).
 Moreover, using the hypotheses, we see that the function $f_{\phi,\psi},$  defined by
 \begin{equation}\label{eq6_Step1_2}
 f_{\phi,\psi}(t,\tau):=f(\phi(t),\psi(\tau)),\qquad  (t,\tau)\in \X\left(\Gamma,\Delta;E,E\right),
 \end{equation} 
satisfies the hypotheses of Theorem \ref{Gonchar_thm}.
By this theorem, let $\hat{f}_{\phi,\psi}$ be the unique function in $ \widehat{\X}\left(\Gamma,\Delta;E,E\right)$
such that
  \begin{equation}\label{eq6_Step1_3}
   (\Alim \hat{f}_{\phi,\psi})(t,\tau)=f_{\phi,\psi}(t,\tau),\qquad  (t,\tau)\in
  \X^{\text{o}}\left(\Gamma, \Delta;E,E\right).
  \end{equation}
  In virtue of (\ref{eq6_Step1_1}) and Theorem \ref{Poletsky} and Lemma \ref{lem_circle}, $(0,0)\in 
   \widehat{\X}\left(\Gamma,\Delta;E,E\right).$
   Then  we can define  the value    of the desired extension function $\hat{f}$ at $(z,w)$   as  follows
\begin{equation}\label{eq6_Step1_formula_hatf}
 \hat{f}(z,w):=\hat{f}_{\phi,\psi}(0,0) .
\end{equation}

 The remaining part of this step is devoted  to showing
  that $\hat{f}$ is well-defined and holomorphic on $ \widehat{W}.$

To this end we fix an arbitrary point $w_0\in G,$ a number $\epsilon_0:$ $0<\epsilon_0< 1-\omega(w_0,B,G),$
 and
 an arbitrary  $\epsilon_0$-candidate $(\psi_0,\Delta_0)$  for   $(w_0,B,G).$
 
 Let  \begin{equation}\label{eq6_Step1_hatW0}
  \widehat{W}_0:=\left\lbrace (z,\tau)\in D\times E:\  \omega(z,A,D)+\omega(\tau,\Delta_0,E)<1
    \right\rbrace.
\end{equation}
 Inspired by   formula (\ref{eq6_Step1_formula_hatf}) we define a function $\hat{f}_0:\   \widehat{W}_0\longrightarrow\C$
  as follows 
 \begin{equation}\label{eq6_Step1_formula_hatf0}
 \hat{f}_0(z,\tau):=\hat{f}_{\phi,\psi_0}(0,\tau) .
\end{equation}
 Here we have used an  $\epsilon$-candidate $(\phi,\Gamma)$    for $(z,A,D),$
 where $\epsilon$ is  arbitrarily  chosen so  that  $0<\epsilon <  1-\omega(z,A,D)-\omega(\tau,\Delta_0,E).$

Using (\ref{eq6_Step1_2})--(\ref{eq6_Step1_3}) and (\ref{eq6_Step1_formula_hatf0}) and arguing as in Part 2) of Lemma \ref{lem_Part1}, 
one can show that $\hat{f}_0$  is well-defined on $\widehat{W}_{0}.$

For all $0<\delta<1$ let
\begin{equation}\label{eq6_Step1_Adelta}
 A_{\delta}:=\left\lbrace z\in D:\
 \omega(z,A,D)<\delta\right\rbrace\ \text{and}\
  E_{\delta}:=\left\lbrace w\in E:\  \omega(w,\Delta_0,E)<1-\delta\right\rbrace.
\end{equation}
Then  by the  construction  in  (\ref{eq6_Step1_formula_hatf0}), we remark that
$\hat{f}_0(z,\cdot)$  is  holomorphic on $E_{\delta}$ for  every fixed $z\in  A_{\delta}.$
 We are able to define a new function  $\tilde{f}_{\delta}$ on $\X\left(A_{\delta},B;D,E_{\delta}
 \right)$ as follows
\begin{equation}\label{eq6_Step1_formula_tildefdelta}
 \tilde{f}_{\delta}(z,\tau):=
\begin{cases}
 \hat{f}_0(z,\tau)
  & \qquad (z,\tau)\in  A_{\delta}\times E_{\delta}, \\
  f(z,\psi_0(\tau)) &   \qquad (z,\tau)\in D\times \Delta_0        .
\end{cases}
\end{equation} 
  Using the hypotheses on $f$ and the previous  remark, we see  that
  $ \tilde{f}_{\delta}\in \mathcal{O}_s\Big(
\X^{\text{o}}\left(A_{\delta}, B;D,E_{\delta}
 \right),\C\Big).$ 
 
 Observe that $A_{\delta}$ is  an open set in $D.$ Consequently,
 $\tilde{f}_{\delta}$ satisfies the hypotheses of Theorem
 \ref{mixedcrossthm}.
 Applying  this theorem  yields a unique function $\hat{f}_{\delta}\in
 \mathcal{O}\Big(\widehat{\X}\left(A_{\delta}, B;D,E_{\delta}
 \right),\C \Big)$
 such that
 \begin{equation*}
 \hat{f}_{\delta}(z,w)= \tilde{f}_{\delta}(z,w),\qquad (z,w)\in  A_{\delta}\times E_{\delta}.
\end{equation*}
This, combined with (\ref{eq6_Step1_formula_tildefdelta}), implies that
 $\hat{f}_0$ is holomorphic on $A_{\delta}\times G_{\delta}.$
On the other hand, it follows from (\ref{eq6_Step1_hatW0}) and (\ref{eq6_Step1_Adelta}) that  
 \begin{equation*}
  \widehat{W}_0=\widehat{\X}\left(A,\Delta_0;D,E
 \right)=\bigcup\limits_{0<\delta<1}A_{\delta}\times G_{\delta}.
  \end{equation*}
 Hence, $\hat{f}_0\in \mathcal{O}(\widehat{W}_0,\C).$
 
 In summary, we have shown that  $\hat{f}_0,$ given by  (\ref{eq6_Step1_formula_hatf0}),
  is well-defined and holomorphic on $\widehat{W}_0.$

 Now we are able to prove that  $\hat{f},$ given by  (\ref{eq6_Step1_formula_hatf}), is well-defined.
 To this end  we fix an arbitrary point $(z_0,w_0)\in \widehat{W},$ an $\epsilon_0:$ $0<\epsilon_0< 1-\omega(z_0,D,G),$
 and
  two arbitrary  $\epsilon_0$-candidates $(\psi_1,\Delta_1)$ and  $(\psi_2,\Delta_2)$  for   $(w_0,B,G).$
  Let  \begin{equation*} 
  \widehat{W}_j:=\left\lbrace (z,\tau)\in D\times E:\  \omega(z,A,D)+\omega(\tau,\Delta_j,E)<1
    \right\rbrace,\qquad j\in\{1,2\}.
\end{equation*}
 Using  formula (\ref{eq6_Step1_formula_hatf0})  define, for $j\in\{1,2\},$ a function $\hat{f}_j:\   \widehat{W}_j\longrightarrow\C$
  as follows 
 \begin{equation}\label{eq6_Step1_formula_hatfj}
 \hat{f}_j(z,\tau):=\hat{f}_{\phi,\psi_j}(0,\tau) .
\end{equation}
 Here we have used any  $\epsilon$-candidate $(\phi,\Gamma)$    for $(z,A,D)$ with a  suitable  $\epsilon>0.$
 Let $\tau_j\in E$ be such that $\psi_j(\tau_j)=w_0,$  $j\in\{1,2\}.$ Then, in virtue of
  (\ref{eq6_Step1_formula_hatf}) and (\ref{eq6_Step1_formula_hatfj}) and the result of the previous paragraph on the 
  well-definedness of $\hat{f}_0,$ the well-defined property
  of $\hat{f}$ is reduced to showing that
  \begin{equation}\label{eq6_Step1_well_defined_hatf}
 \hat{f}_1(\phi(t),\tau_1)=\hat{f}_2(\phi(t),\tau_2) 
\end{equation} 
  for all $t\in E$ and all  $\epsilon$-candidates $(\phi,\Gamma)$    for $(\phi(t),A,D),$ such that
  \begin{equation*}
  \omega(t,\Gamma,A) <\epsilon:=1-\max\limits_{j\in\{1,2\}}\left\lbrace \omega(\tau_1,\Delta_1,E), \omega(\tau_2,\Delta_2,E)
  \right\rbrace.
  \end{equation*}
 Observe that (\ref{eq6_Step1_well_defined_hatf}) follows from an argument based on    Part 2) of Lemma \ref{lem_Part1}. 
Hence, $\hat{f}$  is well-defined on $\widehat{W}.$

As in (\ref{eq6_Step1_Adelta}), for all $0<\delta<1$ let
\begin{equation}\label{eq6_Step1_Adelta_Bdelta}
\begin{split}
 A_{\delta}&:=\left\lbrace z\in D:\
 \omega(z,A,D)<\delta\right\rbrace,\qquad   B_{\delta}:=\left\lbrace w\in G:\
 \omega(w,B,G)<\delta\right\rbrace,\\
 D_{\delta}&:=\left\lbrace z\in D:\  \omega(z,A,D)<1-\delta\right\rbrace  
 ,\qquad G_{\delta}:=\left\lbrace w\in G:\  \omega(w,B,G)<1-\delta\right\rbrace.
 \end{split}
\end{equation}
Now we combine   (\ref{eq6_Step1_Adelta}) and (\ref{eq6_Step1_Adelta_Bdelta})
 and the result that $\hat{f}_0,$       given by  (\ref{eq6_Step1_formula_hatf0}),
  is well-defined and holomorphic on $\widehat{W}_0,$ and the result 
      that  $\hat{f}$  is well-defined on $\widehat{W}.$ Consequently, we obtain that
\begin{equation*}
\hat{f}(\cdot,w)\in\mathcal{O}(D_{\delta},\C),\qquad  w\in B_{\delta},\ 0<\delta<1.
\end{equation*}
Since the  formula  (\ref{eq6_Step1_formula_hatf}) for $\hat{f}$  is symmetric in two variables $(z,w),$
 one also gets that
 \begin{equation*}
\hat{f}(z,\cdot)\in\mathcal{O}(G_{\delta},\C),\qquad  z\in A_{\delta},\ 0<\delta<1.
\end{equation*}
 Since by (\ref{eq6_Step1_Adelta_Bdelta}),
\begin{equation*}
 \widehat{W}=\bigcup\limits_{0<\delta<1} A_{\delta}\times G_{\delta}
 =\bigcup\limits_{0<\delta<1} D_{\delta}\times B_{\delta},
 \end{equation*}
 it follows  from   the previous conclusions that, for all points $(z,w)\in\widehat{W},$
 there is an open neighborhood $U$ of $z$ (resp. $V$ of $w$)  such that
 $f\in\mathcal{O}_s(\X^{\text{o}}(U,V;U,V),\C).$ By the classical Hartogs extension theorem,
 $f\in\mathcal{O}(U\times V,\C).$  Hence,
     $\hat{f}\in\mathcal{O}(\widehat{W},\C).$
     
On the  other  hand, it follows  from    (\ref{eq6_Step1_formula_hatf}) and  the estimate in Theorem \ref{Gonchar_thm} that  
\begin{equation}\label{eq6_Step3_estimate}
\vert \hat{f}\vert_{\widehat{W}}\leq\vert f\vert_{\overline{W}}.
\end{equation} This completes Step 1. \hfill $\square$
     

\smallskip

\noindent{\bf Step 2:}  $f|_{\overline{A}\times\overline{B}}\in\mathcal{C}( \overline{A}\times\overline{B},\C  ).$

\smallskip

\noindent {\it  Proof of Step 2.}
Using the hypotheses  we only need to check the continuity of  $f|_{\overline{A}\times\overline{B}}$
at every point $(a_0,w_0)\in  \overline{A}\times (G\cap \overline{B})$ and
at every point $(z_0,b_0)\in  (D\cap \overline{A})\times  \overline{B}.$
We will  verify the first assertion. To do this let $(a_k)_{k=1}^{\infty}\subset \overline{A}$
and  $(w_k)_{k=1}^{\infty}\subset (G\cap  \overline{B})$ such that $\lim\limits_{k\to\infty} a_k=a_0$ and
 $\lim\limits_{k\to\infty} w_k=w_0.$ We need to show that  
 \begin{equation}\label{eq1_Step2_euclidean}
 \lim\limits_{k\to\infty} f(a_k,w_k)=f(a_0,w_0).
 \end{equation}
Since  $f|_{\overline{W}}$ is locally bounded, we may choose an open connected neighborhood $V$ of $w_0$ such that
 $\sup\limits_{k\geq 1} \vert f(a_k,\cdot)\vert_V 
 <\infty.$ Consequently, by Montel's Theorem, there is a sequence $(k_p)_{p=1}^{\infty}$ such that $(f(a_{k_p},\cdot))$ 
  converges uniformly on compact subsets
 of $V$ to a function $g\in \mathcal{O}(V).$ Equality (\ref{eq1_Step2_euclidean}) is reduced to showing that $g=f(a_0,\cdot)$ on $V.$
 Since $f\in\mathcal{C}_s(\overline{W},\C),$ we deduce that $f(a_0,\cdot) =g$ on  $B\cap V.$   
On the other hand,  $B\cap V$ is non locally pluripolar because $B$ is locally pluriregular and  $w_0\in\overline{B}.$  Hence, we conclude  by
the uniqueness principle  that
 $g=f(a_0,\cdot)$ on $V.$  \hfill $\square$

\smallskip

\noindent{\bf Step 3:}  $\hat{f}$ admits  $\mathcal{A}$-limit $f(\zeta,\eta)$ at all points
  $(\zeta,\eta)\in    W.$
  
  \smallskip
  
\noindent {\it  Proof of Step 3.}
  To this end  we only need to prove that
 \begin{equation}\label{eq6_Step2_Alim}
 \Big (\Alimsup \vert \hat{f}-f(\zeta_0,\eta_0)\vert\Big)(\zeta_0,\eta_0) <\epsilon_0
  \end{equation}
 for   an arbitrary fixed point $(\zeta_0,\eta_0)\in W$  and an arbitrary
  fixed $0<\epsilon_0<1.$
   Suppose without loss of generality that
  \begin{equation}\label{eq6_Step2_supnorm_f}
  \vert f\vert_{\overline{W}}\leq\frac{1}{2}. 
 \end{equation}

First consider $(\zeta_0,\eta_0)\in A\times B.$
     Since  $f\in\mathcal{C}(\overline{A}\times \overline{B},\C),$ one may find an open neighborhood  $U$ of $\zeta_0$ in $\C^n$
   (resp. $V$ of $\eta_0$ in $\C^m$)  so that
 \begin{equation}\label{eq6_Step2_choice_U_V}
 \left\vert f-f(\zeta_0,\eta_0)\right\vert_{(\overline{A}\cap U)\times(\overline{B}\cap V)}<\frac{\epsilon_0^2}{4}.
 \end{equation}
 Consider the open sets
 \begin{equation}\label{eq6_Step2_D'G'}
 D^{'}:=\left\lbrace z\in D:\ \omega(z,A\cap U,D)<\frac{1}{2}  \right\rbrace
 \quad\text{and}\quad  G^{'}:=\left\lbrace w\in G:\ \omega(w,B\cap V,G)<\frac{1}{2}  \right\rbrace.
 \end{equation}
 In virtue of (\ref{eq6_Step2_supnorm_f})--(\ref{eq6_Step2_D'G'}), an application of Theorem \ref{two-constant}
   gives that
 \begin{equation*}
 \vert f(\zeta,w)-f(\zeta,\eta_0)\vert\leq  (\frac{\epsilon_0^2}{4})^{1-\omega(w,B\cap V,G)}
 \leq \frac{\epsilon_0}{2},\qquad \zeta\in A\cap U,\ w\in G^{'}.
 \end{equation*}
 Hence,
 \begin{equation}\label{eq6_Step2_inequality}
 \left\vert f-f(\zeta_0,\eta_0)\right\vert_{\X(\overline{A}\cap U,\overline{B}\cap V;D^{'},G^{'})}\leq
 \frac{\epsilon_0}{2}.
 \end{equation}
  Consider the function $g:\ \X(\overline{A}\cap U,\overline{B}\cap V; D^{'},G^{'}    )\longrightarrow\C,$ given by
 \begin{equation}\label{eq6_Step2_formula_for_g}
 g(z,w):=f(z,w)-f(\zeta_0,\eta_0).
 \end{equation}
 Applying the result of Step 1, we can construct a function $\hat{g}\in\mathcal{O}(\widehat{\X}(A\cap U,B\cap V;
 D^{'},G^{'})
 ,\C)$
 from $g$ in exactly the same way as
 we obtain $\hat{f}\in\mathcal{O}(\widehat{W},\C)$ from $f.$
Moreover, combining    (\ref{eq6_Step1_formula_hatf}) and (\ref{eq6_Step2_formula_for_g}),     we see that
\begin{equation}\label{eq6_Step2_relation_f_and_g}
\hat{g}=\hat{f}-f(\zeta_0,\eta_0)\qquad \text{on}\ \widehat{\X}(A\cap U,B\cap V;D^{'},G^{'}).
\end{equation}
On the other hand, it follows from   formula (\ref{eq6_Step2_formula_for_g}), estimate (\ref{eq6_Step2_inequality}),
and  estimate (\ref{eq6_Step3_estimate})  that
 \begin{equation*}
 \left\vert \hat{g}\right\vert_{\widehat{\X}(A\cap U,B\cap V;D^{'},G^{'})}\leq\frac{\epsilon_0}{2}.
 \end{equation*}
 This, combined with  (\ref{eq6_Step2_relation_f_and_g})  and (\ref{eq6_Step2_D'G'}), implies that
 \begin{equation*}
\Big(\Alimsup\vert \hat{f}(z,w)-f(\zeta_0,\eta_0)\vert\Big )(\zeta_0,\eta_0)\leq\frac{\epsilon_0}{2}.
\end{equation*}
Hence, (\ref{eq6_Step2_Alim}) follows.
In summary, we have shown that
$
\Alim \hat{f}=f
$ on $A\times B.$

Now it remains  to consider  $(\zeta_0,\eta_0)\in A\times G.$ Using   the last  
limit and arguing as in Step 2, one can show that $
\Alim \hat{f}(\zeta_0,\eta_0)=f(\zeta_0,\eta_0).$ 
\hfill $\square$

\smallskip

\noindent{\bf Step 4:}  Proof of the uniqueness of $\hat{f}$ and (\ref{eq_thm_local_euclidean}).\\
\noindent {\it  Proof of Step 4.}  
To prove  the uniqueness of $\hat{f}$ suppose that $\hat{g}\in\mathcal{O}(\widehat{W},\C)$ is a bounded function which admits $\mathcal{A}$-limit 
$f(\zeta,\eta)$ at all points $(\zeta,\eta)\in W.$ Fix an arbitrary point $(z_0,w_0)\in\widehat{W},$
it suffices to show that $\hat{f}(z_0,w_0)=\hat{g}(z_0,w_0).$
Observe that both functions  $\hat{f}(z_0,\cdot)$ and $\hat{g}(z_0,\cdot)$ are bounded and holomorphic
on the $\delta$-level set of $G$ relative to $B:$
\begin{equation*}
G_{\delta,B}:=\left\lbrace w\in G:\ \omega(w, B,G) <1-\omega(z_0,A,D)           \right\rbrace,
\end{equation*}
where $\delta:=\omega(z_0,A,D).$
On the other hand, they admit $\mathcal{A}$-limit $f(z_0,\eta)$ at all points $\eta\in B.$
Consequently, applying  Proposition \ref{prop_levelsets} and Theorem  \ref{unique1}  yields that
 $\hat{f}(z_0,\cdot)=\hat{g}(z_0,\cdot)$  on $G_{\delta,B}.$
 Hence,   $\hat{f}(z_0,w_0)=\hat{g}(z_0,w_0).$
 
To prove   (\ref{eq_thm_local_euclidean})  fix an arbitrary point $(z_0,w_0)\in\widehat{W}.$
 For every $\eta\in B,$  applying Theorem  \ref{two-constant} to $\log \vert f(\cdot,\eta)\vert$
 defined on $D,$ we obtain that
  \begin{equation}\label{eq6_Step3_estimate_2}
 \vert f(z_0,\eta)\vert\leq  \vert f\vert_{A\times B}^{1-\omega(z_0,A,D)}
  \vert f\vert_{W}^{\omega(z_0,A,D)}.
 \end{equation}
 Applying Theorem  \ref{two-constant} again to $\log \vert\hat{f}(z_0,\cdot)\vert$ defined 
 on $G_{\delta,B}$ of the preceeding paragraph, one gets that
  \begin{equation*}
 \vert \hat{f}(z_0,w_0)\vert\leq  \vert f(z_0,\cdot)\vert_{ B}^{1-\omega(w_0,B,G)}
  \vert \hat{f}\vert_{\widehat{W}}^{\omega(w_0,B,G)}.
 \end{equation*}
Inserting  (\ref{eq6_Step3_estimate}) and (\ref{eq6_Step3_estimate_2}) into the right hand side of the latter estimate,
 (\ref{eq_thm_local_euclidean})   follows. Hence Step 4 is finished.
 \hfill $\square$ 

\smallskip

This completes the proof.
\end{proof}

In the sequel we  will need  the  following  refined  version of Theorem  \ref{thm_local_euclidean}.
 \begin{thm}  \label{new_thm_local_euclidean}
 Let $D\subset \C^n,$ $ G\subset \C^m$ be bounded open sets.
  $D$  (resp.  $G$) is equipped with a 
  system of approach regions 
  $\big(\mathcal{A}_{\alpha}(\zeta)\big)_{\zeta\in\overline{D},\  \alpha\in I_{\zeta}}$
  (resp.  $\big(\mathcal{A}_{\alpha}(\eta)\big)_{\eta\in\overline{G},\  \alpha\in I_{\eta}}$).
 Let
  $A,$ $A_0$ (resp. $B,$  $B_0$) be  subsets of  $\overline{ D}$ (resp.
  $\overline{ G}$) such that  $A_0$ and $B_0$ are locally pluriregular
  and that $\overline{A}_0\subset   A^{\ast}$ and   $\overline{B}_0\subset   B^{\ast}.$
  Put 
    \begin{equation*}
    W:= \X(A,B;D,G)\quad\text{and}\quad
     W_0:=\X(A_0,B_0;D,G).
    \end{equation*}
   Then,
   for every bounded function  $f:\ W\longrightarrow\C$ which   satisfies the  following conditions:
\begin{itemize}
     \item[$\bullet$] $f\in\mathcal{C}_s(W,\C)\cap \mathcal{O}_s(W^{\text{o}},\C);$
      \item[$\bullet$] $f|_{A\times B}$ is continuous at all points of $(A\cap\partial D)
      \times (B\cap\partial G),$
\end{itemize}    
 there exists a unique bounded function 
$\hat{f}\in\mathcal{O}(\widehat{W}_0,\C)$ which
  admits  $\mathcal{A}$-limit $f(\zeta,\eta)$ at all points
  $(\zeta,\eta)\in    W_0.$        Moreover, 
\begin{equation}\label{new_eq_thm_local_euclidean}
 \vert \hat{f}(z,w)\vert\leq \vert f\vert_{A_0\times B_0}^{1-\omega(z,A_0,D)-\omega(w,B_0,G)} \vert
 f\vert_{W}^{\omega(z,A_0,D)+\omega(w,B_0,G)},\qquad (z,w)\in\widehat{W}_0.
\end{equation}
\end{thm}
\begin{proof} 
Using the  hypotheses and  applying  Part 1) of Theorem \ref{thm_local_classical}  below
  we  can extend $f$    to a locally bounded   function (still denoted by) $f$ defined
on $\X( A^{\ast}, B^{\ast},D,G)$  such that
$f\in \mathcal{O}_s\big ( \X^{\text{o}}( A^{\ast}, B^{\ast},D,G),\C  \big)$
and that  $f|_{\X(A^{\ast}\cap D, B^{\ast}\cap G;D,G)}$ is  continuous. 
Therefore, the newly defined function $f$ satisfies 
\begin{equation}\label{eq1_new_thm_local_euclidean}
f(a,b)
:=\lim\limits_{k\to\infty}f(a_k,b),
\end{equation}
where $(a,b)$ is an arbitrary point of $ A^{\ast}\times (G\cup  B^{\ast})$  and  $(a_k)_{k=1}^{\infty}\subset  A^{\ast}$ is an  arbitrary sequence
with $\lim\limits_{k\to\infty} a_k=a.$  
Since $f|_W$  is  bounded, it follows that   the  newly defined function $f$ is also bounded.
In virtue of the definition of $A^{\ast}$ and  $B^{\ast}$  we have 
\begin{equation}\label{eq2_new_thm_local_euclidean}
 \partial D\cap A=\partial D\cap  A^{\ast}\quad\text{and}\quad  \partial G\cap B=\partial G\cap  B^{\ast}.
 \end{equation}
Using the  second $\bullet$ in the  hypotheses and formula (\ref{eq1_new_thm_local_euclidean})
  we see that $f|_{ A^{\ast}\times  B^{\ast}}$ is  continuous at  all points  all
$( \partial D\cap A)\times ( \partial G\cap B).$ 
Consequently,  arguing as in the proof of  Step 2 of Theorem   \ref{thm_local_euclidean}
and using  (\ref{eq2_new_thm_local_euclidean}), 
  we  can show that  $f\in\mathcal{C}\big( A^{\ast}\times  B^{\ast},\C\big).$
In summary,  the newly defined function $f$ which is  defined  and  bounded  
on $\X( A^{\ast}, B^{\ast},D,G)$ satisfies
\begin{equation}\label{eq3_new_thm_local_euclidean}
f\in \mathcal{O}_s\big ( \X^{\text{o}}( A^{\ast}, B^{\ast},D,G),\C  \big)
\quad\text{and}\quad f\in\mathcal{C}\big( A^{\ast}\times  B^{\ast},\C\big).
\end{equation}

Observe that $f$ is only  separately continuous  on $\X(A,B;D,G),$ but it is not necessarily so
on the cross $  \X\big ( A^{\ast}, B^{\ast},D,G  \big).$ 
However, we will show that one can  adapt the argument of  Theorem   \ref{thm_local_euclidean}
in order to prove Theorem  \ref{new_thm_local_euclidean}.

We define $\hat{f}$  at an arbitrary point $(z_0,w_0)\in \widehat{W}_0$ as follows:
 Let $\epsilon>0$ be such that
 \begin{equation*}
 \omega(z_0,A_0,D)+\omega(w_0,B_0,G)+2\epsilon<1.
 \end{equation*}
   By Theorem \ref{Poletsky} and Definition \ref{candidate}, there is
 an $\epsilon$-candidate $(\phi,\Gamma)$  (resp. $(\psi,\Delta)$)  for $(z_0,A_0,D)$  (resp.  $(w_0,B,G)$).
To conclude the proof  we only need to  prove  that the function $f_{\phi,\psi},$  defined by
 \begin{equation*}
  f_{\phi,\psi}(t,\tau):=f(\phi(t),\psi(\tau)),\qquad  (t,\tau)\in \X\left(\Gamma,\Delta;E,E\right),
 \end{equation*} 
satisfies the hypotheses of Theorem \ref{Gonchar_thm}. Indeed, having proved this  assertion,
the proof  will follow  along the same lines as  those  given in  Theorem   \ref{thm_local_euclidean}.
This  assertion is again reduced to showing  that for  each fixed $t\in \Gamma,$ the
function  $ f_{\phi,\psi}(t,\cdot)$  admits the angular limit  $f(\phi(t),\psi(\tau))$
for every point $\tau \in \Delta.$  We  will prove  the  last  claim. 

Using the  first $\bullet$ and  Theorem \ref{Poletsky}, we  see that
for  every $a\in A,$ the  function $f(a,\psi(\cdot))\in \mathcal{O}(E,\C)$  admits the angular limit  $f(a,\psi(\tau))$
for every point $\tau \in \Delta.$
Next, using the  hypothesis $\overline{A}_0\subset A^{\ast}$ we may  choose  a sequence  $(a_k)_{k=1}^{\infty}\subset A\cap A^{\ast}$ such that
$\lim\limits_{k\to\infty} a_k=\phi(t)\in \overline{A}_0.$ 
Observe  from (\ref{eq3_new_thm_local_euclidean})  that for  every $k$  the
uniformly bounded  function $f(a_k,\psi(\cdot))\in \mathcal{O}(E,\C)$  admits the angular limit  $f(a_k,\psi(\tau))$
  and that $\lim\limits_{k\to\infty}f(a_k, \psi(\tau))=f(\phi(t), \psi(\tau))$ for every point $\tau \in \Delta.$
Consequently, by  the Khinchin--Ostrowski Theorem (see \cite[Theorem 4, p. 397]{go}),
the above  claim  follows.
\end{proof}

\section{Preparatory results}
The first result of this section shows that   the two definitions 
of plurisubharmonic measure $\widetilde{\omega}(\cdot,A,D),$
given respectively in  Definition  \ref{defi_pluri_measure}
and  in Subsection 2.1 of \cite{nv}, coincide  in the case when $A\subset D.$ 
\begin{prop}\label{prop_coincidence_two_definitions}
 Let $X$  be a complex manifold and 
   $D\subset X$ an   open set.  $D$   is equipped with the
canonical  system $\mathcal{A}$ of approach regions.  Let
  $A$  be a subset of  $D.$  Then 
  $\widetilde{\omega}(z,A,D)=\omega(z,A^{\ast},D).$
\end{prop}
\begin{proof}
Let  $P\in\mathcal{E}(A).$ Then 
by Definition \ref{defi_pluri_measure},
$P\subset  A^{\ast}$ and $P$ is locally pluriregular.
Hence,  $P\subset  ( A^{\ast})^{\ast}=A^{\ast}.$ 
Since   $P\in\mathcal{E}(A)$  is arbitrary, it follows from Definition \ref{defi_pluri_measure} that 
$\widetilde{A}$ is locally pluriregular and  $\widetilde{A}\subset A^{\ast}.$
In particular,
 $ (\widetilde{A})^{\ast}\subset A^{\ast}$ and
 \begin{equation}\label{prop_coincidence_two_definitions_eq1}
\widetilde{\omega}(z,A,D)=\omega(z,\widetilde{A},D)    \geq\omega(z,A^{\ast},D).
\end{equation}
In the sequel we will show that
\begin{equation}\label{prop_coincidence_two_definitions_eq2}
A^{\ast} \subset (\widetilde{A})^{\ast}.
 \end{equation}
 Taking (\ref{prop_coincidence_two_definitions_eq2}) for granted, we have that
 $A^{\ast} = (\widetilde{A})^{\ast}.$
 Consequently, 
\begin{equation*}
\widetilde{\omega}(z,A,D)=\omega(z,\widetilde{A},D)   \leq \omega(z,A^{\ast},D).
\end{equation*}  
This, coupled with (\ref{prop_coincidence_two_definitions_eq1}), completes the proof.
 
To prove    (\ref{prop_coincidence_two_definitions_eq2}) fix an arbitrary point $a\in A^{\ast}$
and an arbitrary but sufficiently small neighborhood $U\subset X$ of $a$ such that $U$ is biholomorphic
to a bounded open set in $\C^n,$ where $n$ is the dimension of $X$ at $a.$
Since $A^{\ast}$ is a Borel subset of $D,$   Theorem 8.5 in  \cite{bt}
provides a subset $P\subset A^{\ast}\cap U $ of type $\mathcal{F}_{\sigma}$\footnote{ This means that
 $P$ is  a countable (or finite) union of relatively closed  subsets of $U.$}
such that  
\begin{equation}\label{prop_coincidence_two_definitions_eq3}
\omega(z,P,U)=\omega(z,A^{\ast}\cap U,U),\qquad   z\in U.
\end{equation}
Write  $P=\bigcup\limits_{n\geq 1} P_n,$  where $P_n$ is  closed.
Observe that $P_n\cap P_n^{\ast}$ is locally pluriregular,
$P_n\setminus (P_n\cap P_n^{\ast})$ is locally pluripolar and  $\overline{P_n\cap P_n^{\ast}}\subset P_n\subset A^{\ast}\cap P.$
Consequently, $ \bigcup\limits_{n\geq 1}(P_n\cap P_n^{\ast})\subset   \widetilde{A}\cap P$ and 
$
P\setminus   \bigcup\limits_{n\geq 1}(P_n\cap P_n^{\ast})
$
is locally pluripolar.
This  implies that
\begin{equation*}
\omega(z,\widetilde{A}\cap U,U)    \leq\omega\left(z,\bigcup\limits_{n\geq 1}(P_n\cap P_n^{\ast})      ,U\right)=\omega(z,P,U), 
\end{equation*}
where the equality holds by applying   Lemma 3.5.3 in \cite{jp1} and  by using the fact that
  $U$ is biholomorphic to a bounded open set in $\C^n.$
This, combined with    (\ref{prop_coincidence_two_definitions_eq3}) and the assumption $a\in A^{\ast},$ implies that
$\omega(a,\widetilde{A}\cap U,U)=0.$
Thus    (\ref{prop_coincidence_two_definitions_eq2})  follows. 
\end{proof}

The main purpose of this and the next sections is to generalize Theorem  \ref{thm_local_euclidean}
to the case where the ``target space" $Z$ is  an arbitrary complex analytic space possessing the Hartogs extension
property.

 \begin{thm}  \label{thm_local_classical}
 Let $D\subset \C^n,$ $ G\subset \C^m$ be two bounded open sets.
  $D$  (resp.  $G$) is equipped with the  canonical
  system of approach regions.  Let $Z$ be a complex analytic space possessing the Hartogs extension
property. 
   Let
  $A$ (resp. $B$) be a subset of  $ D$ (resp.
  $ G$).  Put 
    $W:= \X(A,B;D,G)$ and
     $\widehat{W}:=\widehat{\X}(A,B;D,G).$
    Let $f\in \mathcal{O}_s(W^{\text{o}},Z).$
  \begin{itemize}
  \item[1)] Then $f$ extends  to a mapping (still denoted by) $f$   defined on $\X^{\text{o}}(A\cup A^{\ast},B\cup B^{\ast};D,G)$
such that $f$ is  separately holomorphic   on $\X^{\text{o}}(A\cup A^{\ast},B\cup B^{\ast};D,G)$
and that   $ f|_{\X^{\text{o}}( A^{\ast}, B^{\ast};D,G)}$ is    continuous.
  \item[2)] Suppose in addition that
     $A$ and $B$ are locally pluriregular.
   Then $f$
  extends to a unique mapping 
$\hat{f}\in\mathcal{O}(\widehat{W},Z)$ such that 
  $\hat{f}=f$ on $ W.$      
\end{itemize}
\end{thm}
\begin{proof}
This result has already been  proved in  Th\'eor\`eme 2.2.4 in \cite{az}
starting from Proposition 3.2.1  therein. In the latter proposition Alehyane and Zeriahi
make use of  the method of  doubly orthogonal bases  of Bergman type.
We can avoid this method by simply  replacing every application of  this proposition  by  Theorem \ref{thm_local_euclidean}.
Keeping this change in mind and  using Proposition \ref{prop_coincidence_two_definitions},
 the remaining part of the proof follows along the same
lines as that of Th\'eor\`eme 2.2.4 in \cite{az}.
\end{proof}

\begin{thm}\label{thm_vietanh}
Let  $D,$ $G$ be  complex manifolds,    and let $ A\subset
D,$  $B\subset G$  be open  subsets.  Let $Z$ be a complex analytic space possessing the Hartogs extension
property.  Put $W:=\X(A,B;D,G)$ and
$\widehat{W}:=\widehat{\X}(A,B;D,G).$
Then for any mapping $f\in\mathcal{O}_s(W,Z),$ there is a unique   mapping
$\hat{f}\in\mathcal{O}(\widehat{W},Z)$
such that $\hat{f}=f$ on $W.$
\end{thm}
\begin{proof} It has already been proved in Theorem 5.1 of \cite{nv}. The only places where 
the  method of  doubly orthogonal bases  of Bergman type is involved is the applications
of  Th\'eor\`eme 2.2.4 in \cite{az}. As we already pointed out in Theorem \ref{thm_local_classical},
one can avoid this method by using Theorem \ref{thm_local_euclidean} instead.
\end{proof}

We are ready to formulate a slight generalization of Theorems  \ref{new_thm_local_euclidean}  
   and  \ref{thm_local_classical}.
\begin{thm}  \label{thm_local_localized}
 Let $D\subset \C^n,$ $ G\subset \C^m$ be bounded open sets.
  $D$  (resp.  $G$) is equipped with a 
  system of approach regions 
  $\big(\mathcal{A}_{\alpha}(\zeta)\big)_{\zeta\in\overline{D},\  \alpha\in I_{\zeta}}$
  (resp.  $\big(\mathcal{A}_{\beta}(\eta)\big)_{\eta\in\overline{G},\  \beta\in I_{\eta}}$).
 Let
  $A$  and $A_0$ (resp. $B$ and $B_0$) be two subsets of  $\overline{ D}$ (resp.
  $\overline{ G}$) such that  $A_0$ and $B_0$ are locally pluriregular
  and  that $\overline{A}_0\subset  A^{\ast}$ and  $\overline{B}_0\subset  B^{\ast}.$
   Let $Z$ be a complex analytic space possessing the Hartogs extension
property. Put 
    \begin{equation*}
    W:= \X(A,B;D,G)\quad\text{and}\quad
     W_0:=\X(A_0,B_0;D,G).
    \end{equation*} Then,
   for every bounded mapping  $f:\ W\longrightarrow Z$ which   satisfies the  following conditions:
\begin{itemize}
     \item[$\bullet$] $f\in\mathcal{C}_s(W,Z)\cap \mathcal{O}_s(W^{\text{o}},Z);$
      \item[$\bullet$] $f|_{A\times B}$ is continuous at all points of $(A\cap\partial D)
      \times (B\cap\partial G),$
\end{itemize}    
 there exists a unique bounded mapping 
$\hat{f}\in\mathcal{O}(\widehat{W}_0,\C)$ which
  admits  $\mathcal{A}$-limit $f(\zeta,\eta)$ at all points
  $(\zeta,\eta)\in    W_0.$      

\end{thm}
\begin{proof}
Since $f$ is bounded, one may find an open neighborhood $U$ of $f( \overline{W})$ in $Z$ and
a holomorphic embedding $\phi$ of $U$ into the polydisc $E^k$ of $ \C^k$   such that
$\phi(U)$ is an analytic set in $E^k.$ Now
we are able to apply Theorem    \ref{new_thm_local_euclidean}  to the mapping
$\phi\circ f:\   \overline{W}\longrightarrow \C^k.$ Consequently,
one obtains a unique  bounded mapping $F\in \mathcal{O}(\widehat{W},\C^k)$
 which
  admits  $\mathcal{A}$-limit $(\phi\circ f)(\zeta,\eta)$ at all points
  $(\zeta,\eta)\in    W.$       Using estimate (\ref{new_eq_thm_local_euclidean}) one can show that
  $F\in \mathcal{O}(\widehat{W},E^k).$ Now using Theorem \ref{unique1} it is not difficult to see that
 $F(\widehat{W})\subset\phi(U).$ Consequently, one can define the desired extension mapping
  $\hat{f}$ as follows:
  \begin{equation*}
\hat{f}(z,w):=(\phi^{-1} \circ F)(z,w),\qquad (z,w)\in\widehat{W}. 
  \end{equation*}
\end{proof}

The  following Uniqueness Theorem for holomorphic mappings  generalizes Theorem \ref{unique1}.
\begin{thm}\label{unique2}
Let $X$ be a   complex manifold, $D\subset X$  an open subset
and $Z$  a complex analytic space.  
Suppose that  $D$   is equipped with a 
  system of approach regions 
  $\Big(\mathcal{A}_{\alpha}(\zeta)\Big)_{\zeta\in\overline{D},\  \alpha\in I_{\zeta}}.$
   Let  $A\subset \overline{D}$  be a    locally pluriregular set.
Let $f_1,\ f_2:\   D\cup A\longrightarrow Z$ be locally bounded mappings such that $f_1|_D,f_2|_D\in \mathcal{O}(D,Z)$ and
 $\Alim f_1=\Alim f_2$ on $A.$
 Then $f_1(z)=f_2(z)$ for all $z\in D$ such that
$\omega(z,A,D)\not=1.$
 \end{thm}
We leave the proof to the interested reader. Finally,
we conclude this section with the following  Gluing Lemma.

\begin{lem}  \label{gluing_lem}
 Let $\mathcal{D}$ and $\mathcal{G}$ be open subsets of some complex manifolds
 and  $Z$   a complex analytic space. Suppose that 
   $\mathcal{D}$  (resp.  $\mathcal{G}$) is equipped with a 
  system of approach regions 
  $\big(\mathcal{A}_{\alpha}(\zeta)\big)_{\zeta\in\overline{\mathcal{D}},\  \alpha\in I_{\zeta}}$
  (resp.  $\big(\mathcal{A}_{\beta}(\eta)\big)_{\eta\in\overline{\mathcal{G}},\  \beta\in I_{\eta}}$). 
  Let  $(\mathcal{D}_k)_{k=k_0}^{\infty}$  (resp. $(\mathcal{G}_{k})_{k=k_0}^{\infty}$)
 be  a family of  open subsets of $\mathcal{D}$ (resp. $\mathcal{G}$) equipped with the induced system of  approach regions. 
Let $(\mathcal{P}_k)_{k=k_0}^{\infty}$  (resp. $(\mathcal{Q}_k)_{k=k_0}^{\infty}$)
be a family of locally pluriregular  subsets of $\overline{\mathcal{D}}$ (resp. 
 $\overline{\mathcal{G}}$).
Suppose,  in addition, that
\begin{itemize}
\item[(i)] 
 $\mathcal{P}_{k}\subset \mathcal{P}_{k_0},$   $ \mathcal{D}_{k_0}\subset
 \mathcal{D}_{k},$ and $\mathcal{P}_{k}$
 is locally pluriregular relative to $ \mathcal{D}_{k_0}.$
 Similarly,
   $\mathcal{Q}_{k}\subset \mathcal{Q}_{k_0},$  $ \mathcal{G}_{k_0}\subset
 \mathcal{G}_{k},$ and   $\mathcal{Q}_{k}$ is locally pluriregular relative to    $ \mathcal{G}_{k_0}. $
\item[(ii)]   There are  a family of locally bounded mappings  $(f_k)_{k=k_0}^{\infty}$ such that
$f_k:\  \X^{\text{o}}\left(\mathcal{P}_{k} ,\mathcal{Q}_{k};\mathcal{D}_{k},\mathcal{G}_{k}\right)\longrightarrow Z$
verifies $f_k=f_{k_0}$ on 
$\X^{\text{o}}\left(\mathcal{P}_{k} ,\mathcal{Q}_{k};\mathcal{D}_{k_0},\mathcal{G}_{k_0}\right),$
and a family of holomorphic mappings  $(\hat{f}_k)_{k=k_0}^{\infty}$       such that
   $\hat{f}_{k} \in   \mathcal{O}\Big(\widehat{\X}
 \left(\mathcal{P}_{k} ,  \mathcal{Q}_{k}  ;   \mathcal{D}_{k}, \mathcal{G}_{k}\right),Z\Big), $
  and
 \begin{equation*} 
(\Alim\hat{f}_{k})(z,w)= f_k(z,w) ,\qquad  (z,w)\in
\X^{\text{o}}\left(\mathcal{P}_{k} ,\mathcal{Q}_{k};\mathcal{D}_{k_0},\mathcal{G}_{k_0}\right).
\end{equation*}
\item[(iii)] 
 There are  open subsets  $U$  of $\mathcal{D}$  and $V$ of $\mathcal{G}$ such that 
  $
 \widetilde{\omega}(z, \mathcal{P}_{k},\mathcal{D}_{k_0})+
 \widetilde{\omega}(w, \mathcal{Q}_{k},\mathcal{G}_{k_0})<1
 $
for  all $(z,w)\in U\times V$ and $k\geq k_0.$
\end{itemize}
 Then   $
\hat{f}_{k}(z,w)= \hat{f}_{k_0}(z,w) $ for all $  (z,w)\in
    U\times  V $ and  $k\geq k_0.$
\end{lem}
\begin{proof}
By (iii), we have that 
\begin{equation}\label{eq_1_gluing_lem}
U\times V \subset H:= \widehat{\X}\left(\mathcal{P}_{k} ,\mathcal{Q}_{k};\mathcal{D}_{k_0},\mathcal{G}_{k_0}\right).
\end{equation}
 On the other hand, using (i)   we see that
\begin{equation} \label{eq_2_gluing_lem}
 H \subset   \widehat{\X}\left(\mathcal{P}_{k} ,\mathcal{Q}_{k};\mathcal{D}_{k},\mathcal{G}_{k}
 \right)\cap
   \widehat{\X}\left(\mathcal{P}_{k_0} ,\mathcal{Q}_{k_0};\mathcal{D}_{k_0},\mathcal{G}_{k_0}\right).
\end{equation}
Fix  arbitrary   $(z_0,w_0)\in H$ and $k\geq k_0.$  
Observe that both  mappings $\hat{f}_{k}(\cdot,w_0) $ and $\hat{f}_{k_0} (\cdot,w_0)$ are
defined on $\left\lbrace z\in \mathcal{D}_{k_0}:\  \omega(z,\mathcal{P}_{k},\mathcal{D}_{k_0})
<1- \omega(w_0,\mathcal{Q}_{k},\mathcal{G}_{k_0})\right\rbrace.$
Using (ii) and  Proposition \ref{prop_levelsets}, we may apply Theorem  \ref{unique2} to these mappings
and conclude that   $\hat{f}_{k}(z_0,w_0) =\hat{f}_{k_0} (z_0,w_0).$
 \end{proof}
\section{Local and semi-local versions of Theorem A}

The aim of this section is to generalize  Theorem  \ref{new_thm_local_euclidean} to some cases where the ``target space" $Z$ is 
 a complex analytic  space possessing the  Hartogs extension property. Our philosophy is the following: we first
  apply  Theorem  \ref{new_thm_local_euclidean} 
 locally in order to obtain various local extension mappings, then we glue them together.   
The gluing  process needs the following
\begin{defi}\label{defi_collection}
Let $\mathcal{M}$ be a complex manifold and  $Z$  a complex space. Let $(U_j)_{j\in J}$ be
a family of open subsets of $\mathcal{M},$ and 
 $(f_j)_{j\in J}$ 
a family of   mappings such that $f_j\in \mathcal{O}(U_j,Z).$
We say that the family $(f_j)_{j\in J}$ is {\bf collective} if, for any $j,k\in J,$    $f_j=f_k$ on $U_j\cap U_k.$
The unique holomorphic mapping $f: \bigcup\limits_{j\in J}U_j\longrightarrow Z,$ defined by
$f:=f_j$ on $U_j,$  $j\in J,$ is called  {\bf the collected mapping} of $(f_j)_{j\in J}.$
\end{defi}

We arrive at the following local version of Theorem A. 
 \begin{thm}  \label{thm_local}
  Let $D\subset \C^p,$ $ G\subset \C^q$ be bounded open sets and  $Z$  
 a complex analytic  space possessing the  Hartogs extension property.
  $D$  (resp.  $G$) is equipped with a 
  system of approach regions 
  $\big(\mathcal{A}_{\alpha}(\zeta)\big)_{\zeta\in\overline{D},\  \alpha\in I_{\zeta}}$
  (resp.  $\big(\mathcal{A}_{\beta}(\eta)\big)_{\eta\in\overline{G},\  \beta\in I_{\eta}}$).   
   Let
  $A,$ $A_0$ (resp. $B,$  $B_0$) be  subsets of  $\overline{ D}$ (resp.
  $\overline{ G}$) such that  $A_0$ and $B_0$ are locally pluriregular
  and that $\overline{A}_0\subset   A^{\ast}$ and   $\overline{B}_0\subset   B^{\ast}.$
  Put 
    \begin{equation*}
    W:= \X(A,B;D,G)\quad\text{and}\quad
     W_0:=\X(A_0,B_0;D,G).
    \end{equation*}
   Then,
   for every  mapping  $f:\ W\longrightarrow  Z$ which   satisfies the  following conditions:
\begin{itemize}
     \item[$\bullet$] $f\in\mathcal{C}_s(W,Z)\cap \mathcal{O}_s(W^{\text{o}},Z);$
    \item[$\bullet$] $f$  is  locally  bounded  along  $\X\big (A\cap \partial D, B\cap \partial G;D,G\big);$
      \item[$\bullet$] $f|_{A\times B}$ is continuous at all points of $(A\cap\partial D)
      \times (B\cap\partial G),$
\end{itemize}    
    there exists a unique  mapping 
$\hat{f}\in\mathcal{O}(\widehat{W}_0,Z)$ which
  admits  $\mathcal{A}$-limit $f(\zeta,\eta)$ at all points
  $(\zeta,\eta)\in    W_0.$       
\end{thm}

Theorem  \ref{thm_local} generalizes Theorem  \ref{new_thm_local_euclidean} to the case where the ``target space" $Z$ is 
 an arbitrary  complex analytic  space possessing the  Hartogs extension property.
Since  the  proof is somewhat technical, the  reader  may  skip it at the  first  reading.
\begin{proof}
Recall that for $a\in\C^k$ and $r>0,$ $\B(a,r)$ denotes the open ball centered  at $a$ with radius $r.$
For $0<\delta<1$ and $0<r$   put
 \begin{equation}\label{eq1_thm_local}
 \begin{split}
 D_{a,\delta,r}&:=\left\{z\in  D\cap \B(a,r):\ \omega(A_0\cap\B(a,r),D\cap\B(a,r))<\delta\right\}  
        ,\qquad a\in\overline{A}_0,\\
  G_{b,\delta,r}&:=\left\{w\in  G\cap \B(b,r):\ \omega(B_0\cap\B(b,r),G\cap\B(b,r))<\delta\right\} ,\qquad
 b\in\overline{B}_0.
 \end{split}
 \end{equation}
 Applying Part 1) of Theorem   \ref{thm_local_classical} and using the  hypotheses on $f,$
we  see that $f$ extends  to a  mapping  defined on  $\X(A\cup A^{\ast},B\cup B^{\ast};D,G)$ such that
$f$ is separately holomorphic   on $\X^{\text{o}}(A\cup A^{\ast},B\cup B^{\ast};D,G)$  and that
$f|_{\X( A^{\ast}, B^{\ast};D,G)}$ is  locally bounded.

Therefore, using the compactness of $\overline{A}_0$  and $\overline{B}_0,$ 
one may find a real number $r_0>0$ such that 
 \begin{equation}\label{eq_thm_local_fab}
 f_{a,b}:=f|_{\X\left(A_0\cap\B(a,r),B_0\cap \B(b,r);D \cap\B(a,r)  ,G\cap\B(b,r)\right)}
 \end{equation}
  is  bounded  for all $0<r\leq r_0$  and $a\in\overline{A}_0,$  $ b\in\overline{B}_0.$    Applying Theorem  \ref{thm_local_localized} to $f_{a,b}$ , one obtains a mapping
 \begin{equation}\label{eq_thm_local_hatf_ab}
 \hat{f}_{a,b}\in \mathcal{O}\left(\widehat{\X} \Big(A_0\cap\B(a,r),B_0\cap \B(b,r);D\cap\B(a,r)  ,G\cap\B(b,r)
 \Big),Z\right) 
 \end{equation}
 which admits $\mathcal{A}$-limit $f$ on $\X\Big(A_0\cap\B(a,r),B_0\cap \B(b,r);D\cap\B(a,r),G\cap\B(b,r)\Big).$
 
  Fix  $0<\delta_0<\frac{1}{2}.$ Then it follows from  
      (\ref{eq1_thm_local}) 
       that for $ 0<r\leq r_0,$ $a\in\overline{A}_0,$  $ b\in\overline{B}_0.$
    \begin{equation*}
 D_{a,\delta_0,r} \times G_{b,\delta_0,r}  \subset
 \widehat{\X} \Big(A_0\cap\B(a,r),B_0\cap \B(b,r);D\cap\B(a,r)        ,G\cap\B(b,r)
 \Big) .
\end{equation*}
 This, combined with (\ref{eq_thm_local_hatf_ab}), implies that
  \begin{equation}\label{eq_thm_local_hatf_ab_new}
 \hat{f}_{a,b}\in \mathcal{O}\left(D_{a,\delta_0,r} \times G_{b,\delta_0,r}         ,Z\right),
 \qquad  0<r\leq r_0,\ a\in\overline{A}_0,\ b\in\overline{B}_0.
 \end{equation}
 Next we fix a finite covering $(A_0\cap\B(a_m,r))_{m=1}^M$ of $A_0$
 and  $(B_0\cap\B(b_n,r))_{n=1}^N$ of $B_0,$ where $(a_m)_{m=1}^M\subset \overline{A}_0$ and
   $(b_n)_{n=1}^N\subset \overline{B}_0.$ 

 We divide the proof into two steps.
 
 \smallskip
  
\noindent {\bf  Step 1:} {\it Fix an open set $G^{'}\Subset G.$
Then there exists $r_1$: $0<r_1<r_0$  with the following property: for every 
$a\in \overline{A}_0$ 
there exist
an open subset
$A_a$  of $D$ 
 and 
 a mapping 
 \begin{equation*}
 \hat{f}=\hat{f}_{a}\in  \mathcal{O}\left(A_a \times\left(G^{'}\cup\bigcup\limits_{n=1}^{N}  G_{b_n,\delta_0,r_0}\right)
          ,Z\right)
\end{equation*}
such that  
\begin{equation*}
\hat{f}(z,w)=\hat{f}_{a,b_n}(z,w),\qquad (z,w)\in (A_a\cap D_{a,\delta_0,r_0})\times  G_{b_n,\delta_0,r_0},
\qquad n=1,\ldots,N;
\end{equation*}
and that $A_a$ is of the form $\{z\in D\cap \B(a,r_1):\  \omega(z,A_0\cap \B(a,r_1),D\cap \B(a,r_1))<\delta_a \} $
for some  $0<\delta_a<\delta_0.$
 }
 
 \smallskip
 
 \noindent {\it Proof of Step 1.} Fix  an  arbitrary point $a_0\in\overline{A}_0.$
  First we claim that there are a sufficiently small number $r_1:$ $0<r_1<r_0$  and 
 a finite number of open subsets   $(V_n)_{n=1}^{N_0}$ of $G$  with the following properties:
 \begin{itemize}
 \item[(a)] 
 $V_1=  G_{b_1,\delta_0,r_0}$ and     $\left( G_{b_n,\delta_0,r_0}\right)_{n=1}^N\subset (V_n)_{n=1}^{N_0}$
 (see the notation in (\ref{eq1_thm_local}));
\item[(b)] $f|_{(\overline{A}_0\cap \B(a,r_1))\times V_n}$ is bounded, $n=1,\ldots,N_0;$ 
 \item[(c)]  $G^{'}\Subset \bigcup\limits_{n=1}^{N_0} V_n;$
\item[(d)] $V_n\cap V_{n+1}\not=\varnothing,$  $n=1,\ldots,N_0-1.$
 \end{itemize}
 
 Indeed,  we first start with the test $r_1:=r_0$ and $N_0:=N$ and 
  $ (V_n)_{n=1}^{N_0}:=\left( G_{b_n,\delta_0}\right)_{n=1}^N.$ In virtue of 
  (\ref{eq_thm_local_fab}) we see that 
 our choice 
 satisfies (a)--(b). If (c)--(d) are  satisfied then  we are done.
 Otherwise, we will make the following  procedure.
 
 Fix a point $w_0\in G^{'}.$
For $n=1,\ldots,N,$ let $\gamma_n:\ [0,1]\rightarrow G$ be
a continuous one-to-one map such that
\begin{equation*}
\gamma_n(0)=w_0\qquad\text{and}\qquad \gamma_n(1)\in  G_{b_n,\delta_0,r_0}.
\end{equation*} 
Since $f$ is locally bounded,
 there exist  
  sufficiently small numbers $r_1, s:$ $0<r_1\leq r_0$ and $0<s$ such that
  $f|_{(\overline{A}_0\cap \B(a,r_1))\times \B(w,s)}$ is bounded for all $a\in \overline{A}_0$ and $w\in
   \overline{G^{'}}\bigcup\limits_{n=1}^N \gamma_n([0,1]).$
 Therefore, we may add  to  the starting  collection  $ (V_n)_{n=1}^{N}$
     some balls of the form $ \B(w,s),$ where  $w\in
   \overline{G^{'}}\bigcup\limits_{n=1}^N \gamma_n([0,1]),$ and the new collection  $ (V_n)_{n=1}^{N_0}$ still satisfies (a)--(b).
   Now it remains to show that by adding a finite  number of   suitable balls   $ \B(w,s),$  (c)--(d) are also satisfied. But this
   assertion follows from an almost obvious
 geometric argument. In fact, we may  
   renumber  the collection  $(V_n)$  if necessary.
 Hence, the above claim has been shown.

 Using (c)--(d)  above we may  fix  open sets $U_n\Subset V_n$  for  $n=1,\ldots,N_0,$ such that
 \begin{equation}\label{eq_thm_local_choice_Un}
 G^{'}\Subset \bigcup\limits_{n=1}^{N_0} U_n\quad\text{and}\quad  U_n\cap U_{n-1}\not=\varnothing,\quad 1<n\leq N_0.
 \end{equation}
 In what follows we will find the desired set $A_{a_0}$ and the desired holomorphic mapping  $\hat{f}$ after $N_0$ steps.
  Namely, after the $n$-th step $(1\leq n\leq N_0)$, we construct  an open subset $A_{n}$ of $D$ in the form
  $D_{a_0,\delta_n,r_1}$ for a suitable  $\delta_n>0,$
  and 
  a    mapping  $\hat{f}_n\in \mathcal{O}\Big( A_{n}\times\big(\bigcup\limits_{p=1}^{n} U_{p}\big),Z\Big).$
  Finally, we obtain $A_{a_0}:=A_{N_0}$ and $\hat{f}:=\hat{f}_{N_0}.$ Now we carry out this construction.

  In the first step, using (\ref{eq1_thm_local}), (\ref{eq_thm_local_hatf_ab}), (\ref{eq_thm_local_hatf_ab_new})
   and  (a), we define 
  \begin{equation*}
 \delta_1:=\delta_0,\ A_1:=D_{a_0,\delta_1,r_1}\qquad\text{and}\qquad \hat{f}_1(z,w):=\hat{f}_{a_0,b_1}(z,w),\qquad  (z,w)\in A_{1}\times U_1.
 \end{equation*}
Suppose  that we have constructed an open subset $A_{n-1}$ of $D$ and 
a mapping  $\hat{f}_{n-1}\in  \mathcal{O}\Big(A_{n-1} \times \big(\bigcup\limits_{p=1}^{n-1} U_{p}\big)
,Z\Big)$ for some $n:\ 2\leq n\leq N_0.$
We wish to construct   an open subset  $A_{n}$ of $D$  and 
a mapping  $\hat{f}_{n}\in  \mathcal{O}\Big(A_n \times  \big(\bigcup\limits_{p=1}^{n} U_{p}\big)       ,Z\Big).$
There are two cases to consider.

\noindent {\bf Case} $V_n=G_{b_m,\delta_0}$ for some $1\leq m\leq N.$

In this case let $\delta_n:=\delta_{n-1}$ and  $A_{n}:=A_{n-1}=D_{a_0,\delta_{n-1},r_1} ,$ and
 \begin{equation*}
 \hat{f}_{n}:=
\begin{cases}
 \hat{f}_{n-1},
  & \qquad\text{on}\  A_{n}\times\Big(\bigcup\limits_{l=1}^{n-1} U_{l}\Big), \\
  \hat{f}_{a_0,b_m}, &   \qquad\text{on}\ A_{n}\times U_n        .
\end{cases}
\end{equation*}

\noindent {\bf Case} $V_n\not\in \Big(G_{b_m,\delta_0}\Big)_{m=1}^N.$ 

By (\ref{eq_thm_local_choice_Un}) fix a nonempty open set $K\Subset U_n\cap U_{n-1}.$ Then by the induction,
 $\hat{f}_{n-1}\in  \mathcal{O}\left(A_{n-1} \times K,Z\right).$
 Recall from (b) that $f:\ (\overline{A}_0\cap \B(a_0,r_1))\times V_n\longrightarrow Z$ is bounded.
 Since  $f$ is  locally bounded, by decreasing $r_1>0$ (if necessary)  we may  assume that
 \begin{equation*}
 g:=f|_{\X(\overline{A}_0\cap \B(a_0,r_1), K; D\cap \B(a_0,r_1) , V_n)}
 \end{equation*}
  is bounded.
 Applying  Theorem \ref{thm_local_localized} to $g,$   we obtain
  \begin{equation*}
  \hat{g}\in   \mathcal{O}\Big(\widehat{\X}(A_0\cap \B(a_0,r_1), K; D\cap\B(a_0,r_1) , V_n),Z\Big)
  \end{equation*}
  which extends $g.$ 
  Since $U_n\Subset V_n,$ we may choose $\delta_n$ such that
 $ 0<\delta_n<1-\sup\limits_{w\in U_n} \omega(w, K,V_n).$
Using this and  (\ref{eq1_thm_local}), 
it follows that
\begin{equation*}
D_{a_0,\delta_n,r_1}\times U_n\subset \widehat{\X}( A_0\cap \B(a_0,r_1), K; D  \cap \B(a_0,r_1)        , V_n              ).
\end{equation*}
Therefore, let  $A_n:=  D_{a_0,\delta_n,r_1}  $ and define
 \begin{equation*}
 \hat{f}_{n}:=
\begin{cases}
 \hat{f}_{n-1},
  & \qquad\text{on}\  A_n\times    \Big(\bigcup\limits_{l=1}^{n-1} U_{l}\Big) , \\
  \hat{g}, &   \qquad\text{on}\ A_{n}\times U_n        .
\end{cases}
\end{equation*}

This completes our construction in the $n$-step. 
Finally, we put $A_{a_0}:=A_{N_0}$ and $\hat{f}_{a_0}:=\hat{f}_{N_0}.$
Using this and (\ref{eq_thm_local_hatf_ab}) and (\ref{eq_thm_local_choice_Un}) and (a),
the desired conclusion  of Step 1 follows.  \hfill $\square$

 \smallskip
  
\noindent {\bf  Step 2:} {\it Completion of the proof.} 
  
\smallskip

\noindent{\it  Proof of Step 2.}
Fix a sequence of relatively compact open subsets $(D^{'}_k)_{k=1}^{\infty}$  of $D$
(resp.  $(G^{'}_k)_{k=1}^{\infty}$  of $G$)  such that
$D^{'}_k\nearrow D$ and  $G^{'}_k\nearrow G$  as $k\nearrow\infty.$ 
Put
 \begin{equation}\label{eq_thm_local_Dk_Gk}
 D_k:=D^{'}_k\cup\bigcup\limits_{m=1}^{M}  D_{a_m,\delta_0,r_0},\quad G_k:=G^{'}_k\cup\bigcup\limits_{n=1}^{N}  G_{b_n,\delta_0,r_0},\qquad k\geq 1.
\end{equation}
Using the result of Step 1,  
we may find,  for every $k,$ a number $0<r_{k} < r_0$  with the following properties:
\begin{itemize}
\item[$\bullet$] for every $a\in A_0,$ 
there is  $0<\delta_{a,k}<\delta_0$  
such that  by considering the open set
 \begin{equation*}
 A_{a,k}:=\{z\in D\cap \B(a,r_k):\  \omega\left(z,A_0\cap \B(a,r_k),D\cap \B(a,r_k)\right)<\delta_{a,k} \} 
 \end{equation*}
 one can find  
 a mapping 
 $
 \hat{f}_{a,k}\in  \mathcal{O}\left(A_{a,k} \times G_k,Z\right)
$
satisfying  
\begin{equation}\label{eq_thm_local_property_fak}
\hat{f}_{a,k}=\hat{f}_{a,b_n}\qquad\text{on}\ (A_{a,k}\cap D_{a,\delta_0,r_k})\times  G_{b_n,\delta_0,r_k},\qquad n=1,\ldots,N;
\end{equation}
 \item[$\bullet$]   for every $b\in B,$ 
there is  $0<\delta_{b,k}<\delta_0$ such that  by considering the open set
\begin{equation*}
B_{b,k}:=\{w\in G\cap \B(b,r_k):\  \omega\left(z,B_0\cap \B(b,r_k),G\cap \B(b,r_k)\right)<\delta_{b,k} \}    
\end{equation*}
one can find a mapping 
 $
 \hat{f}_{b,k}\in  \mathcal{O}\left(D_k\times B_{b,k} ,Z\right)
$
satisfying
\begin{equation}\label{eq_thm_local_property_fbk}
\hat{f}_{b,k}=\hat{f}_{a_m,b}\qquad\text{on}\  D_{a_m,\delta_0,r_k}
\times (B_{b,k}\cap G_{b,\delta_0,r_k}),\qquad m=1,\ldots,M.
\end{equation}
\end{itemize}

Next using the compactness of $\overline{A}_0$ and $\overline{B}_0,$  one may find, for every $k,$ 
two finite coverings $(A_0\cap\B(a_m^{'},r_{k}))_{m^{'}=1}^{M_k}$ of $A_0$
 and  $(B_0\cap\B(b_{n^{'}},r_{k}))_{n^{'}=1}^{N_k}$ of $B_0,$ where $(a_{m^{'}})_{m^{'}=1}^{M_k}\subset \overline{A}_0$ and
   $(b_{n^{'}})_{n^{'}=1}^{N_k}\subset \overline{B}_0.$ Put 
\begin{equation}\label{eq_thm_local_Ak_Bk}
A_{k}:=\bigcup\limits_{m^{'}=1}^{M_k} A_{a_{m^{'}},k}\quad\text{and}\quad B_{k}:=\bigcup\limits_{n^{'}=1}^{N_k}
 B_{b_{n^{'}},k},\quad k\geq 1 .
\end{equation}
 In virtue of (\ref{eq_thm_local_Dk_Gk})--(\ref{eq_thm_local_Ak_Bk})
and (\ref{eq_thm_local_fab})--(\ref{eq_thm_local_hatf_ab_new}), the family
$(\hat{f}_{a_{m^{'}},k})_{m^{'}=1}^{M_k} \bigcup (\hat{f}_{b_{n^{'},k}})_{n^{'}=1}^{N_k}$ is collective for every $k\geq 1.$
Let 
\begin{equation}\label{eq_thm_local_hatfk}
\hat{f}_k\in\mathcal{O}\Big(\X(A_k,B_k;D_k,G_k),Z\Big)
\end{equation}
 denote the collected mapping of this  family.
 
 Next, we show that
 \begin{equation}\label{eq_thm_local_limit}
 \lim\limits_{k\to\infty}\omega(z,A_0,D_k)=\omega(z,A_0,D)\ \text{and}\
  \lim\limits_{k\to\infty}\omega(w,B_0,G_k)=\omega(z,B_0,G),\qquad  z\in D,\  w\in G. 
 \end{equation}
It is sufficient to prove the first identity in (\ref{eq_thm_local_limit}) since the proof of the
second one is similar.
Observe that there is  $u\in\mathcal{PSH}(D)$ such that
$\omega(\cdot,A_0,D_k)\searrow u$ as $k\nearrow\infty$
and $u\geq \omega(\cdot,A_0,D)$ on $D.$
Therefore, the proof of (\ref{eq_thm_local_limit}) will be complete if one  can show that
 $u\leq \omega(\cdot,A_0,D)$  on $D.$
 
 To this end observe that for every $a\in A_0$  there is $1\leq m \leq M$ such that
 $a\in \B(a_m,r_0).$ Consequently, using (\ref{eq_thm_local_Dk_Gk}),
 \begin{equation*}
 (\Alimsup u)(a)\leq \Big(\Alimsup \omega(\cdot,A_0\cap \B(a_m,r_0), D_{a_m,\delta_0,r_0})\Big)(a)
 =0,
 \end{equation*}
 where the equality follows from an application of Proposition \ref{prop_levelsets}.
 This, combined with the obvious inequality $u\leq 1,$  implies that $u\leq  \omega(\cdot,A_0,D).$  Hence,
   (\ref{eq_thm_local_limit}) follows.
 
 We are now in the  position to define the desired extension  mapping $\hat{f}.$
Indeed, one  glues
$\left(\hat{f}_{k}\right)_{k=1}^{\infty}$ given in (\ref{eq_thm_local_hatfk}) together to obtain
$\hat{f}$ in the following way
\begin{equation*}
\hat{f}:=\lim\limits_{k\to \infty} \hat{f}_{k}\qquad \text{on}\
 \widehat{W}_0.
\end{equation*}
One needs to check that the last limit 
 exists and possesses all the required
properties.  In virtue of (\ref{eq_thm_local_property_fak})--(\ref{eq_thm_local_limit}), and the Gluing Lemma 
\ref{gluing_lem},
the proof will be complete if we can show the following\\
{\bf Claim.} {\it For every
$(z_0,w_0)\in\widehat{W}_0,$ there are an open neighborhood
$U\times V$ of  $(z_0,w_0)$ and   $\delta_0>0$ such that the hypotheses of Lemma \ref{gluing_lem}
is fulfilled with 
\begin{equation*}
\mathcal{D}:=D,\ \mathcal{G}:=G,\ \mathcal{P}_{k}:=A_{k},\  \mathcal{Q}_{k}:=B_k,\
\mathcal{D}_{k}:=D_k,\ \mathcal{G}_{k}:=G_k,\qquad k\geq 1.
\end{equation*}}
  
 To this end let  
 \begin{equation*} 
 \delta_{0}:=\frac{1-\omega(z_0,A_0,D)- \omega(w_0,B_0,G)}{2},
\end{equation*}
 and let $U\times V$ be an open  neighborhood of  $(z_0,w_0)$ such that
 \begin{equation*} 
   \omega(z,A_0,D)+ \omega(w,B_0,G)< \omega(z_0,A_0,D)+ \omega(w_0,B_0,G)+\delta_0.
\end{equation*}
Then   using  these inequalities and (\ref{eq_thm_local_limit}),
 we see that there is a sufficiently big $q_0\in \N$ such that
 for $q_0\leq q \leq p$  and  $(z,w)\in U\times V,$  
  \begin{equation*}
\begin{split}
 \omega(z,A_{p},D_q)+ \omega(w,B_p,D_q)&\leq
 \omega(z,A_0,D_q)+\omega(w,B_0,G_q)\\
&\leq\omega(z,A_0,D)+ \omega(w,B_0,G) +\delta_0
<1.
\end{split}
\end{equation*}
This proves the above claim. Hence, the proof of the theorem is finished.
 \end{proof}
 
Now we are  able to  formulate the following semi-local result. 
 \begin{thm}  \label{thm_semi_local}
  Let $D$  be an open subset of a complex manifold and $ G\subset \C^m$ a bounded open set and  $Z$  
 a complex analytic  space possessing the  Hartogs extension property.
  $D$  (resp.  $G$) is equipped with the canonical system of approach regions 
  (resp. the system of  approach regions   $\big(\mathcal{A}_{\beta}(\eta)\big)_{\eta\in\overline{ G},\  \alpha\in I_{\eta}}$).
 Let
$A$ be an open subset of $D$ and 
  let $B,$  $B_0$ be  subsets of  
  $\overline{ G}$ such that  $B_0$ is locally pluriregular
   and   $\overline{B}_0\subset   B^{\ast}.$
  Put 
    \begin{equation*}
    W:= \X(A,B;D,G)\quad\text{and}\quad
     W_0:=\X(A,B_0;D,G).
    \end{equation*}
   Then,
   for every  mapping  $f:\ W\longrightarrow Z$ which   satisfies the  following conditions:
\begin{itemize}
     \item[$\bullet$] $f\in\mathcal{C}_s(W,Z)\cap \mathcal{O}_s(W^{\text{o}},Z);$
    \item[$\bullet$] $f$  is  locally  bounded  along  $
     D\times (B\cap \partial G),$
\end{itemize}    
    there exists a unique  mapping 
$\hat{f}\in\mathcal{O}(\widehat{W}_0,Z)$ which
  admits  $\mathcal{A}$-limit $f(\zeta,\eta)$ at all points
  $(\zeta,\eta)\in    W_0.$  
   \end{thm}
\begin{proof}
First, applying Part 1) of Theorem   \ref{thm_local_classical} and using the  hypotheses on $f,$
we  see that $f$ extends  to a  mapping (still denoted  by) $f$ defined on  $\X(A,B\cup B^{\ast};D,G)$ 
such that $f$ is separately holomorphic   on $\X^{\text{o}}(A,B\cup B^{\ast};D,G)$
and  that $f|_{ \X(A, B^{\ast};D,G)}$ is  locally bounded.
 
 We define $\hat{f}$  at a point $(z_0,w_0)\in \widehat{W}_0$ as follows:
 Let $\epsilon>0$ be such that
 \begin{equation}\label{eq1_thm_semi_local}
 \omega(z_0,A,D)+\omega(w_0,B_0,G)+\epsilon<1.
 \end{equation}
   By Theorem \ref{Rosaythm} and Proposition \ref{prop_Rosay}, there is
   a holomorphic disc $\phi\in\mathcal{O}(\overline{E}, D)$ 
such that  $\phi(0)=z_0$   and  
\begin{equation}\label{eq2_thm_semi_local}
1-\frac{1}{2\pi}\cdot\mes (\phi^{-1}(A)\cap\partial E )<\omega(z_0,A,D)+\epsilon.
\end{equation} 
Moreover, using the hypotheses, we see that the mapping $f_{\phi},$  defined by
 \begin{equation}\label{eq3_thm_semi_local}
 f_{\phi}(t,w):=f(\phi(t),w),\qquad  (t,w)\in \X\left(\phi^{-1}(A)\cap\partial E,B
 ;E,G\right),
 \end{equation} 
satisfies the hypotheses of Theorem \ref{thm_local}.
By this theorem, let $\hat{f}_{\phi}$ be the unique mapping in $ \widehat{\X}\left(\phi^{-1}(A)\cap\partial E
,B_0;E,G\right)$
such that
  \begin{equation}\label{eq4_thm_semi_local}
   (\Alim \hat{f}_{\phi})(t,w)=f_{\phi}(t,w),\qquad  (t,w)\in
  \X\left(\phi^{-1}(A)\cap\partial E, B_0;E,G\right).
  \end{equation}
  In virtue of (\ref{eq1_thm_semi_local})--(\ref{eq2_thm_semi_local}), $(0,w_0)\in 
   \widehat{\X}\left(\phi^{-1}(A)\cap \partial E,B_0;E,G\right).$
   Then  the value at $(z_0,w_0)$ of the desired extension mapping $\hat{f}$ is given by
\begin{equation*}  
 \hat{f}(z_0,w_0):=\hat{f}_{\phi}(0,w_0) .
\end{equation*}
Using this and  (\ref{eq3_thm_semi_local})--(\ref{eq4_thm_semi_local}), 
and arguing as in Part 2) of Lemma \ref{lem_Part1}, 
one can show that $\hat{f}$  is well-defined on $\widehat{W}_0.$

To show that  $\hat{f}$ is holomorphic, one argues as in Step 1 of the proof of   Theorem \ref{thm_local_euclidean}.
To show that $\hat{f}$ 
  admits  $\mathcal{A}$-limit $f(\zeta,\eta)$ at all points
  $(\zeta,\eta)\in    W_0$ and that it is uniquely defined,
 one proceeds as in Step 2--4 of the proof of Theorem \ref{thm_local_euclidean}
 making the obviously necessary changes and adaptations.
  Hence, the proof  is finished.
\end{proof}


%
%
%
%
%

\section{The proof of  Theorem A}
First we need a variant of Definition \ref{defi_pluri_measure}.
For a set $A\subset \overline{D},$
Let $\widetilde{\mathcal{E}}(A)$ be  the set of all elements    $P\in \mathcal{E}(A)$
with the property  that there is an open neighborhood $U\subset X$ of $P$ such that
$U$ is biholomorphic to a domain in some $\C^n.$  
Then it can be checked that 
\begin{equation}\label{eq1_Proof_Theorem_A}
\widetilde{A}:=\bigcup\limits_{P\in \widetilde{\mathcal{E}}(A)} P.
\end{equation}
This identity will allow us to pass from ``local informations" to ``global extensions".

 For the proof we need to develop some preparatory results.
 
In virtue of (\ref{eq1_Proof_Theorem_A}), for any $P\in  \widetilde{\mathcal{E}}(A)          $ (resp. $Q\in 
\widetilde{\mathcal{E}}(B)$) fix an
open neighborhood $U_{P}$ of  $P$  (resp.  $V_Q$ of $Q$) such that $U_{P} $
(resp. $V_Q$) is  biholomorphic to a domain in $\C^{d_P}$  (resp. in
$\C^{d_Q}$),
where $d_P $ (resp.   $d_Q$) is the dimension of $D$ (resp. $G$) at points of  $P$ (resp. $Q$).
 For any $0<\delta\leq\frac{1}{2}$  define
\begin{equation}\label{eq_UV_Proof_Theorem_A}
\begin{split}
U_{P,\delta}&:=\left\lbrace z\in U_{P}:\ \omega(z, P,   U_P)<\delta  \right\rbrace,\qquad
P\in           \widetilde{\mathcal{E}}(A) ,\\
V_{Q,\delta}&:=\left\lbrace w\in V_{Q}:\ \omega(w, Q,   V_Q)<\delta  \right\rbrace,\qquad
Q\in      \widetilde{\mathcal{E}}(B) ,\\
A_{\delta}&:=\bigcup\limits_{P\in   \widetilde{\mathcal{E}}(A)} U_{P,\delta},\qquad
B_{\delta}:=\bigcup\limits_{Q\in   \widetilde{\mathcal{E}}(B)} V_{Q,\delta},\\
 D_{\delta}&:=\left\lbrace z\in D:\
  \widetilde{\omega}(z,A,D)<1-\delta\right\rbrace,\quad
  G_{\delta}:=\left\lbrace w\in G:\  \widetilde{\omega}(w,B,G)<1-\delta\right\rbrace.
\end{split}
\end{equation}

\begin{lem}\label{lem_comparison}
 We keep the above notation.
Then:
\begin{itemize}
\item [(1)] For every $\zeta\in \widetilde{A}$ and $\alpha\in I_{\zeta},$
there is an open neighborhood $\mathcal{U}$ of $\zeta$ such that
$\mathcal{U}\cap \mathcal{A}_{\alpha}(\zeta)\subset A_{\delta}.$
\item[(2)] $A_{\delta}$ is an open subset of $D$ and  $ A_{\delta}\subset D_{1-\delta} \subset D_{\delta}.$
\item[(3)] $ \widetilde{\omega}(z,A,D)-\delta
\leq  \omega(z, A_{\delta},D) \leq\widetilde{\omega}(z,A,D),$ $z\in D .$
\end{itemize}
\end{lem}

\smallskip

\noindent {\it Proof of Lemma \ref{lem_comparison}.} 
To prove Part (1) fix, in view of   (\ref{eq1_Proof_Theorem_A})--(\ref{eq_UV_Proof_Theorem_A}),
  $P\in \widetilde{\mathcal{E}}(A),$   $\zeta\in  P$ and $\alpha\in I_{\zeta}.$
Using   the definition
of local pluriregularity, we see that
$\limsup\limits_{z\to\zeta,\ z\in\mathcal{A}_{\alpha}(\zeta)}\omega(z,P,U_P)=0.$
Hence, Part (1) follows.

The assertion that $A_{\delta}$ is open follows immediately from (\ref{eq_UV_Proof_Theorem_A}).
 Since $0<\delta\leq\frac{1}{2},$ the second  inclusion in Part (2)
 is clear. To prove the first
inclusion let $z$ be an arbitrary point of
$A_{\delta}.$ Then there is   $P\in   \widetilde{\mathcal{E}}(A)$   such that $z\in
U_{P,\delta}.$ Using (\ref{eq_UV_Proof_Theorem_A})  and Definition \ref{defi_pluri_measure}  we
 obtain
\begin{equation}\label{eq1_lem_comparison}
 \widetilde{\omega}(z,A,D)     =\omega(z,\widetilde{A},D)
\leq  \omega(z,P,U_P)<\delta .
\end{equation}
Hence, $z\in D_{1-\delta},$ which in turn implies that $A_{\delta}\subset
D_{1-\delta}.$

It follows from Part (1) that
 \begin{equation*}
 \omega(z,A_{\delta},D)
\leq   \omega(z,\widetilde{A},D) = \widetilde{\omega}(z,A,D),\qquad z\in D,
\end{equation*}
which proves the second estimate in  Part (3).
To complete the proof  let   $P\in   \widetilde{\mathcal{E}}(A) $  and
$0<\delta\leq\frac{1}{2}.$ We deduce from (\ref{eq1_lem_comparison})  that $
 \widetilde{\omega}(z,A ,D)-\delta\leq 0$ for  $   z\in U_{P,\delta}.$ Hence, by
(\ref{eq_UV_Proof_Theorem_A}),
\begin{equation*}
 \widetilde{\omega}(z,A ,D)-\delta\leq 0, \qquad z\in
A_{\delta}.
\end{equation*}
On the other hand,  $ \widetilde{\omega}(z,A ,D)-\delta<1,$ $z\in D.$ Recall from Part (2) that
$A_{\delta}$ is an open subset of $D_{\delta}.$
Consequently, the first estimate of Part (3)
follows.  
\hfill $\square$

Now we are able to  to prove Theorem  A  in the following
special case.
\begin{prop}\label{prop_special_case}
 Let $D$ be an open subset of a complex manifold and  $ G $ a  bounded open subset of $\C^m$ and  $Z$  
 a complex analytic  space possessing the  Hartogs extension property.
  $D$  (resp.  $G$) is equipped with a 
  system of approach regions 
  $\big(\mathcal{A}_{\alpha}(\zeta)\big)_{\zeta\in\overline{D},\  \alpha\in I_{\zeta}}$
  (resp.  $\big(\mathcal{A}_{\beta}(\eta)\big)_{\eta\in\overline{G},\  \beta\in I_{\eta}}$).    Let
  $A$ be a  subset of $\overline{D},$  
  let $B,$  $B_0$ be  subsets of  
  $\overline{ G}$ such that  $B_0$ is locally pluriregular
   and   $\overline{B}_0\subset   B^{\ast}.$
  Put 
    \begin{eqnarray*}
    W&:=& \X(A,B;D,G),\ \
     W_0:=\X(A,B_0;D,G),\ \  
       \widetilde{W}^{\text{o}} :=\big ( (D\cup \widetilde{A})\times B_0\big)\bigcup \big( \widetilde{A}\times(G\cup B_0)\big),       \\
     \widehat{W}^{\text{o}}&:=&\left\lbrace (z,w)\in D\times G:\ \widetilde{\omega}(z,A,D)+\omega(w,B_0,G)<1\right\rbrace.
    \end{eqnarray*}
   Then,
   for every  mapping  $f:\ W\longrightarrow Z$ which   satisfies the  following conditions:
\begin{itemize}
     \item[$\bullet$] $f\in\mathcal{C}_s(W,Z)\cap \mathcal{O}_s(W^{\text{o}},Z);$
    \item[$\bullet$]  $f$ is locally bounded along $\X\big (A\cap\partial D,B\cap\partial G;D,G\big);$
   \item[$\bullet$]    $f|_{A\times B}$ is continuous at all points of 
   $(A\cap\partial D)\times (B\cap \partial G),$
    \end{itemize}    
  there exists a unique mapping 
$\hat{f}\in\mathcal{O}(\widehat{W}^{\text{o}},Z)$ which
  admits  $\mathcal{A}$-limit $f(\zeta,\eta)$ at all points
  $(\zeta,\eta)\in    \widetilde{W}^{\text{o}}.$     
\end{prop}

\smallskip

\noindent {\it Proof of  Proposition \ref{prop_special_case}.}
First, applying Part 1) of Theorem   \ref{thm_local_classical} and using the  hypotheses on $f,$
we  see that $f$ extends  to a  mapping (still  denoted  by  $f$) defined on  $\X(A\cup A^{\ast},B\cup B^{\ast};D,G)$ 
such that $f$ is  separately holomorphic   on $\X^{\text{o}}(A\cup A^{\ast},B\cup B^{\ast};D,G)$ and
that    $f|_{\X( A^{\ast},B^{\ast};D,G)}$ is locally bounded.
 
For each $P\in  \widetilde{\mathcal{E}}(A)    ,$  $U_P$ (resp. $G$) is biholomorphic to an open set in $\C^{d_P}$ (resp. in $\C^m$).
 Consequently,  the  mapping $    f_P:=f|_{\X\left( \overline{P} ,B;  U_{P},G\right)}$
   satisfies the hypotheses of   Theorem
\ref{thm_local}. Hence, we obtain a unique mapping
   $\hat{f}_{P} \in   \mathcal{O}\big(\widehat{\X}
 \left(P ,B_0; U_{P},G\right),Z\big) $
such that
\begin{equation}\label{eq1_prop_special_case}
(\Alim\hat{f}_{P})(z,w)=f_P(z,w)=f(z,w),\qquad  (z,w)\in
 \X\left( P,B_0;  U_{P},G\right) .
\end{equation}
Let $0<\delta\leq \frac{1}{2}$ 
 and $G^{'}_{\delta}:=\{w\in G:\  \omega(w,B_0,G)<1-\delta\}.$
 We will show that  the family 
$\left(\hat{f}_{P}|_{U_{P,\delta}\times G^{'}_{\delta}} \right)_{P\in  \widetilde{\mathcal{E}}(A)}$ is collective
in the sense of Definition \ref{defi_collection},
where  $U_{P,\delta}$  is given  in   (\ref{eq_UV_Proof_Theorem_A}).

 To prove this assertion let $P_1,\ P_2 $ be arbitrary  elements of
 $ \widetilde{\mathcal{E}}(A)    .$
  By (\ref{eq1_prop_special_case}), we have    
\begin{equation} \label{eq2_collecting_lem}
(\Alim \hat{f}_{P_1}) (z,w)=f(z,w)=(\Alim \hat{f}_{P_2})(z,w),\qquad (z,w)\in (U_{P_1}\cap U_{P_2})\times
 B_0.
 \end{equation}
The  assertion  is  reduced to showing  that 
    \begin{equation}\label{eq3_collecting_lem}
\hat{f}_{P_1} (z,w)=\hat{f}_{P_2}(z,w),\qquad (z,w)\in\widehat{\X} \left(P_1 ,B_0; U_{P_1},G\right)
\cap\widehat{\X}\left( P_{2} ,B_0; U_{P_2},G\right).
 \end{equation}
 To this end fix $(z_0,w_0) \in \widehat{\X} \left(P_1 ,B_0; U_{P_1},G\right)
\cap\widehat{\X}\left( P_{2} ,B_0; U_{P_2},G\right).$ 
Observe that both mappings $w\mapsto \hat{f}_{P_1}(z_0,w)$ and  $w\mapsto \hat{f}_{P_2}(z_0,w)$
belong to $\mathcal{O}(\mathcal{G},Z),$ where
$\mathcal{G}$ is  the connected component which contains $w_0$ of the following open set
\begin{equation*}
\left\lbrace  w\in G:\   \omega(w, B_0,G)<1-\max\limits_{j\in\{1,2\}}\omega(z_0,P_j,U_j) \right\rbrace.
\end{equation*}
  Applying Theorem \ref{unique2} to these mappings using   (\ref{eq2_collecting_lem}), Proposition \ref{prop_levelsets}
  and
  (\ref{eq3_collecting_lem}),  the above assertion follows.

In virtue of (\ref{eq_UV_Proof_Theorem_A}) let     
\begin{equation}\label{eq2_prop_special_case}
\tilde{\tilde{f}}_{\delta}\in \mathcal{O}(A_{\delta}\times
G^{'}_{\delta},Z)
\end{equation}
denote the collected mapping of the family $\left(\hat{f}_{P}|_{U_{P,\delta}\times G^{'}_{\delta}} \right)_{P\in  \widetilde{\mathcal{E}}(A)}$.
 In virtue of  (\ref{eq1_prop_special_case}) and (\ref{eq2_prop_special_case}),
 we are able to define a new   mapping  $\tilde{f}_{\delta}$ on $\X\left(A_{\delta},B;D,G^{'}_{\delta}
 \right)$ as follows
\begin{equation*}
 \tilde{f}_{\delta}:=
\begin{cases}
 \tilde{\tilde{f}}_{\delta},
  & \qquad\text{on}\  A_{\delta}\times G^{'}_{\delta}, \\
  f, &   \qquad\text{on}\ D\times B        .
\end{cases}
\end{equation*}
Using this  and  (\ref{eq1_prop_special_case})--(\ref{eq2_prop_special_case}), we  see that
\begin{equation}\label{eq3_prop_special_case}
 \Alim \tilde{f}_{\delta}=f\qquad\text{on}\ \X(A\cap \widetilde{A},B_0;D,G^{'}_{\delta}).
\end{equation}
 Since $A_{\delta}$ is an  open subset of $X$ and $G^{'}_{\delta}$ is a  bounded open set in $\C^m,$
   we are able to  apply  Theorem  \ref{thm_semi_local}
to  $ \tilde{f}_{\delta}$ in order to obtain a mapping
  $\hat{f}_{\delta}\in \mathcal{O}\big(  \widehat{\X}\left(A_{\delta}, B_0;D,G^{'}_{\delta}
 \right),Z \big)$  such that
\begin{equation}\label{eq4_prop_special_case}
\Alim  \hat{f}_{\delta}=  \tilde{f}_{\delta}\qquad\text{on}\ \X(A_{\delta},
  B_0;D,G^{'}_{\delta}).
\end{equation}

We are now in a position to define the desired extension  mapping $\hat{f}.$
Indeed, one  glues
$\left(\hat{f}_{\delta}\right)_{0<\delta\leq\frac{1}{2}}$ together to obtain
$\hat{f}$ in the following way
\begin{equation*}
\hat{f}:=\lim\limits_{k\to \infty} \hat{f}_{\frac{1}{k}}\qquad \text{on}\
 \widehat{W}^{\text{o}}  .
\end{equation*}
One needs to check that the last limit
 exists and possesses all the required
properties.  In virtue of (\ref{eq3_prop_special_case})--(\ref{eq4_prop_special_case})
and Lemma \ref{gluing_lem},
the proof will be complete if one can show that for every
$(z_0,w_0)\in\widehat{W}^{\text{o}},$ there are an open neighborhood
$U\times V$ of  $(z_0,w_0)$ and   $\delta_0>0$ such that  hypothesis (iii) of Lemma \ref{gluing_lem}
is fulfilled with 
\begin{equation*}
\mathcal{D}:=D,\ \mathcal{G}:=G,\ \mathcal{P}_{k}:=A_{\frac{1}{k}},\  \mathcal{Q}_{k}:=B_0,\
\mathcal{D}_{k}:=D,\ \mathcal{G}_{k}:=G^{'}_{\frac{1}{k}},\qquad k>2.
\end{equation*}
  To this end let  
 \begin{equation*}
 \delta_{0}:=\frac{1-\widetilde{\omega}(z_0,A,D)- \omega(w_0,B_0,G)}{2},
\end{equation*}
 and let $U\times V$ be an open  neighborhood of  $(z_0,w_0)$ such that
 \begin{equation*} 
   \widetilde{\omega}(z,A,D)+ \omega(w,B_0,G)< \widetilde{\omega}(z_0,A,D)+ \omega(w_0,B_0,G)+\delta_0.
\end{equation*}
Then for $k> \frac{1}{\delta_0}$ and for $(z,w)\in U\times V,$    using 
the last inequality,  and applying   Part (3) of Lemma 
\ref{lem_comparison}   and  Proposition \ref{prop_levelsets}, we see that
\begin{equation*}
\begin{split}
 \widetilde{\omega}(z,A_{\frac{1}{k}},D)+ \omega(w,B_0,G^{'}_{\delta_0})&\leq
 \widetilde{\omega}(z,A,D)+\frac{ \omega(w,B_0,G)}{1-\delta_0}\\
&\leq\frac{ \widetilde{\omega}(z,A,D)+ \omega(w,B_0,G) }{1-\delta_0}
<1.
\end{split}
\end{equation*}
This proves the above assertion. Hence, the proof of the proposition is finished.
\hfill  $\square$

\smallskip

We now arrive at

\smallskip

\noindent {\bf Proof of  Theorem A.} 
First, applying Part 1) of Theorem   \ref{thm_local_classical} and using the  hypotheses on $f,$
we  see that $f$ extends  to a  mapping (still  denoted  by) $f$  defined on  $\X(A\cup A^{\ast},B\cup B^{\ast};D,G)$ 
such that $f$ is  separately holomorphic   on $\X^{\text{o}}(A\cup A^{\ast},B\cup B^{\ast};D,G)$ and that
  $f|_{\X( A^{\ast}, B^{\ast};D,G)}$ is locally
bounded.

For each $P\in  \widetilde{\mathcal{E}}(A)    ,$  $U_P$  is biholomorphic to an open set in $\C^{d_P}.$ 
 Consequently,  the  mapping $ f_P:=f|_{\X\left( \overline{P} ,B;  U_{P},G\right)}$
   satisfies the hypotheses  of  Proposition \ref{prop_special_case}.  Hence, we obtain a unique mapping
     $\hat{f}_{P} \in   \mathcal{O}\big(\widehat{\X}^{\text{o}}
 \left(P ,B; U_{P},G\right),Z\big) $\footnote{
 Here  $\widehat{\X}^{\text{o}}
 \left(P ,B; U_{P},G\right):=\left\lbrace  (z,w)\in U_P\times G:\ \omega(z,P,U_P)+\widetilde{\omega}(w,B,G)
 <1 \right\rbrace.$}
such that
\begin{equation}\label{eq1_thmA_special_case}
(\Alim\hat{f}_{P})(z,w)=f(z,w),\qquad  (z,w)\in
 \X\big(  P,\widetilde{B}\cap B;  U_P,G\big) .
\end{equation}
 Let $0<\delta\leq\frac{1}{2}.$ Using  (\ref{eq1_thmA_special_case}) and  arguing as in  the proof of  Proposition \ref{prop_special_case},
 we may collect  the family
 $\left(\hat{f}_{P}|_{U_{P,\delta}\times G_{\delta}} \right)_{P\in  \widetilde{\mathcal{E}}(A)       }$
in order to obtain the collected   mapping $\tilde{f}^A_{\delta}\in \mathcal{O}(A_{\delta}\times
G_{\delta},Z).$

Similarly,  for each $Q\in   \widetilde{\mathcal{E}}(B) ,$   one obtains a unique mapping
 $\hat{f}_{Q} \in   \mathcal{O}\big(\widehat{\X}^{\text{o}}
 \left(A, Q ; D,V_{Q}\right),Z\big) $\footnote{
 Here  $\widehat{\X}^{\text{o}}
 \left(A ,Q; D,V_{Q}\right):=\left\lbrace  (z,w)\in D\times V_Q:\ \widetilde{\omega}(z,A,D)+\omega(w,Q,V_Q)
 <1 \right\rbrace.$}
such that
\begin{equation}\label{eq2_thmA_special_case}
(\Alim\hat{f}_{Q})(z,w) =f(z,w),\qquad (z,w)\in
 \X\big( A\cap  \widetilde{A},  Q;  D,V_{Q}\big) .
\end{equation}
Moreover, one can collect  the family
 $\left(\hat{f}_{Q}|_{D_{\delta}\times V_{Q,\delta}} \right)_{ Q\in   \widetilde{\mathcal{E}}(B)             }$
in order to obtain the collected   mapping $\tilde{f}^B_{\delta}\in \mathcal{O}(D_{\delta}\times
B_{\delta},Z).$

Next, we prove that
\begin{equation} \label{eq3_thmA_special_case}
 \tilde{f}^A_{\delta}=\tilde{f}^B_{\delta}\qquad \text{on}\  A_{\delta}\times
B_{\delta}.
\end{equation}
Indeed, in virtue of  (\ref{eq1_thmA_special_case})--(\ref{eq2_thmA_special_case})  it suffices to show that
 for any $ P\in   \widetilde{\mathcal{E}}(A)$ and $ Q\in   \widetilde{\mathcal{E}}(B)$ 
 and any $0<\delta\leq\frac{1}{2},$  
\begin{equation}\label{eq4_thmA_special_case}
\hat{f}_{P} (z,w)=\hat{f}_{Q}(z,w),\qquad (z,w)\in U_{P,\delta}\times V_{Q,\delta}.
 \end{equation}
 Observe that in virtue of  (\ref{eq1_thmA_special_case})--(\ref{eq2_thmA_special_case})  
 one has that
\begin{equation*} 
(\Alim \hat{f}_{P})(z,w)=(\Alim\hat{f}_{Q})(z,w)=f(z,w), \qquad (z,w)\in   \X\left( P 
, Q; U_{P},V_Q\right).
\end{equation*}
  Recall that $U_P$ (resp. $V_Q$) is biholomorphic to a domain
in $\C^{d_P}$ (resp. $\C^{d_Q}$). Consequently, applying the uniqueness of Theorem \ref{thm_local} yields that
 \begin{equation*}
\hat{f}_{P}(z,w)= \hat{f}_{Q}(z,w) , \qquad (z,w)\in \widehat{\X}\left(P ,Q;
U_{P},V_Q\right).
\end{equation*}
   Hence, the proof of  (\ref{eq4_thmA_special_case}) and then the proof of (\ref{eq3_thmA_special_case})          are finished.

 In virtue of  (\ref{eq3_thmA_special_case}),
 we are able to define a new  mapping $\tilde{f}_{\delta}:\
 \X^{\text{o}}\left(A_{\delta}, B_{\delta};D_{\delta},
 G_{\delta}\right)\longrightarrow Z$ as follows
\begin{equation}\label{eq5_thmA_special_case}  
 \tilde{f}_{\delta}:=
\begin{cases}
\tilde{f}^A_{\delta},
  & \qquad\text{on}\  A_{\delta}\times G_{\delta}, \\
  \tilde{f}^B_{\delta}, &   \qquad\text{on}\ D_{\delta}\times B_{\delta}        .
\end{cases}
\end{equation}
Using formula (\ref{eq5_thmA_special_case}) it can be readily checked that
$\tilde{f}_{\delta}\in \mathcal{O}_s\Big(\X^{\text{o}}\left(A_{\delta}, B_{\delta};D_{\delta},G_{\delta}
 \right),Z\Big).$ Since  we know from  Part (2) of Lemma \ref{lem_comparison}
  that  $A_{\delta}$
 (resp.  $B_{\delta}$) is an open subset of $D_{\delta}$  (resp.
 $G_{\delta}$), we are able to apply  Theorem  \ref{thm_vietanh}
to  $ \tilde{f}_{\delta}$  for every $0<\delta\leq\frac{1}{2}.$ Consequently,
one
  obtains a unique  mapping $\hat{f}_{\delta}\in \mathcal{O}\Big(
 \widehat{\X}\left(A_{\delta}, B_{\delta};D_{\delta},G_{\delta}
 \right),Z \Big)$ such that
\begin{equation}\label{eq6_thmA_special_case}  
  \hat{f}_{\delta}=  \tilde{f}_{\delta}\qquad\text{on}\ \X^{\text{o}}\left(A_{\delta},
  B_{\delta};D_{\delta},G_{\delta}\right).
\end{equation}
It follows from   (\ref{eq1_thmA_special_case})--(\ref{eq2_thmA_special_case}) and
   (\ref{eq5_thmA_special_case})--(\ref{eq6_thmA_special_case})  that
 \begin{equation}\label{eq7_thmA_special_case}
 \Alim \hat{f}_{\delta}=f\qquad\text{on}\
\X\left(A \cap \widetilde{A}, B\cap \widetilde{B};D_{\delta},G_{\delta}\right).
 \end{equation}
 In addition, for any $0<\delta\leq\delta_0\leq\frac{1}{2},$ and any $(z,w)\in
 A_{\delta}\times B_{\delta},$ there is   $P\in       \widetilde{\mathcal{E}}(A) 
          $ 
  such that $z\in U_{P,\delta_0}.$ 
  Therefore, it
 follows from the construction of $\tilde{f}^A_{\delta},$
    (\ref{eq5_thmA_special_case}) and (\ref{eq6_thmA_special_case}) that
\begin{equation*}
  \hat{f}_{\delta}(z,w)=  \hat{f}_P(z,w)=\hat{f}_{\delta_0}(z,w) .
\end{equation*}
This proves that $ \hat{f}_{\delta} =  \hat{f}_{\delta_0}$ on $
  A_{\delta}\times B_{\delta}$ for $ 0<\delta\leq\delta_0\leq\frac{1}{2}.$
  Hence,
\begin{equation}\label{eq8_case}
  \hat{f}_{\delta} =  \hat{f}_{\delta_0} \qquad\text{on}\
 \X( A_{\delta},B_{\delta};D_{\delta_0},G_{\delta_0}),\ 0<\delta\leq\delta_0\leq\frac{1}{2}
   .
\end{equation}

We are now in a position to define the desired extension  mapping $\hat{f}.$
 \begin{equation*}
\hat{f}:=\lim\limits_{k\to \infty} \hat{f}_{\frac{1}{k}}\qquad \text{on}\
 \widehat{\widetilde{W}}.
\end{equation*}
 To prove that $\hat{f}$  satisfies   the desired conclusion of the theorem one
  proceeds as in the end of the proof of Proposition \ref{prop_special_case}.   In virtue of
   (\ref{eq7_thmA_special_case})--(\ref{eq8_case}) and Lemma \ref{gluing_lem},
the proof will be complete if we can  verify that for every
$(z_0,w_0)\in\widehat{W},$ there are an open neighborhood
$U\times V$ of  $(z_0,w_0)$ and   $\delta_0>0$ such that  hypothesis (iii) of Lemma \ref{gluing_lem} is fulfilled with 
\begin{equation*}
\mathcal{D}:=D,\ \mathcal{G}:=G,\ \mathcal{P}_{k}:=A_{\frac{1}{k}},\  \mathcal{Q}_{k}:=B_{\frac{1}{k}},\
\mathcal{D}_{k}:=D_{\frac{1}{k}},\ \mathcal{G}_{k}:=G_{\frac{1}{k}},\qquad k>2.
\end{equation*}
  Since the verification follows along almost the same lines as 
  that of  Proposition \ref{prop_special_case},
  it is, therefore, left to the interested reader. 
 
Hence, the proof of Theorem A  is finished.
 \hfill $\square$

\section{Applications} \label{section_Applications}
In this  section we  give various  applications  of Theorem A using  different systems of approach regions
defined in Subsection  \ref{Subsection_Examples}. 
\subsection{Canonical  system of  approach  regions}

For every open subset $U\subset \R^{2n-1}$ and every continuous function $h:\  U\longrightarrow\R,$
the graph  
$$\left\lbrace z=(z^{'},z_n)=(z^{'},x_n+iy_n)\in \C^n:\ (z^{'},x_n)\in U\quad\text{and}\quad y_n=h(z^{'},x_n)     \right\rbrace$$
is called a {\it topological hypersurface in $\C^n$.} 

 Let  $X$ be a complex manifold of dimension $n.$
 A subset $A\subset X$ is said to be  a  {\it topological hypersurface} if,
 for every point $a\in A,$ there is a    local chart  $(U,\phi:\ U\rightarrow\C^n)$ around $a$
 such that  $\phi(A\cap U)$ is a topological hypersurface in $\C^n$

Now let    $D\subset X$  be an open subset and let  $A\subset\partial D$
be an open subset
(with respect to the topology induced on $\partial D$). Suppose in addition that
$A$ is a topological hypersurface.
A point  $a\in A$ is said to be {\it   of type 1  (with respect to $D$)}
if, for every neighborhood $U$ of $a$ there is an open  neighborhood $V$ of $a$ such that $V\subset U$ and  $V\cap D$ is  
a domain.  Otherwise, $a$ is said to be {\it   of type 2}.
We see easily that if $a$ is of type 2, then     for every neighborhood $U$ of $a,$   
   there are an open neighborhood $V$ of $a$ and
two  domains $V_1,$ $V_2$  such that $V\subset U,$    $V\cap D=V_1\cup V_2$ 
 and all points in $A\cap V$ are of type 1 with respect to $V_1$ and $V_2.$
 
 In virtue of  Proposition 3.7 in  \cite{pn3} we  have the following  
 \begin{prop}\label{prop_application2_1}
    Let  $X$ be a complex manifold and 
     $D$   an open subset of $X.$
  $D$ is  equipped  with  the canonical  system  of  approach  regions.
  Suppose that  $A\subset\partial D$  is an open boundary subset
which is also a topological hypersurface.
Then $A$ is locally  pluriregular and  $A\subset \widetilde{A}.$
\end{prop}

This, combined with Theorem A, implies the  following result.
\begin{thm}\label{application2_2}
Let $X,\ Y$  be two complex manifolds,
  and $D\subset X,$ $ G\subset Y$  two nonempty open sets. $D$ (resp.  $G$) is  equipped  with
  the canonical  system of  approach  regions.
   Let
  $A$ (resp. $B$) be a nonempty open subset of  $\partial D$ (resp.
  $\partial G$) which is also a
   topological hypersurface.
   Let $Z$  be a     complex  analytic space    possessing the  Hartogs extension property.
   Define 
   \begin{eqnarray*}
   W  &:= &\X(A,B;D,G),\\
    \widehat{W}  &:= &\left\lbrace  (z,w)\in D\times G:\  \omega(z,A,D)+\omega(w,B,G)<1     \right\rbrace.
   \end{eqnarray*}
 Let  $f:\ W\longrightarrow Z$ be   such that:
\begin{itemize}
\item[ (i)]  $f\in\mathcal{C}_s(W,Z)\cap \mathcal{O}_s(W^{\text{o}},Z);$
\item[(ii)]   $f$ is locally bounded on $W;$ 
\item[ (iii)]   $f|_{A\times B}$ is continuous. 
\end{itemize}

Then  there exists a unique mapping
$\hat{f}\in\mathcal{O}(\widehat{W})$
such that  
\begin{equation*}
\lim\limits_{\widehat{W}\ni(z,w)\to(\zeta,\eta)}\hat{f}(z,w)=f(\zeta,\eta),\qquad (\zeta,\eta)\in W.
\end{equation*} 

 If, moreover,  $Z=\C$ and $\vert f\vert_W<\infty,$ then
\begin{equation*}
 \vert \hat{f}(z,w)\vert\leq \vert f\vert_{A\times B}^{1-\omega(z,w)} \vert
 f\vert_W^{\omega(z,w)},\qquad (z,w)\in\widehat{W}.
\end{equation*}
\end{thm}

The special case where $Z=\C$ has been proved in  \cite{pn3}.
\subsection{System of angular approach  regions}

We will use  the  terminology and the   notation 
 in Paragraph  3 of Subsection \ref{Subsection_Examples}.
  More  precisely, if $D$ is  an open set of a  Riemann surface  such that  $D$ is  good on a  nonempty part of $\partial D,$ we equip $D$
    with the  system of angular  approach  regions 
 supported on this  part. Moreover, the     notions such as set of positive  length,
 set of zero length, locally pluriregular point 
 which exist on $\partial E$ can be  transferred  to
 $\partial D$   using   conformal mappings  in a local  way  (see \cite{pn2} for  more details). 

  \begin{thm} \label{application1_2}
  Let  $X,\ Y$  be  Riemann  surfaces  
  and $D\subset X,$ $ G\subset Y$ open subsets and
  $A$ (resp. $B$) a subset of   $\partial D$ (resp.
  $\partial G$)  such that
      $D$ (resp. $G$) is  good on $A$ (resp. $B$) and that both
  $A$ and $B$ are of   positive length.  Let $Z$  be a     complex  analytic space    possessing the  Hartogs extension property.
 Define  
 \begin{eqnarray*}
 W&:=&\X(A,B;D,G),\qquad W^{'}:=\X(A^{'},B^{'};D,G),\\
  \widehat{W}  &:= &\left\lbrace  (z,w)\in D\times G:\  \omega(z,A,D)+\omega(w,B,G)<1     \right\rbrace,\\
   \widehat{W^{'}}  &:= &\left\lbrace  (z,w)\in D\times G:\  \omega(z,A^{'},D)+\omega(w,B^{'},G)<1     \right\rbrace,
 \end{eqnarray*}
 where  $A^{'}$  (resp. $B^{'}$) is the set of   points  at  which  $A$ (resp. $B$) is locally pluriregular
 with respect to the system of angular  approach  regions  supported on $A$  (resp.  $B$), and  $\omega(\cdot,A,D),$
 $\omega(\cdot,A^{'},D)$   (resp.  $\omega(\cdot,B,G),$
 $\omega(\cdot,B^{'},G)$) are calculated using the canonical  system of approach regions. 
 
 Then for every mapping  $f:\ W\longrightarrow Z$ which satisfies the  following conditions:
\begin{itemize}
\item[ (i)]   $f\in \mathcal{C}_s(W,Z)\cap \mathcal{O}_s(W^{\text{o}},Z);$
\item[(ii)] $f$ is locally  bounded;
\item[(iii)]    $f|_{A\times B}$ is continuous,
\end{itemize}
there  exists a  unique  mapping   
$\hat{f}\in\mathcal{O}(\widehat{W^{'}},Z)$
which  admits the angular limit $f$ at all points of $W\cap W^{'}.$ 

If $A$ and $B$  are  Borel  sets  or if $X=Y=\C$ then 
 $ \widehat{W} = \widehat{W^{'}}.$
 
If  $Z=\C$ and $\vert f\vert_W<\infty,$ then
\begin{equation*}
 \vert \hat{f}(z,w)\vert\leq \vert f\vert_{A\times B}^{1-\omega(z,A^{'},D)-\omega(w,B^{'},G)} \vert
 f\vert_W^{\omega(z,A^{'},D)+\omega(w,B^{'},G)},\qquad (z,w)\in\widehat{W^{'}}.
\end{equation*}
\end{thm}

Theorem \ref{application1_2} generalizes, in some  sense, the  result of \cite{pn2}.

In the  above  theorem we have  used  the  equality  $\widehat{\widetilde{W}}=\widehat{W^{'}}$  when either $A$ and $B$  are  Borel  sets  or 
 $X=Y=\C.$
This  follows from   the  identity $\omega(\cdot,A,D)=\widetilde{\omega}(\cdot,A,D)$ when either $A$ is  a Borel set or $D\subset \C$ 
(see Theorem  4.6 in \cite{pn2}). On the other hand, we can sharpen Theorem \ref{application1_2}  further, namely,
 hypothesis (i) can be  replaced  by a weaker hypothesis (i')  as  follows:
\begin{itemize}
\item[ (i')]  {\it for  any $a\in A$ the mapping $f(a,\cdot)|_G$ is  holomorphic and has  angular limit $f(a,b)$
at all points $b\in B,$  and   for  any $b\in B$ the mapping  $f(\cdot,b)|_D$ is  holomorphic and has  angular limit $f(a,b)$
at all points $a\in A.$}      
\end{itemize}
To  see this  it suffices to observe that the hypotheses of Theorem \ref{Poletsky} and Theorem \ref{thm_local_euclidean}
can be  weakened considerably  when the bounded open set $D$    therein  is just  one-dimensional.   

\subsection{System of conical  approach  regions} 
The  remaining part  of this section is  devoted  to  two important  applications of  Theorem A:
a boundary cross theorem and  a mixed  cross theorem.
 In order  to formulate them,  we need  to introduce some   terminology and  notation.

Let  $X$  be  an arbitrary  complex manifold  and $D\subset  X$  an open subset.
 We say  that  a set   $A\subset \partial D$ is {\it locally contained in a  generating manifold} if  there  exist  an (at most countable) 
 index set $J\not=\varnothing,$ a  family  of open subsets 
 $(U_j)_{j\in J}$  of $X$ and  a  family  of 
   {\it generating manifolds}\footnote{ A differentiable  submanifold  $\mathcal{M}$ of a complex manifold
  $X$ is  said to be a {\it generating
  manifold}
  if  for all  $\zeta\in\mathcal{M},$   every  complex vector subspace of  $T_{\zeta}X$ containing   $T_{\zeta}\mathcal{M}$
coincides  with 
    $T_{\zeta}X.$}   $(\mathcal{M}_j)_{j\in J}$ such that 
    $A\cap U_j\subset \mathcal{M}_j,$  $j\in J,$ and that $A\subset  \bigcup_{j\in J}  U_j.$
The dimensions of $\mathcal{M}_j$ may  vary according  to $j\in J.$
Given   a set   $A\subset \partial D$ which is locally contained in a  generating manifold, we say that $A$ is {\it of positive size} if under the  above  notation  $\sum_{j\in J}\mes_{ \mathcal{M}_j}(A\cap U_j)>0,$
where  $\mes_{\mathcal{M}_j}$ denotes the    Lebesgue measure on  $\mathcal{M}_j.$ 
 A point $a\in A$ is  said  to be  a {\it density point} of $A$  if it is a density point of $A\cap U_j$ on $\mathcal{M}_j$ for some $j\in J.$
 Denote  by  $A^{'}$ the set of density points   of $A.$
 
 Suppose now that  $A\subset \partial D$ is  of positive size.
 We  equip $D$  with the  system of conical approach regions   supported on   $A.$
  Using the   work of  B. J\"{o}ricke  (see, for  example, Theorem 3, pages 44--45 in \cite{jr}),
 one can show that\footnote{  A complete proof will be  available in \cite{nv2}.}   $A$  is  locally  pluriregular  at all density points 
  of $A.$ 
   Observe  that    $\mes_{\mathcal{M}_j}\big((A\setminus A^{'})\cap U_j\big)=0$  for $j\in J.$
Therefore, it is  not difficult  to  show that   $A^{'}$  is  locally pluriregular. Choose an increasing sequence  
$(A_n)_{n=1}^{\infty}$  of subsets of $A$  such that 
  $A_n\cap U_j$  is  closed and  $   
\mes_{\mathcal{M}_j} \Big( (A\setminus \bigcup\limits_{n=1}^{\infty} A_n)\cap U_j\Big)=0$ for $j\in J.$
Observe that  $A^{'}_n$ is locally  pluriregular,  $\overline{A^{'}_n\cap U_j}\subset A$ for $j\in J$
 and that  $\widehat{A}:= \bigcup\limits_{n=1}^{\infty} A^{'}_n$  is locally pluriregular
and that    $\widehat{A}$ is   locally pluriregular at all points of  $A^{'}.$  Consequently, it  follows 
from  Definition \ref{defi_pluri_measure}  that 
\begin{equation*}
\widetilde{\omega}(z,A,D)\leq  \omega(z,A^{'},D),\qquad  z\in D.
\end{equation*}
This estimate, combined with  Theorem A, implies the  following result  which is  a generalization 
in higher  dimensions of   Theorem  \ref{application1_2}.

  \begin{thm} \label{application1_3}  
  Let $X,\ Y$  be  two complex  manifolds,
  let $D\subset X,$ $ G\subset Y$ be  two open sets, and let
  $A$ (resp. $B$) be  a  subset of  $\partial D$ (resp.
  $\partial G$). 
   $D$  (resp.  $G$) is  equipped with a 
  system of conical approach  regions 
  $\big(\mathcal{A}_{\alpha}(\zeta)\big)_{\zeta\in\overline{ D},\  \alpha\in I_{\zeta}}$
  (resp.  $\big(\mathcal{A}_{\beta}(\eta)\big)_{\eta\in\overline{ G},\  \beta\in I_{\eta}}$)  supported on  $A$  (resp. on $B$).
   Suppose  in  addition that
    $A$  and  $B$  are  of positive size.
   Let $Z$  be a     complex  analytic space    possessing the  Hartogs extension property.
   Define 
   \begin{eqnarray*}
   W^{'}  &:= &\X(A^{'},B^{'};D,G),\\
    \widehat{W^{'}}  &:= &\left\lbrace  (z,w)\in D\times G:\  \omega(z,A^{'},D)+\omega(w,B^{'},G)<1     \right\rbrace,
   \end{eqnarray*}
     where $A^{'}$  (resp.  $B^{'}$) is  the set of  density  points   of $A$  (resp.  $B$).

 Then,
   for  every  mapping   $f:\ W\longrightarrow Z$
  which satisfies the   following  conditions: 
   \begin{itemize}
   \item[$\bullet$]    $f\in\mathcal{C}_s(W,Z)\cap \mathcal{O}_s(W^{\text{o}},Z);$
    \item[$\bullet$] $f$ is locally bounded;
    \item[$\bullet$]          $f|_{A\times B}$is    continuous,
    \end{itemize}
     there exists a  unique mapping  
$\hat{f}\in\mathcal{O}(\widehat{W^{'}},Z)$ which
 admits  $\mathcal{A}$-limit $f(\zeta,\eta)$ at  every point
  $(\zeta,\eta)\in  W\cap   W^{'}.$     

 If,  moreover, $Z=\C$ and  $\vert f\vert_W<\infty,$ then
\begin{equation*}
 \vert \hat{f}(z,w)\vert\leq \vert f\vert_{A\times B}^{1-\omega(z,A^{'},D)-\omega(w,B^{'},G)} \vert
 f\vert_W^{\omega(z,A^{'},D)+\omega(w,B^{'},G)},\qquad (z,w)\in\widehat{W^{'}}.
\end{equation*}
\end{thm}

The second  application is  a very general mixed cross   theorem. 

 \begin{thm} \label{application2_3}  
  Let $X,\ Y$  be  two complex  manifolds,
  let $D\subset X,$ $ G\subset Y$ be  open sets,  let
  $A$ be a subset of  $\partial D,$   and let  $B$  be  a subset of
  $ G$. 
   $D$ is  equipped with the
  system of conical approach regions 
  $\big(\mathcal{A}_{\alpha}(\zeta)\big)_{\zeta\in\overline{ D},\  \alpha\in I_{\zeta}}$  supported on  $A$ 
   and  $G$  is equipped with the canonical   system of 
   approach  regions  
  $\big(\mathcal{A}_{\beta}(\eta)\big)_{\eta\in\overline{ G},\  \beta\in I_{\eta}}$.
   Suppose in addition that 
    $A$  is  of positive size.
    Let $Z$  be a     complex  analytic space    possessing the  Hartogs extension property. 
   Define 
    \begin{eqnarray*}
   W^{'}  &:= &\X(A^{'},B^{\ast};D,G),\\
    \widehat{W^{'}}  &:= &\left\lbrace  (z,w)\in D\times G:\  \omega(z,A^{'},D)+\omega(w,B^{\ast},G)<1     \right\rbrace,
   \end{eqnarray*}
    where $A^{'}$  is the set of  density points  of $A$ 
    and  $B^{\ast}$ denotes, as usual (see Subsection \ref{Subsection_approach_regions}  above), 
    the set  of points in  $\overline{B}\cap  G$  at  which $B$  is  locally pluriregular.

 Then,
   for  every mapping   $f:\ W\longrightarrow Z$
  which satisfies  the  following conditions: 
   \begin{itemize}
   \item[$\bullet$]    $f\in\mathcal{C}_s(W,Z)\cap \mathcal{O}_s(W^{\text{o}},Z);$
    \item[$\bullet$] $f$ is  locally bounded  along  $A\times G,$
    \end{itemize}
     there  exists a unique  mapping    
$\hat{f}\in\mathcal{O}(\widehat{W^{'}},Z)$ which
 admits  $\mathcal{A}$-limit $f(\zeta,\eta)$ at every point
  $(\zeta,\eta)\in  W\cap   W^{'}.$        

 If,  moreover, $Z=\C$ and  $\vert f\vert_W<\infty,$ then
\begin{equation*}
 \vert \hat{f}(z,w)\vert\leq \vert f\vert_{A\times B}^{1-\omega(z,A^{'},D)-\omega(w,B^{\ast},G)} \vert
 f\vert_W^{\omega(z,A^{'},D)+\omega(w,B^{\ast},G)},\qquad (z,w)\in\widehat{W^{'}}.
\end{equation*}
\end{thm}

\medskip

\noindent{\bf Concluding remarks.}
In ongoing  joint-works  with Pflug \cite{pn5,pn6}  we   develop  
 new cross theorems with singularities. On the other hand, in \cite{pn4} the problem  of optimality
 of the envelope of holomorphy $\widehat{\widetilde{W}}$ in Theorem A  has been investigated. 

\end{document}